\documentclass[reqno,10pt]{amsart}
\usepackage{amssymb,amsmath,amsthm,}
\usepackage{a4wide}
\usepackage{amscd}
\usepackage{amsfonts}
\usepackage{amssymb}
\usepackage{latexsym}
\usepackage{color}
\usepackage{esint}

\newtheorem{theorem}{Theorem}

\newtheorem{corollary}[theorem]{Corollary}
\newtheorem{definition}[theorem]{Definition}
\newtheorem{lemma}[theorem]{Lemma}
\newtheorem{proposition}[theorem]{Proposition}
\newtheorem{remark}[theorem]{Remark}

\let\d=\delta

\let\p=\partial

\let\O=\Omega

\newcommand{\be}{\begin{equation}}
\newcommand{\bm}{\begin{multline}}
\newcommand{\ee}{\end{equation}}
\newcommand{\dd}{\mathrm{d}}

\def\p{\partial}

\def\O{\Omega}

\def\d{\mathrm{d}}
\def\B{\begin{equation}}
\def\E{\end{equation}}
\def\BN{\begin{eqnarray*}}
\def\EN{\end{eqnarray*}}

\begin{document}


\title{Uniform Regularity for free-boundary navier-stokes equations with surface tension}
\author{tarek elgindi and donghyun lee}
\maketitle


\begin{abstract}
We study the zero-viscosity limit of free boundary Navier-Stokes equations with surface tension in $\mathbb{R}^3$ thus extending the work of Masmoudi and Rousset [1] to take surface tension into account. Due to the presence of boundary layers, we are unable to pass to the zero-viscosity limit in the usual Sobolev spaces. Indeed, as viscosity tends to zero, normal derivatives at the boundary should blow-up. To deal with this problem, we solve the free boundary problem in the so-called Sobolev co-normal spaces (after fixing the boundary via a coordinate transformation). We prove estimates which are uniform in the viscosity. And after inviscid limit process, we get the local existence of free-boundary Euler equation with surface tension. One of the main differences between this work and the work [1] is our use of time-derivative estimates and certain properties of the Dirichlet-Neumann operator.
\end{abstract}


\tableofcontents


\section{Introduction}

	The water-wave problem has been studied for several decades from several different points of view. First the local existence for free boundary problem of Navier-Stokes equation without surface tension was shown by Beale \cite{TB}.  Allain \cite{GA} and Tani \cite{AT} proved local-existence for the free-boundary Navier-Stokes equation with surface tension in the case of two dimensions and three dimensions, respectively. Moreover, with surface tension, global regularity was also studied by Beale \cite{TB}. 

In the case where the fluid is assumed to be inviscid and irrotational, the problem can be reduced to the boundary. Recently, global regularity was achieved by S. Wu \cite{SW} and by Germain, Masmoudi and Shatah \cite{PGNMJS} for small data. In the general case (where the vorticity may be non-zero) local well-posedness was proven by a number of different authors first by Christodoulou and Lindblad \cite{DCHL} and Lindblad \cite{HL}, then Coutand and Shkoller \cite{DCSS}, Lannes \cite{DL}, Shatah and Zheng \cite{JSCZ}, and Masmoudi and Rousset \cite{NMFR}.

In this paper we consider the vanishing viscosity limit for the water wave problem when surface tension is taken into account. The inviscid limit problem of the free-surface Navier-Stokes equation \emph{without} surface tension was studied by Masmoudi and Rousset in [1]. When surface tension is not taken into account, the boundary, $h,$ has same regularity as the velocity, $u$,(say $H^m$). In the process of doing high order energy estimates, one loses half a derivative due to some commutators. That commutator comes from $D^m \nabla\varphi$, where $\varphi$ is harmonic extension of $h$ to the interior of the domain, which is $\frac{1}{2}$ more regular than $h$. The main idea of the paper \cite{NMFR} is to use Alinhac's good unknown which reduces the order of loss via a critical cancellation. A second important component of \cite{NMFR}, which sets it apart from the rest of the works on the free-boundary Euler equations, is that the authors are forced to use spaces which only measure co-normal regularity. This is due to the presence of boundary layers which form during the process of sending the viscosity $\varepsilon>0$ to zero. Indeed, because of the boundary layer, we expect that near the boundary $u^\varepsilon$ behaves like $$u^\varepsilon \sim u(t,x) + \sqrt{\varepsilon}U(t,y,\frac{z}{\sqrt{\varepsilon}}),$$ where  $u$ is the solution of the free-boundary Euler equation, $U$ is a some profile, $y$ is 2-d horizontal variable, and $z$ is 1-d vertical variable. So, for high order Sobolev space, we cannot hope to get interval of time independent of $\varepsilon$, which is crucial to get stong compactness of solution sequences. Hence we consider a Sobolev co-normal space, in which we expect to maintain boundedness of the Lipschitz norm as well as boundedness of higher order co-normal derivatives on an interval of time independent of kinematic viscosity $\varepsilon > 0$.  \\ 
\indent Now let us consider the similar case with surface tension. We will still use Sobolev co-normal spaces like in \cite{NMFR} because boundary layers are still present. However, in this case the boundary is more regular so we will not need Alinhac's good unknown. Our main problem comes from the fact that the pressure term in the Euler equations becomes significantly less regular when surface tension is introduced. We thus encounter commutators with order $m+\frac{3}{2}$, which we cannot control. For this reason, we decided to do energy estimates using space \emph{and} time derivatives. This helps because time derivatives actually count for $\frac{3}{2}$ space derivatives on the boundary (this is deduced by studying the properties of the Dirichlet to Neumann mapping). Using this fact, we can derive local existence for a time interval independent to $\varepsilon$. Finally, we deduce the solution of free-boundary Euler equation (subject to surface tension) as $\varepsilon\rightarrow 0$, using a strong compactness argument. Also note that upper index $(\cdot)^b$ means function value on free-surface, $(\cdot)\vert_{\text{free-surface}}$. 

\subsection{The Free-boundary Navier-Stokes equations}
We solve the incompressible free-boundary Navier-Stokes equations under the effect of gravity in an unbounded domain. Assume that above the free-surface of the fluid is a vacuum.  The system we get is:

\begin{equation} \label{1.1}
\partial_t u + u\cdot\nabla u + \nabla p = \varepsilon\Delta u,\,\,\,\,\,\,(x,y,z)\in \Omega_t,\quad t>0, \quad \varepsilon>0,
\end{equation}
\begin{equation} \label{1.2}
\nabla\cdot u = 0,\,\,\,\,\,\,\,(x,y,z)\in \Omega_t,
\end{equation}
where $\Omega_t$ is free domain, occupied by fluid, and $\varepsilon>0$ is kinematic viscosity. We write the fluid boundary as $h$, a function of horizontal variables $x$ and $y$, so that
$$
\Omega_t := \left\{(x,y,z)\in \mathbb{R}^3,\,\,\,\, z < h(t,x,y)\right\},
$$

Our first boundary condition is the moving boundary condition (or called kinematic boundary condition), which roughly says that the boundary moves with the fluid:
\begin{equation} \label{1.3}
\p_t h = u (t,x, y, h(t,x,y) )\cdot \mathbf{N},\,\,\,\,(x_1,x_2)\in \mathbb{R}^2,
\end{equation}
where $\mathbf{N} = (-\p_1 h, \p_2 h, 1) := (-\nabla_{y} h, 1)$ is outward normal direction vector on the boundary.
Our second boundary condition is the continuity of stress tensor on the boundary.\
\begin{equation} \label{1.4}
p^b\mathbf{n} - 2\varepsilon (Su)^b \mathbf{n} = gh\mathbf{n} - \nabla\cdot\frac{\nabla h}{\sqrt{1+|\nabla h|^2}}\mathbf{n} ,\,\,\,\,\,\,x \in \partial\Omega_t,
\end{equation}
where $\mathbf{n} = \mathbf{N}/|\mathbf{N}|$, outward unit normal vector, and $Su$ is the symmetric part of the gradient of $u$
$$
Su := \frac{(\nabla u) + (\nabla u)^T}{2}.
$$
Note that $g$ is gravitational constant. In this paper, we consider system (\ref{1.1})-(\ref{1.4}) and its vanishing viscosity limit $\varepsilon\rightarrow 0$. 

\subsection{Parametrization to a Fixed domain}
The first step is to rewrite the problem on the fixed domain $S := \{(x,y,z)| \  z < 0 \}$. This can be done by a diffeomorphism $\Phi(t,\cdot)$,
\begin{equation} \label{1.5}
\begin{split}
\Phi(t,\cdot): S=\mathbb{R}^2 &\times (-\infty,0) \rightarrow \Omega_t, \\
(x,y,z)\mapsto (x,y,\varphi(t,x,y,z)) &:= (x,y,Az+\eta(t,x,y,z)).
\end{split}
\end{equation}
We use function $v$ and $q$ to denote the velocity and pressure on the fixed domain $S$.
\begin{equation} \label{1.6}
v(t,x,y,z)=u(t,\Phi(t,x,y,z)),\,\,q(t,x,y,z)=p(t,\Phi(t,x,y,z)).
\end{equation}

There are several different choices for $\Phi$ and we have to decide which one is optimal for our purposes. 
We will need to find $\varphi(t,\cdot)$ so that $\Phi(t,\cdot)$ is a diffeomorphism (Surely, $\partial_z\varphi \geq 0$, because it is diffeomorphism). 
One easy option is to set $\varphi(t,x,y,z)=z + h(t,x,y)$. But it is more useful to take a function $\Phi$ which is actually \emph{more} regular than $h$. If one thinks about using a harmonic extension, we see that it is possible for $\Phi$ to be $\frac{1}{2}$ of a derivative more regular than $h$.  
We take a smoothing diffeomorphism as was done in [1]:
\begin{equation*} \label{eta def}
\varphi(t,x,y,z) := Az+\eta(t,x,y,z).
\end{equation*} 
To ensure that $\Phi(0,\cdot)$ is a diffeomorphism, $A$ should be chosen so that 
\begin{equation*}
\partial_z\varphi(0,x,y,z)\geq 1,\,\,\,\,\forall (x,y,z)\in S,
\end{equation*}
and $\eta$ is given by the extension of $h$ to the domain $S$, defined by 
\begin{equation*} \label{extension}
\hat{\eta}(\xi,z)=\chi(z\xi)\hat{h}(\xi). \\
\end{equation*}

We want to rewrite equations (\ref{1.1})-(\ref{1.5}) on the moving domain $\Omega_t$ as equations on the fixed domain $S$. Using change of variables, we get
\begin{equation*} 
\begin{split}
(\partial_i u)(t,x,y,\varphi) &= (\partial_i v - \frac{\partial_i\varphi}{\partial_z\varphi}\partial_z v)(t,x,y,z),\quad i=t,1,2,  \\
(\partial_3 u)(t,x,y,\varphi) &= (\frac{1}{\partial_z\varphi}\partial_z v)(t,x,y,z).
\end{split}
\end{equation*}
So we define the following operators on fixed domain $S$ :

\begin{definition} \label{varphipartial}
We define the following differential operators. 
\begin{equation*}
\begin{split}
\partial_i^\varphi := \partial_i - \frac{\partial_i\varphi}{\partial_z\varphi}\partial_z,\,\,\,\,\,\, i=t,1,2, \,\,\,\,\,\,\,\,\,\,\ \partial_3^\varphi := \frac{1}{\partial_z\varphi}\partial_z,   \\
\nabla^\varphi := (\p^\varphi_{1}, \p^\varphi_{2}, \p^\varphi_{3}),\quad \Delta^\varphi := \p^\varphi_{11} + \p^\varphi_{22} + \p^\varphi_{33} . 
\end{split}
\end{equation*}
With these definition, $\partial_i u \circ \Phi = \partial_i^\varphi v, \,\,\,\,\,\,\,i=t,1,2,3$. And each $\p_i^{\varphi}$'s commute each other. \\
\end{definition}

\noindent Hence our equations in $S$ are,
\begin{equation} \label{1.12}
\partial_t^\varphi v + v\cdot\nabla^\varphi v + \nabla^\varphi q = \varepsilon\Delta^\varphi v,\quad \text{in} \quad(x,y,z)\in S,\quad \varepsilon>0,
\end{equation}
\begin{equation} \label{1.13}
\nabla^\varphi\cdot v = 0,\quad \text{in} \quad (x,y,z)\in S,
\end{equation}
\begin{equation} \label{1.14}
\partial_t h = v^b\cdot \mathbf{N} ,\quad \text{on} \quad \{z=0\},  \\
\end{equation}
\begin{equation} \label{1.15}
q^b\mathbf{n} - 2\varepsilon (S^\varphi v)^b \mathbf{n} = gh\mathbf{n} - \nabla\cdot\frac{\nabla h}{\sqrt{1+|\nabla h|^2}}\mathbf{n},\quad\text{on} \quad \{z=0\},  \\
\end{equation}

\noindent In this paper, we use extended definition for normal vector $\mathbf{N}$. 
\begin{definition} \label{extendnormal}
	We define the folloiwng extended normal vector in $S$ and on $\p S$.
	\[
		\mathbf{N} := (-\p_1 \eta, -\p_2 \eta, 1) = (-\nabla_{x,y}\eta,1),\quad \mathbf{n} = \frac{\mathbf{N}}{|\mathbf{N}|}.
	\]
	On the boundary $\{z=0\}$, $\mathbf{N}$ and $\mathbf{n}$ are just standard outward vectors on the boundary as we defined in (\ref{1.3}).
\end{definition}

\subsection{Functional Framework and Notations}
We introduce co-normal spaces and some function spaces that are tailored to our problem. From boundary layer behavior and boundary condition (\ref{1.4}), we cannot perform normal derivative $\p_3$ more than one on the boundary. So we redefine normal derivative so that it is equivalent to standard normal derivatives $\p_3$, away from the boundary and vanishes on the boundary. We multiply function $g(z)$ to $\p_3$ so that
\[
	g(0)=0\quad\text{and}\quad \Big| \frac{\dd^{k}}{\dd z^{k}} g(z) \Big| \leq \text{uniform bound}\,\, C_{k},
\]
for each $k\in \{0\}\cup\mathbb{N}$ and $z\leq 0$. One possibility is to pick $g(z) = \frac{z}{1-z}$.

\begin{definition} \label{space def}
We define Sobolev co-normal derivatives on $S$ as:
\begin{equation} \label{conormal derivatives}
Z_1 = \partial_1,\,\,\,\,\,Z_2 = \partial_2,\,\,\,\,\,Z_3 = \frac{z}{1-z}\partial_z,\,\,\,\,\,Z^\alpha := Z^{(\alpha_1,\alpha_2,\alpha_3)} := Z_{1}^{\alpha_1} Z_{2}^{\alpha_2} Z_{3}^{\alpha_3}.
\end{equation}
We also use the following symbol:
\begin{equation} \label{m conormal}
Z^{m} := \p_t^k Z_{1}^{\alpha_1} Z_{2}^{\alpha_2} Z_{3}^{\alpha_3},\quad \text{for some}\quad k,\alpha_1,\alpha_2,\alpha_3 \in\{0\}\cup\mathbb{N}, \quad k + |\alpha| := k + \alpha_1 + \alpha_2 + \alpha_3 = m .
\end{equation}
There are many cases of $(\alpha_t,\alpha_1,\alpha_2,\alpha_3)$ for fixed $m$ but we will sum all those cases later, so we do not have to distinguish each case. Using co-normal derivatives, we define Sobolev co-normal spaces, similar as standard Sobolev spaces. Here we use $\p S$ to denote upper boundary $\{z=0\}$.
\begin{eqnarray*}
	H^{m}_{co}(S) &:=& \{ f\in L^2(S) \,\vert\, Z_{1}^{\alpha_1} Z_{2}^{\alpha_2} Z_{3}^{\alpha_3} f \in L^2(S),\quad \alpha_1 + \alpha_2 + \alpha_3 \leq m  \},  \\
	W^{m,\infty}_{co}(S) &:=& \{ f\in L^{\infty}(S) \,\vert\, Z_{1}^{\alpha_1} Z_{2}^{\alpha_2} Z_{3}^{\alpha_3} f \in L^{\infty}(S),\quad \alpha_1 + \alpha_2 + \alpha_3 \leq m  \},  \\
	H^{m}_{tan}(S) &:=& \{ f\in L^2(S) \,\vert\, Z_{1}^{\alpha_1} Z_{2}^{\alpha_2} f \in L^2(S),\quad \alpha_1 + \alpha_2 \leq m  \},  \\
	W^{m,\infty}_{tan}(S) &:=& \{ f\in L^{\infty}(S) \,\vert\, Z_{1}^{\alpha_1} Z_{2}^{\alpha_2} f \in L^{\infty}(S),\quad \alpha_1 + \alpha_2 \leq m  \},  \\
	H^{m}_{co}(\p S) &:=& \{ f\in L^2(\p S) \,\vert\, Z_{1}^{\alpha_1} Z_{2}^{\alpha_2} f \in L^2(\p S),\quad \alpha_1 + \alpha_2 \leq m  \},  \\
	W^{m,\infty}_{co}(\p S) &:=& \{ f\in L^{\infty}(\p S) \,\vert\, Z_{1}^{\alpha_1} Z_{2}^{\alpha_2} f \in L^{\infty}(\p S),\quad \alpha_1 + \alpha_2 \leq m  \},
\end{eqnarray*}
with following norms in each function spaces.
\begin{eqnarray*}
\left\|f\right\|_{H_{co}^m(S)}^2 &:=& \left\|f\right\|_{m}^2 := \sum_{|\alpha |\leq m} \left\|Z^\alpha f\right\|_{L^2(S)}^2,  \\
\left\|f\right\|_{W_{co}^{m,\infty}(S)} &:=& \left\|f\right\|_{m,\infty} := \sum_{|\alpha|\leq m}\left\|Z^\alpha f\right\|_{L^\infty(S)},   \\
\left\|f\right\|_{H_{tan}^m(S)}^2 &:=& \left\|f\right\|_{m,tan}^2 := \sum_{|\alpha |\leq m, \alpha_3=0} \left\|Z^\alpha f\right\|_{L^2(S)}^2,  \\
\left\|f\right\|_{W_{tan}^{m,\infty}(S)} &:=& \left\|f\right\|_{m,\infty,tan} := \sum_{|\alpha|\leq m,\alpha_3=0}\left\|Z^\alpha f\right\|_{L^\infty(S)},   \\
\left|f\right|_{H^m(\p S)}^2 &:=& \left|f\right|_{m}^2 := \sum_{|\alpha |\leq m,\alpha_3=0} \left|Z^\alpha f\right|_{L^2(\p S)}^2,  \\
\left|f\right|_{W^{m,\infty}(\p S)} &:=& \left|f\right|_{m,\infty} := \sum_{|\alpha |\leq m,\alpha_3=0}\left|Z^\alpha f\right|_{L^\infty(\p S)}.
\end{eqnarray*} 
Note that for boundary functions, co-normal space and standard Sobolev spaces are identical, since there is no $Z_3$. When $s=0$, we write $\|\cdot\| := \|\cdot\|_{L^2} $. Also note that $H^m(S) \hookrightarrow H^{m}_{co}(S) \hookrightarrow H^{m}_{tan}(S)$ by property of $\frac{z}{1-z}$ and definition of $H^{m}_{tan}$.
\end{definition}


Because we will consider time derivatives in energy estimate, it would be convenient to define function spaces containing time derivatives.
\begin{definition} \label{XY space}
Let $v$ and $h$ are smooth functions defined in $S$ and on $\p S$ respectively. We define function spaces $X^{m,s}$, $Y^{m,s}$, $\mathcal{H}^{m,s}$, and $\mathcal{K}^{m,s}$.
\begin{eqnarray*}
	X^{m,s}(\p S) &:=& \{ f\in L^2(\p S) \,\vert\, \p_t^{k} Z_{1}^{\alpha_1} Z_{2}^{\alpha_2} f \in H^{s}(\p S),\quad k + \alpha_1 + \alpha_2 \leq m  \},  \\
	X^{m,s}(S) &:=& \{ f\in L^2(S) \,\vert\, \p_t^{k} Z_{1}^{\alpha_1} Z_{2}^{\alpha_2} Z^{\alpha_3}_{3}f \in H^{s}_{co}(S),\quad k + \alpha_1 + \alpha_2 + \alpha_3 \leq m  \},  \\
	Y^{m,s}(\p S) &:=& \{ f\in L^{\infty}(\p S) \,\vert\, \p_t^{k} Z_{1}^{\alpha_1} Z_{2}^{\alpha_2} f \in W^{s,\infty}(\p S),\quad k + \alpha_1 + \alpha_2 \leq m  \},  \\
	Y^{m,s}(S) &:=& \{ f\in L^{\infty}(S) \,\vert\, \p_t^{k} Z_{1}^{\alpha_1} Z_{2}^{\alpha_2} Z^{\alpha_3}_{3}f \in W^{s,\infty}_{co}(S),\quad k + \alpha_1 + \alpha_2 + \alpha_3 \leq m  \},  \\
	\mathcal{H}^{m,s}(S) &:=& \{ f\in L^{2}(S) \,\vert\, \p_t^{k} Z_{1}^{\alpha_1} Z_{2}^{\alpha_2} \p^{\alpha_3}_{3}f \in H^{s}(S),\quad k + \alpha_1 + \alpha_2 + \alpha_3 \leq m  \},  \\
	\mathcal{K}^{m,s}(S) &:=& \{ f\in L^{\infty}(S) \,\vert\, \p_t^{k} Z_{1}^{\alpha_1} Z_{2}^{\alpha_2} \p^{\alpha_3}_{3}f \in W^{s,\infty}(S),\quad k + \alpha_1 + \alpha_2 + \alpha_3 \leq m  \}.  \\
\end{eqnarray*}
Especially, to avoid confusion we use $|\cdot|$ for functions defined on the boundary $\p S$ and $\|\cdot\|$ for functions defined on the interior $S$. Using this norm symbols, we define norms of above function spaces $X^{m,s}$, $Y^{m,s}$, $\mathcal{H}^{m,s}$, and $\mathcal{K}^{m,s}$.
\begin{eqnarray*}
\left|h\right|_{X^{m,s}}^2 &:=& \sum_{(k,\alpha),k+|\alpha |\leq m}
 |\partial_t^k D_{h}^\alpha h |_{s}^2,\quad\text{where $h$ is defined on}\quad \p S,  \\
\left\|v\right\|_{X^{m,s}}^2 &:=& \sum_{(k,\alpha),k+|\alpha |\leq m}
 \|\partial_t^k Z^\alpha v \|_{s}^2,\quad\text{where $v$ is defined on}\quad S,   \\
\left|h\right|_{Y^{m,s}} &:=& \sum_{(k,\alpha),k+|\alpha |\leq m}
 |\partial_t^k D_{h}^\alpha h |_{s,\infty},\quad\text{where $h$ is defined on}\quad \p S,  \\
\left\|v\right\|_{Y^{m,s}} &:=& \sum_{(k,\alpha),k+|\alpha |\leq m}
 \|\partial_t^k Z^\alpha v \|_{s,\infty},\quad\text{where $v$ is defined on}\quad S,  \\
\left\|v\right\|_{\mathcal{H}^{m,s}}^2 &:=& \sum_{(k,\alpha),k+|\alpha |\leq m}
\|\partial_t^k \p_1^{\alpha_1} \p_2^{\alpha_2} \p_3^{\alpha_3} v \|_{H^{s}}^2,\quad\text{where $v$ is defined in}\quad S,  \\
\left\|v\right\|_{\mathcal{K}^{m,s}}^2 &:=& \sum_{(k,\alpha),k+|\alpha |\leq m}
\|\partial_t^k \p_1^{\alpha_1} \p_2^{\alpha_2} \p_3^{\alpha_3} v \|_{W^{s,\infty}}^2,\quad\text{where $v$ is defined in}\quad S,
\end{eqnarray*}
where $D_{h}^\alpha$ means horizontal derivatives $\p_1^{\alpha_1}\p_2^{\alpha_2}$. This definition means $m+s$ order (including both time and space) Sobolev co-normal spaces which contains at least $s$ number of spatial co-normal derivatives. $\mathcal{H}^{m,s}$ is defined similar as $X^{m,s}$, except every derivatives are standard derivatives, not co-normal. Note that 
\[
	\mathcal{H}^{m,s}\hookrightarrow X^{m,s},\quad \mathcal{K}^{m,s}\hookrightarrow Y^{m,s},\quad  X^{m+s,0} \hookrightarrow X^{m,s}\quad\text{and}\quad Y^{m+s,0} \hookrightarrow Y^{m,s},
\]
by definition.
\end{definition}

\noindent In co-normal spaces, there is no proper notion of half derivatives, so we define the following for convenience.
\begin{definition} \label{bracket}
We use the following definition for fractional order in Sobolev co-normal derivatives.
\begin{equation*}
[\frac{m}{2}] := \begin{cases}
\frac{m}{2} \in \{0\}\cup\mathbb{N},\,\,\,\,m\,\,\,\,\text{even}, \\
\frac{m-1}{2} \in \{0\}\cup\mathbb{N},\,\,\,\,m\,\,\,\,\text{odd} .\\
\end{cases}
\end{equation*}
\end{definition}

\subsection{Main Result}
We state main results of this paper. The following theorem states the local energy estimate uniform in kinematic viscosity $\varepsilon>0$.

\begin{theorem} \label{theorem-uniform}
For fixed sufficiently large $m\geq 6$, let initial data be given so that
\begin{equation}
	I_{m}(0) := \|v_{0}\|_{H^{m}(S)} + \sum_{k=1}^{m} \varepsilon^{k}\|v_{0}\|_{H^{m+k}(S)} + |h_{0}|_{H^{\frac{3}{2}m + 1}(\p S)} \leq R
\end{equation}
and satisfy compatibility conditions
\begin{equation} \label{compatible}
\Pi \big( S^\varphi v^\varepsilon \mathbf{n} \big) \big\vert_{t=0,z=0}= 0, \quad \Pi := {I} - \mathbf{n}^{b}\otimes\mathbf{n}^{b}.
\end{equation}

Then for $\forall \varepsilon \in (0,1]$, there exist time independent $T>0$ and some $C>0$, such that there exist a unique solution $(v^\varepsilon,h^\varepsilon)$ on $[0,T]$, and the following energy estimate hold.
\begin{equation} \label{main est}
\begin{split}
\mathcal{N}_m (T) &:= \sup_{t\in[0,T]}\Big( \|v(t)\|_{X^{m,0}}^2 + |h(t)|_{X^{m,1}}^2 + \|\partial_z v(t)\|_{X^{m-2,0}}^2 + \|\partial_z v(t)\|_{Y^{[\frac{m}{2}],0}}^2 \Big)  \\
&\quad + \|\partial_z v\|_{L_{T}^4 X^{m-1,0}}^2 + \varepsilon\int_0^T \Big( \|\nabla v (t)\|^2_{X^{m,0}} + \|\nabla\partial_z v(t)\|^2_{X^{m-2,0}} \Big) dt < C.
\end{split}
\end{equation}
\end{theorem}

\noindent Using the result of Theorem \ref{theorem-uniform}, we send $\varepsilon>0$ to zero to get a unique solution of free-boundary Euler equations with surface tension, i.e. vanishing viscosity limit.
\begin{theorem} \label{inviscid limit}
Let us assume that
\begin{equation}
\lim_{\varepsilon\rightarrow 0}\Big( \left\|v_0^\varepsilon - v_0\right\|_{L^2(S)} + \left\|h_0^\varepsilon - h_0\right\|_{H^1(\partial S)} \Big) = 0,
\end{equation}
where $(v_0^\varepsilon, h_0^\varepsilon)$ and $(v_0, h_0)$ satisfy the assumptions of Theorem \ref{theorem-uniform}. Then there exist $(v,h)$ satisfying 
\begin{equation}
v \in L^\infty([0,T],H_{co}^m(S)),\,\,\partial_z v \in L^\infty([0,T],H_{co}^{m-2}(S)),\,\,h \in L^\infty([0,T],H_{co}^{m+1}(\mathbb{R}^2))
\end{equation}
and
\begin{equation}
\lim_{\varepsilon\rightarrow 0}\sup_{[0,T]} \Big( \left\|v^\varepsilon - v\right\|_{L^2(S)} + \left\|v^\varepsilon - v\right\|_{L^\infty(S)} + \left\|h^\varepsilon - h\right\|_{H^1(\mathbb{R}^2)} + \left\|h^\varepsilon - h\right\|_{W^{1,\infty}(\mathbb{R}^2)} \Big) = 0.
\end{equation}
Moreover, $(v,h)$ is the unique solution of free-boundary Euler equation,
\begin{equation}
\begin{split}
	&\partial_t^\varphi v + (v\cdot\nabla^\varphi)v + \nabla^\varphi q = 0, \quad\text{in}\quad S  \\	 
	&\nabla^\varphi\cdot v = 0, \quad\text{in}\quad S  \\
	&\partial_t h = v^b\cdot N, \quad\text{on}\quad \p S \\
	&q = gh - \eta\nabla\cdot\Big(\frac{\nabla h}{\sqrt{1+|\nabla h|^2}} \Big), \quad\text{on}\quad \p S.	
\end{split}
\end{equation}
\end{theorem}

\subsection{Scheme of the Proof}
We briefly explain main idea of this paper in several steps.
\begin{remark}In this paper $\Lambda(\cdot,\cdot)$ denotes an increasing continuous function in all its arguments and $\Lambda_0 = \Lambda(\frac{1}{c_0})$. Both may vary from line to line. Also $C$ means a constant independent to $\varepsilon$ and also vary from line to line.
\end{remark}

\subsubsection{Energy estimate of v and h}
Let us apply $Z^m := Z_1^{\alpha_1} Z_2^{\alpha_2} Z_3^{\alpha_3}$ to (\ref{1.12}). Then we get energy estimate as following. 
\begin{equation} \label{E0}
\begin{split}
E_0(t) &:= \|v\|^2_{H_{co}^m} + |h|^2_{H^{m+1}} + \varepsilon\int_0^t\|\nabla v\|^2_{H_{co}^m} \\
&\leq C_0 + \Lambda(R)\int_0^t \Big( E_0(s) + \left\|\nabla v\right\|^2_{H_{co}^{m-1}} + |h|^2_{H^{m+\frac{3}{2}}} \Big) ds , 
\end{split}
\end{equation}

where $C_0$ is some terms depending on the initial data, $R$ contains $E_0$ and $L^\infty$-type terms with order of $[\frac{m}{2}]$. The problem is that $|h|_{H^{m+\frac{3}{2}}} $, on the right hand side, cannot be controlled by $E_0$. This term comes from the pressure estimates involving the surface tension term. To estimate $|h|_{H^{m+\frac{3}{2}}} $, we use Dirichlet-Neumann estimates, time-derivatives, and a special decomposition of the pressure term. We decompose the pressure $q$ into three pieces, $q = q^{E} + q^{NS} + q^{S}$. One piece, $q^S$ (which is the pressure due to the surface tension term) sovles, 
\begin{equation} \label{press eq}
\Delta q^S = 0,\quad q^S|_{z=0} = -\nabla\cdot\frac{\nabla h}{\sqrt{1+|\nabla h|^2}}. 
\end{equation} 
Applying $\p_t$ on kinematic boundary condition (\ref{1.14}), we get
$$
h_{tt} = v_t^b \cdot \mathbf{N} + v^b\cdot \mathbf{N}_t .
$$
And using Naiver-Stokes (\ref{1.12}), 
\begin{equation*}
h_{tt} = -(\nabla q^S)^b \cdot N + v^b\cdot N_t + (other\,\,terms).
\end{equation*}
Since $(q^S)^b \sim \Delta h$ by surface tension, we can get $h_{tt} \sim \nabla\Delta h$, so \emph{heuristically}, $\partial_t h \sim \partial_x^{\frac{3}{2}} h$. Showing $|\p_{t} h| \lesssim |\p_{x}^{3/2}h|$ is easy by direct computation. However, what we need is to claim $|\p_{x}^{3/2} h| \lesssim |\p_{t}h|$. From Proposition \ref{dirineumann}, we have a kind of Garding's inequality which gives inequality with inverse direction. We will see that (roughly),
\begin{equation} \label{simple DN}
|Z^{m-1}h|^2_{L^2 H^{\frac{5}{2}}} \lesssim |Z^{m-1}\p_t h|^2_{L^2 H^{1}} + {\theta} |Z^{m-1}\p_t h(t)|^2_{L^{2}} + \sqrt{T} |Z^{m-1}\p_t \nabla h(t)|^2_{L^{2}}.
\end{equation} 
Especially we note that above scheme $\p_{t} h \sim \p_{x}^{3/2}h$ is valid only in the sense of $L^{2}$ in time, not pointwise in time. Therefore, $|h|^{2}_{m+\frac{3}{2}}$ in the RHS of (\ref{E0}) is controlled by above (\ref{simple DN}). To control right hand sides of (\ref{simple DN}), we should perform energy estimate of which energy is gained by applying $Z^m := \p_t Z^{\alpha}, \ (|\alpha|=m-1)$, $m-1$ spatial conormal derivatives and one time derivative. Then we get 
\begin{equation*}
\begin{split}
E_1 (t) &:= \|\partial_t v\|^2_{H_{co}^{m-1}} + |\partial_t h|^2_{H^{m}} + \varepsilon\int_0^t\|\nabla\partial_t v\|^2_{H_{co}^{m-1}}  \\
&\leq C_0 + \Lambda(R)\int_0^t \Big( E_1(s) + \left\|\nabla\partial_t v\right\|^2_{H_{co}^{m-2}} + |\partial_t h|^2_{H^{m+\frac{1}{2}}} \Big) ds ,
\end{split}
\end{equation*}
where in this case, $R$ contains $E_0$ and $E_1$ and $L^\infty$-type terms up to $[\frac{m}{2}]$ order. For sufficiently small $\theta$ and time $T$, right hand side terms in (\ref{simple DN}) are absorbed or controlled by $L^{\infty}$ type energy terms in $E_{1}(t)$. However, similar situation happens on the right hand side. To cotrol $|\p_{t}h|_{H^{m+\frac{1}{2}}}$, we should use property of Dirichlet-Neumann operator (Proposition \ref{dirineumann}) again. Therefore, we need the energy of $E_2(t)$, by one more time derivative and one less space derivatives. For each steps, we should use dirichelt-Neumann operator estimate Proposition \ref{dirineumann}. This process will be repeated until the last step, $Z^m := \p_t^{m-1}Z$. 
\begin{equation*}
\begin{split}
E_2 (t) &:= \|\partial_t^2 v\|^2_{H_{co}^{m-2}} + |\partial_t^2 h|^2_{H^{m-1}} + \varepsilon\int_0^t\|\nabla\partial_t^2 v\|^2_{H_{co}^{m-2}}   \\
&\leq C_0 + \Lambda(R)\int_0^t \Big( E_2(s) + \left\|\nabla\partial_t^2  v\right\|^2_{H_{co}^{m-3}} + |\partial_t^2  h|^2_{H^{m-\frac{1}{2}}} \Big) ds  ,
\end{split}
\end{equation*}
\,\,\,\,\,\,\,\,\,\,\,\,\,\,\,\,\,\,\,\,\,\,\,\,\,\,\,\,\,\,\,\,\,\,\,\,\,\,\,\,\,\,\,\,\,\,\,\,\,\,\,\,\,\,\,\,\,\,\,\,\,\,\,\,\,\,\,\,\,\,\,\,\,\,\,\,\,\,\,\,\,\,\,\,\,\,\,\,\,\,\,\,\,\,\,\,\,\,\,\,\,\,\,\,\,\,\,\,\,\,\,\,\,\,\,\,\,\,\,\,\,\,\,\,\vdots
\\
\\
\begin{equation*}
\begin{split}
E_k(t) &:= \|\partial_t^k v\|^2_{H_{co}^{m-k}} + |\partial_t^k h|^2_{H^{m-k+1}} + \varepsilon\int_0^t\|\nabla\partial_t^k v\|^2_{H_{co}^{m-k}}   \\
&\leq C_0 + \Lambda(R)\int_0^t \Big( E_k(s) + \|\nabla\partial_t^k  v\|^2_{H_{co}^{m-k-1}} + |\partial_t^k  h|^2_{H^{m-k+\frac{3}{2}}} \Big) ds,\,\,\,\,1 \leq k \leq  m-1  ,
\end{split}
\end{equation*}
\,\,\,\,\,\,\,\,\,\,\,\,\,\,\,\,\,\,\,\,\,\,\,\,\,\,\,\,\,\,\,\,\,\,\,\,\,\,\,\,\,\,\,\,\,\,\,\,\,\,\,\,\,\,\,\,\,\,\,\,\,\,\,\,\,\,\,\,\,\,\,\,\,\,\,\,\,\,\,\,\,\,\,\,\,\,\,\,\,\,\,\,\,\,\,\,\,\,\,\,\,\,\,\,\,\,\,\,\,\,\,\,\,\,\,\,\,\,\,\,\,\,\,\,\vdots
\\
\begin{equation*}
\begin{split}
E_{m-1}(t) &:= \|\partial_t^{m-1} v\|^2_{H_{co}^{1}} + |\partial_t^{m-1} h|^2_{H^{2}} + \varepsilon\int_0^t\|\nabla\partial_t^{m-1} v\|^2_{H_{co}^{1}}   \\
&\leq C_0 + \Lambda(R)\int_0^t \Big( E_{m-1}(s) + \left\|\nabla\partial_t^{m-1} v\right\|^2_{H_{co}^{1}} + |\partial_t^{m-1} h|^2_{H^{\frac{5}{2}}} \Big) ds .
\end{split}
\end{equation*}
If we sum the above $m$ estimates, then $E_0 + E_1 + \cdots E_{m-2} + E_{m-1}$ controls every high order term of $h$, except $|\partial_t^{m-1} h|^2_{L^2 H^{\frac{5}{2}}}$ of the $m-1$ step estimate.

\subsubsection{Energy estimate for all-time derivatives}
To close energy estimate we consider $Z^m=\p_t^m$ case. Then we get the following energy,
$$
E^{m}(t) :=  \|\partial_t^{m} v\|^2_{L^2} + |\partial_t^{m} h|^2_{H^{1}} + \varepsilon\int_0^t\|\nabla\partial_t^{m} v\|^2_{L^2} .
$$
In this step, we should claim that the commutator $|\partial_t^m h|_{H^{3/2}}$ \emph{does not} appear. (We cannot apply $\p_{t}^{m+1}$ to our equation!) Since bad commutator $|\partial_t^m h|_{H^{3/2}}$ comes from pressure estimate, let us investigate pressure term in energy estimate. \\

When we apply $Z^{m}=\p_{t}^{m}$, pressure term $q^{S}$ defined in (\ref{press eq}) in the energy estimate looks like
\begin{equation} \label{cancel scheme}
\int_{0}^{t}\int_S \p_t^m v \cdot \nabla^\varphi\partial_t^m q^S  \sim \int_{0}^{t}\int_{\p S} \p_{t}^{m}u\cdot\mathbf{N} \p_{t}^{m}q  - \int_{0}^{t}\int_S (\nabla^\varphi \cdot \partial_t^m v) \partial_t^m q^S ,\quad\text{ by integration by parts}.
\end{equation}
We use kinematic boundary condition (\ref{1.14}) and divergence free condition for $v$ (\ref{1.13}) to see 
\begin{equation*}
\begin{split}
&\nabla^\varphi \cdot \partial_t^m v \sim m\partial_t\mathbf{N} \partial_t^{m-1} \p_{z}v + \text{low order terms} , \\
&\p_{t}^{m}v\cdot\mathbf{N} \sim m \p_{t} \mathbf{N} \p_{t}^{m-1} v + \text{low order terms} . \\
\end{split}
\end{equation*}
So (\ref{cancel scheme}) becomes,
\begin{equation} \label{cancel scheme 2}
\begin{split}
	\sim \int_{0}^{t}\int_{\p S} m \p_{t} \mathbf{N} \p_{t}^{m-1} v \p_{t}^{m}q  - \int_{0}^{t}\int_S m\partial_t\mathbf{N} \partial_t^{m-1} \p_{z}v \partial_t^m q^S ,\quad\text{ by integration by parts}.
\end{split}
\end{equation}
Meanwhile, from the equation (\ref{press eq}), we cannot derive an estimate for $\|\p_{t}^{m} q\|_{L^{2}}$ in the second term of the RHS. Instead we can estimate $\|\p_{t}^{m-1} \p_{z}q\|_{H^{m-1}}$. We perform integration by parts for $\p_{z}$ and $\p_{t}$ to exchange $\p_{t}^{m}q$ into $\p_{t}^{m-1}\p_{z}q$.  \\
\begin{equation} \label{cancel scheme 3}
\begin{split}
(\ref{cancel scheme 2}) \sim \underbrace{ \int_{0}^{t}\int_{\p S} m \p_{t} \mathbf{N} \p_{t}^{m-1} v \p_{t}^{m}q  - \int_{0}^{t}\int_{\p S} m\partial_t\mathbf{N} \partial_t^{m-1} v \partial_t^m q }_{\text{cancellation}} + \underbrace{ \int_{0}^{t}\int_{S} m\partial_t\mathbf{N} \partial_t^{m} v \partial_t^{m-1}\p_{z} q }_{(*)}.
\end{split}
\end{equation}
We have cancellation for the boundary integrals those contain $\p_{t}^{m}q$. Note that pressure part in $(*)$ has been changed into the same form as previous step, i.e. $Z^{m}=\p_{t}Z^{\alpha}, \ (|\alpha|=m-1)$ case. So applying Dirichlet-Neumann estimate resolve $|\p_{t}^{m-1}h|_{\frac{5}{2}}$ issue.  \\

Meanwhile, when we perform integration by parts in time, we should treat some terms evaluated at time $0$ and $t$, without time integration. Since integrand of those term has one lesser time derivatives, we should treat (roughly)
\begin{equation} \label{time bdry}
	\int_{S} \p_{t}^{m-1}v \p_{t}^{m-1}\p_{z} q \bigg\vert_{t}  \ \sim  \ \|\p_{t}^{m-1}v(t)\|^{2} + |\p_{t}^{m-1}h(t)|^{2}_{\frac{5}{2}} + (\text{other terms}) , 
\end{equation}
without time integration. We want to apply Dirichlet-Neumann operator estimate to control right hand side, but (\ref{simple DN}) requires $L^{2}$ in time to change $\p_{x}^{3/2}$ into $\p_{t}$. Instead of above simple esimate, we use the fact that $\|\p_{t}^{m-1}v(t)\|^{2}$ is not optimal. Instead of above (\ref{time bdry}), we split as
\[
\int_{S} \p_{t}^{m-1}v \p_{t}^{m-1}\p_{z} q \bigg\vert_{t}  \ \sim  \ \|\p_{t}^{m-1}|\nabla| v(t)\|^{2} + \||\nabla|^{-1}\p_{t}^{m-1}\p_{z}q(t)\|^{2} + (\text{other terms}) , 
\] 
where $|\nabla|$ is Fourier muliplier. Then, the first term on the right hand side becomes $\|\p_{t}^{m-1}\p_{z}v(t)\|$. However, terms with normal derivative is not controlled by conormal Sobolev spaces, and following sections will explain that optimal regularity of $\p_{z}v(t)$ is $L^{\infty}_{T}X^{m-2,0}$ or $L^{4}_{T}X^{m-1,0}$. Therefore, we cannot close energy estimate with $\|\p_{t}^{m-1}|\nabla_{y}| v(t)\|^{2}$ without time integration.  \\
Instead, we sharply split above term as
\begin{equation} \label{sharp time bdry}
\int_{S} \p_{t}^{m-1}v \p_{t}^{m-1}\p_{z} q \bigg\vert_{t}  \ \sim  \ \|\p_{t}^{m-1}|\nabla_{y}| v(t)\|^{2} + \||\nabla_{y}|^{-1}\p_{t}^{m-1}\p_{z}q(t)\|^{2} + (\text{other terms}) , 
\end{equation}
where $|\nabla_{y}|$ is horizontal Fourier muliplier. The first terms is controlled by $\|v\|^{2}_{X^{m-1,1}}$, so is absorbed by the conormal energy terms. However, $|\nabla_{y}|^{-1}$ could be insufficient to reduce order of the second term $\||\nabla_{y}|^{-1}\p_{t}^{m-1}\p_{z}q(t)\|^{2}$ because of $\p_{z}$. So we use the fact that we can split pressure $q$ into Euler part $q^{E}$, Navier-Stokes part $q^{NS}$, and Surface tension part $q^{S}$. These are defined in (\ref{q_E}), (\ref{q_NS}), and (\ref{q_S}), respectively. We can apply (\ref{time bdry}) for $q^{E}$ because estimate of $q^{E}$ is not harmful in terms of regularity of $h$. For $q^{NS}$ and $q^{S}$, we use sharp split (\ref{sharp time bdry}) and the fact that $q^{NS}$ and $q^{S}$ solve harmonic equation in the lower half plane. To explain this idea, let us consider simplified the model problem in the lower half plane which resembles $q^{S}$ in (\ref{press eq}).
\begin{equation}
\Delta g = 0, \quad g\vert_{z=0} = \Delta_{y} h, \quad \text{where}\quad  z<0,
\end{equation}
where $\Delta_{y}$ is horizontal laplacian which is simplified term for surface tension. By horizontal Fourier transform, $\hat{g}$ solves
\begin{equation}
\begin{split}
|\xi|^{2}\hat{g} &= \p_{zz} \hat{g}, \\
\hat{g}  &= \hat{g} (\xi,z=0) e^{|\xi| z},  \\
\p_{z} \hat{g} &= |\xi|^{3} \hat{h} e^{\xi z} 
\end{split}
\end{equation}
This implies $\p_{z}$ is changed into horizontal derivatives in the case of harmonic function. And we do not need full fourier multiplier $|\nabla|$ to reduce Sobolev order as following. We have
\begin{equation*}
\begin{split}
\||\nabla_{y}|^{-1} \p_{t}^{m-1}\p_{z} g\|^{2} &\leq  \int_{-\infty}^{0}\int_{\p S} \frac{1}{|\xi|^{2}} |\xi|^{6}|\p_{t}^{m-1} \hat{h}|^{2} e^{2|\xi| z}  \\
&\lesssim \big| |\xi|^{3/2}\p_{t}^{m-1}\hat{h} \big|_{L^{2}(\mathbb{R}^{2})}^{2} \\
&\lesssim |\p_{t}^{m-1}h|_{\frac{3}{2}}^{2} \sim |h|^{2}_{X^{m,\frac{1}{2}}},  \\
\end{split}
\end{equation*}
which is one lesser order than estimate $\|\p_{t}^{m-1}\p_{z} g\|^{2} \sim |h|^{2}_{X^{m,\frac{3}{2}}}$. We also have to perform similar estimate for $q^{NS}$ to make $\varepsilon\|\nabla \p_{t}^{m-1}v(t)\|$ as small dissipation type,
\begin{equation} 
\begin{split}
\||\nabla_{y}|^{-1} \p_{t}^{m-1}\p_{z} {q}^{NS}(T)\|
&\leq \varepsilon T \int_{0}^{T} \|\nabla \p_{t}^{m} v(t)\|^{2} dt + (\text{good terms}),\quad T \ll 1.  \\
\end{split}
\end{equation}
By combining $E_{0}(t), E_{1}(t), \cdots, E_{m-1}(t)$, and $E_{m}(t)$ properly (with some small weight for $E_{m}(t)$), we can close energy estimates in terms of regularity of $h$ for sufficiently small $T \ll 1$.  \\

\subsubsection{$L^2$-type normal estimate}
The problems from bad commutators in terms of $h$ was resolved by above steps. However we will see another bad commutators which contain $\left\|\partial_z v\right\|_{X^{m-1,0}}$. This cannot be controlled, since conormal Sobolev space $H^{1}_{co}$ is weaker than standard Sobolev space $H^{1}$. Therefore we should construct another energy estimate for $\nabla v$. However equation of $\nabla v$ is not clear. Instead, we make an energy estimates of new variable: 
\[
 S_n := \Pi (S^\varphi v \mathbf{n} ),\,\,\,\,\,\,\text{where} \,\,\,\,\,\,\,\Pi := I - \mathbf{n} \otimes \mathbf{n},\quad \text{tangential projection onto free-boundary}.
\] 
We can show that this quantity is equivalent to $\p_z v$. And $S_n$ has very good boundary condition $S_n\vert_{z=0} = 0$ by boundary condition (\ref{1.15}). Therefore, if we get energy estimate for $\|S_n\|_{m-1}$, we can close energy estimate. Unfortunately, equation of $S_n$ yields $(D^{\varphi})^2 q$ as pressure term. Therefore, we cannot use advantage of divergence free of velocity in energy estimate and we cannot reduce the order of pressure. Finally, the optimal regularity of $S_n$ will be $m-2$, not $m-1$. Considering pressure estimate, we get the estimate of $S_n$:
\begin{equation*}
\left\|S_n (T)\right\|^2_{X^{m-2,0}} + \varepsilon\int_0^T \|\nabla S_n \|^2_{X^{m-2,0}} dt \leq C_0 + \Lambda(R)\int_0^T \big(E(t) + \|S_n\|^{2}_{X^{m-2,0}} + |h|^{2}_{X^{m-2,\frac{7}{2}}} \big) dt.
\end{equation*} 
Since the last term $|h|^{2}_{X^{m-2,\frac{7}{2}}}$ is not in energy $E(t)$, we use Dirichlet-Neumann estimate again to close energy estimate in terms of $h$-regularity. We fail to get $m-1$ estimate for $\partial_z v$. But, above estimate can be used to control $L^\infty$ type estimate.

\subsubsection{$L^\infty$-type normal estimate}
Next, we consider estimates for $L^\infty$-type terms, which is included in $R$ above. In Sobolev conormal space, we cannot use standard Sobolev embedding, because of its normal derivative vanishes on the boundary. We use $S_n$ instead of $\partial_z v$. The main difficulty is the commutator between $Z_3$ and Laplacian. We consider a thin layer near the boundary and reparameterize so that $\partial_z^\varphi\partial_z^\varphi$ look like $\partial_{zz}$. And then, we change the advection term as
$$
\partial_t\rho + (w_y(t,y,0),zw_3(t,y,0))\cdot\nabla\rho - \varepsilon\partial_{zz}\rho = l.o.t.
$$
We do not apply a simple maximum principle for convection-diffusion equations. We apply Duhamel's principle using Green's evolution kernel. Then we can conclude
$$
\|\partial_z v\|_{Y^{[\frac{m}{2}],0}}^2 \lesssim \|\partial_v(0)\|_{Y^{[\frac{m}{2}],0}}^2 + \Lambda(R)\int_0^t \varepsilon\|\partial_{zz} v\|_{X^{m-2,0}}^2,
$$
for sufficiently large $m$.

\subsubsection{Vorticity estimate.}
We still have to control $\|\partial_z v\|_{X^{m-1,0}}$, since we failed to control by $S_{n}$. We use critical idea of \cite{NMFR}, vorticity estimate, which is also equivalent to $\partial_z v$. The biggest advantage of vorticity is that taking curl removes pressure term in Navier-Stokes equation. This makes us possible to get $m-1$ order estimate. However, because we are considering general vorticity data on the free-boundary, standard energy estimate yields high order boundary integral, which cannot be controlled. Instead, in \cite{NMFR}, $L^4_t$ type estimate for vorticity was developed, based on the heat equation with general boundary data and paradifferential calculus. This idea makes us possible to a get an estimate of $\|\omega\|_{L^4 X^{m-1}}$, where $\omega := \nabla^{\varphi}\times v$. Although this is not $L^{\infty}_{T}$ type energy, but we are suffice to get $L^{p}_{T}X^{m-1,0}, \ (p>2)$ type estimate, because what we should control is $\|\nabla v\|_{L^{2}_{T}X^{m-1,0}}$ from the right hand side of (\ref{E0}), for instance.

\subsubsection{Uniform Existence, Uniqueness, and Vanishing viscosity limit}
At this point we  have made all the necessary estimates to close the main energy estimate (\ref{main est}). In particular, the right hand side of the energy estimate is independent of kinematic viscosity $\varepsilon$, provided the energy remains bounded. So, using the preliminary existence result of A.Tani \cite{AT} and strong compactness arguments, we get local existence. For uniqueness, it is suffice to do $L^2$-estimate for the difference of two solutions  and we conclude by Gronwall's inequality. For vanishing viscosity limit, we have $L^{2}$ weak convergence from compactness argument. Moreover, $L^{2}$ norm convergence is gained by zero order $L^{2}$ energy estimate. Finally we have we have strong $L^{2}$ convergence and the unique limit solves free-boundary Euler equation with surface tension.

\subsection{Comparing the problem with and without surface tension}

Surface tension is, overall, a regularizing force in the water wave problem; however, it introduces several (perhaps unexpected) difficulties. Here we want to elaborate upon the differences between the paper of Masmoudi-Rousset \cite{NMFR} (the case where no there is no surface tension) and our result (where surface tension is taken into account. 
In terms of the basic functional framework, both works use Sobolev co-normal spaces due to the presence of boundary layers. However, there are big differences between these two works. First, let us look at a scheme of \cite{NMFR} (no surface tension case). When we have no surface tension, $m$-order energy estimate contains $|h|_m$ as its energy. The main problem which the authors faced in \cite{NMFR} is the presence of certain high order commutators. To get around this problem, the authors made use of Alinhac's good unknown which allowed them to close the energy estimates. They use the good unknowns: $V^\alpha = Z^\alpha v - \partial_z^\varphi v Z^\alpha \eta$ and $Q^\alpha = Z^\alpha q - \partial_z^\varphi q Z^\alpha \eta$, because, with this new variable, the bad commutator $Z^\alpha N$ disappears.  \\ 
\indent Meanwhile, in the surface tension case (this paper), the $m^{th}$-order energy estimate has $|h|_{m+1}$ in its energy. So now one does not need Alinhac's good unknown. Nevertheless, we also lack $\frac{1}{2}$ order ( $|h|_{m+\frac{3}{2}}$ appears in the commutators) because of the pressure. Because we use co-normal spaces,
$$
\int_S Z^m v\cdot \nabla^\varphi Z^m q 
$$
make high order commutator about pressure $q$ in $S$ (which vanishes in case of standard Sobolev space derivatives $D^m$, by divergence free condition). Since $q^b \sim \Delta h$, $q \sim \partial_x^\frac{3}{2} h$. As mentioned in above scheme, it is bounded by taking time derivatives. The crucial point is that when we only take time derivatives of the equation the worst commutator does not show up. Meanwhile, lack of regularity of $S_n$ also appear in our case. 
Therefore we should use the idea of vorticity estimate in \cite{NMFR} also.  \\
\indent Regarding $L^\infty$ estimates for $S_n$, \cite{NMFR} requires $\varepsilon
\|\partial_{zz} v\|_{L^\infty}$. But we do not need $\varepsilon\|\partial_{zz} v\|_{Y^{k,0}}$. This is because, $\varepsilon
\|\partial_{zz} v\|_{L^\infty}$ appears by Alinhac's unknown which include $\partial_z v Z^\alpha \eta$.

\section{Basic Propositions}
\subsection{Basic propositions}
We construct some propositions to estimate commutators.
\begin{proposition} \label{commutator est}
Let $u,v \in X^{m,0} \cap Y^{[\frac{m}{2}],0}(S)$. For $m \geq 2$, we get the following estimates.
\begin{eqnarray*} 
\left\| Z^{m}(uv)\right\| &:=& \left\| Z^{m}(uv)\right\|_{L^2(S)} \\
&\lesssim& 
\left\|u\right\|_{X^{m,0}} \left\|v\right\|_{Y^{[\frac{m}{2}],0}} + 
\left\|v\right\|_{X^{m,0}} \left\|u\right\|_{Y^{[\frac{m}{2}],0}},    \\
\left\|[Z^{m},u]v\right\| &:=& \left\|Z^{m}(uv) - u (Z^m v)\right\|_{L^2(S)}\\ &\lesssim& 
\left\|u\right\|_{X^{m,0}}\left\|v\right\|_{Y^{[\frac{m}{2}],0}} +
\left\|v\right\|_{X^{m-1,0}}\left\|u\right\|_{Y^{[\frac{m}{2}],0}},   \\
\left\|[Z^{m},u,v]\right\| &:=& \left\|Z^{m}(uv) - (Z^m u)v - u(Z^m v)\right\|_{L^2(S)} \\
&\lesssim& 
\left\|u\right\|_{X^{m-1,0}}\left\|v\right\|_{Y^{[\frac{m}{2}],0}} +
\left\|v\right\|_{X^{m-1,0}}\left\|u\right\|_{Y^{[\frac{m}{2}],0}},
\end{eqnarray*}
where $Z^m$ is defined in (\ref{m conormal}) and function space $X$ and $Y$ are defined in Definition \ref{XY space}. Also note that $[\frac{m}{2}]$ is defined in Defintion \ref{bracket}.
\end{proposition}

\begin{proof} 
We cannot use general Leibnitz Rule since $Z_3 = \frac{1}{1-z}\partial_z$, but every order of derivatives of $\frac{1}{1-z}$ is uniformly bounded for $z < 0$.
\begin{eqnarray*}
	\left\| Z^m(uv) \right\|_{L^2} &\leq& \sum_{ \tiny \begin{array}{c}
	(k_{1},\beta,k_{2},\gamma), \\ k_{1}+|\beta|+k_{2}+|\gamma|\leq m \end{array} } \| \p_t^{k_{1}} Z^\beta u \ \p_t^{k_{2}} Z^\gamma v \|_{L^2}   \\
	&=& \sum_{k_{1}+|\beta|\geq k_{2}+|\gamma|} \| \p_t^{k_{1}}Z^\beta u \ \p_t^{k_{2}} Z^\gamma v \|_{L^2}  + \sum_{k_{1}+|\beta|< k_{2}+|\gamma|} \| \p_t^{k_{1}}Z^\beta u \ \p_t^{k_{2}} Z^\gamma v \|_{L^2}    \\
	&\leq& \sum_{k_{1}+|\beta|\geq k_{2}+|\gamma|} \| \p_t^{k_{1}}Z^\beta u \|_{L^2} \| \p_t^{k_{2}} Z^\gamma v \|_{L^{\infty}}  + \sum_{k_{1}+|\beta|< k_{2}+|\gamma|} \| \p_t^{k_{1}}Z^\beta u \|_{L^{\infty}} \| \p_t^{k_{2}} Z^\gamma v \|_{L^2}   \\
	&\lesssim& \left\|u\right\|_{X^{m,0}} \left\|v\right\|_{Y^{[\frac{m}{2}],0}} + 
	\left\|v\right\|_{X^{m,0}} \left\|u\right\|_{Y^{[\frac{m}{2}],0}},  \\
\end{eqnarray*}	
\begin{eqnarray*}
	\left\| [Z^m,u]v \right\|_{L^2} &\leq& \sum_{ \tiny \begin{array}{c}
			(k_{1},\beta,k_{2},\gamma), \\ k_{1}+|\beta|+k_{2}+|\gamma|\leq m, \\ k_{1}+|\beta| \geq 1 \end{array} } \| \p_t^{k_{1}} Z^\beta u \ \p_t^{k_{2}} Z^\gamma v \|_{L^2}   \\
	&=& \sum_{ \tiny \begin{array}{c}
	k_{1}+|\beta|\geq k_{2}+|\gamma|, \\ k_{1}+|\beta| \geq 1 \end{array} } \| \p_t^{k_{1}}Z^\beta u \ \p_t^{k_{2}} Z^\gamma v \|_{L^2}  + \sum_{ \tiny \begin{array}{c}
	k_{1}+|\beta|\geq k_{2}+|\gamma|, \\ k_{1}+|\beta| \geq 1 \end{array} } \| \p_t^{k_{1}}Z^\beta u \ \p_t^{k_{2}} Z^\gamma v \|_{L^2}    \\
	&\lesssim& \left\|u\right\|_{X^{m,0}} \left\|v\right\|_{Y^{[\frac{m}{2}],0}} + 
	\left\|v\right\|_{X^{m-1,0}} \left\|u\right\|_{Y^{[\frac{m}{2}],0}},   \\
\end{eqnarray*}	
\begin{eqnarray*}	
	\left\| [Z^m,u,v] \right\|_{L^2} &\leq& \sum_{ \tiny \begin{array}{c}
			(k_{1},\beta,k_{2},\gamma), \\ k_{1}+|\beta|+k_{2}+|\gamma|\leq m, \\ k_{1}+|\beta| \geq 1, \\ k_{2}+|\gamma| \geq 1 \end{array} } \| \p_t^{k_{1}} Z^\beta u \ \p_t^{k_{2}} Z^\gamma v \|_{L^2}   \\
	&=& \sum_{ \tiny \begin{array}{c}
			k_{1}+|\beta|\geq k_{2}+|\gamma|, \\ k_{1}+|\beta| \geq 1, \\ k_{2}+|\gamma| \geq 1 \end{array} } \| \p_t^{k_{1}}Z^\beta u \ \p_t^{k_{2}} Z^\gamma v \|_{L^2}  + \sum_{ \tiny \begin{array}{c}
			k_{1}+|\beta|\geq k_{2}+|\gamma|, \\ k_{1}+|\beta| \geq 1, \\ k_{2}+|\gamma| \geq 1 \end{array} } \| \p_t^{k_{1}}Z^\beta u \ \p_t^{k_{2}} Z^\gamma v \|_{L^2}    \\
			&\lesssim& \left\|u\right\|_{X^{m-1,0}} \left\|v\right\|_{Y^{[\frac{m}{2}],0}} + 
			\left\|v\right\|_{X^{m-1,0}} \left\|u\right\|_{Y^{[\frac{m}{2}],0}}.
\end{eqnarray*}
\end{proof}

\begin{remark}
The idea is that for each bilinear term, we put the $L^2$ norm on the term with higher derivatives and  the $L^\infty$ norm to low order term. In co-normal derivatives, there is no proper notion of fractional derivatives, so $Z_3^{1/2}$ does not make sense. We deal when $m$ is even, so that $\frac{m}{2}$ is also a integer, but our result also work for odd $m$, because it suffices to give $\frac{m-1}{2}$ orders to $L^\infty$ and $\frac{m+1}{2}$ orders to $L^2$.  
\end{remark}

The following proposition is for anisotropic embeddings and trace properties in co-normal spaces.
\begin{proposition} \label{aniso embed}
1) For $s_1\geq 0$, $s_2\geq 0$ such that $s_1+s_2 >2$ and $u$ such that $u\in H^{s_1}_{tan}(S)$, $\partial_zu\in H^{s_2}_{tan}(S)$, we have the anisotropic Sobolev embedding:
\begin{equation*}
\left\|u\right\|_{L^\infty}^2
\lesssim \left\|\partial_zu\right\|_{H^{s_2}_{tan}}\left\|u\right\|_{H^{s_1}_{tan}}.
\end{equation*}
2) For $u\in H^1(S)$, we have the trace estimate :
\begin{equation*}
\left|u(\cdot,0)\right|_{H^s(\mathbb{R}^2)}
\lesssim \left\|\partial_zu\right\|^{\frac{1}{2}}_{H^{s_2}_{tan}}
\left\|u\right\|^{\frac{1}{2}}_{H^{s_1}_{tan}},
\end{equation*}
with $s_1+s_2=2s\geq 0$, where $H^{m}_{tan}$ is defined in Definition \ref{space def}. Especially,
\begin{equation*}
	|u^b|_{H^{s+\frac{1}{2}}} \lesssim \|\p_z u\|^2_{H^{s+1}_{tan}} + \|u\|^2_{H^{s}_{tan}}
\end{equation*}
\end{proposition}
\begin{proof}
For 1) and 2), see Proposition 2.2 in \cite{NMFR}. 
\end{proof}

\begin{proposition} \label{horizon est}
We state classical commutator estimate on $\mathbb{R}^{2}$. Let $\Lambda$ is the Fourier mulitplier $(1+|\xi|^{2})^{\frac{s}{2}}$. Then we have the following estimate.
\begin{equation}
\begin{split}
	|\Lambda^{s}(fg)| &\leq C_{s}\big( |f|_{L^{\infty}(\mathbb{R}^{2})} |g|_{H^{s}(\mathbb{R}^{2})} + |g|_{L^{\infty}(\mathbb{R}^{2})} |f|_{H^{s}(\mathbb{R}^{2})} \big),  \\
	|[\Lambda^{s},f]\nabla g|_{L^{2}(\mathbb{R}^{2})} &\leq C_{s}\big( |\nabla f|_{L^{\infty}(\mathbb{R}^{2})} |g|_{H^{s}(\mathbb{R}^{2})} + |\nabla g|_{L^{\infty}(\mathbb{R}^{2})} |f|_{H^{s}(\mathbb{R}^{2})} \big),  \\
	|uv|_{\frac{1}{2}} &\lesssim |u|_{1,\infty} |v|_{\frac{1}{2}},  \\
	|uv|_{-\frac{1}{2}} &\lesssim |u|_{1,\infty} |v|_{-\frac{1}{2}}.  \\
\end{split}
\end{equation}
\begin{proof}
	First three estimates can be found in \cite{NMFR}. The last inequality comes from duality argument of third one.
\end{proof}
\end{proposition}

\subsection{Estimate of $\eta$}
In this subsection, we investigate regularity of mapping $\eta$, defined in (\ref{eta def}) and (\ref{extension}). As explained before, the reason we choose such smoothing diffeomorphism is that regularity of $\eta$ is $\frac{1}{2}$ better than $h$. This fact is crucial later, because this term can accommodate an extra $\frac{1}{2}$ derivative in bilinear estimates. For example, in the pressure estimates: 
$$
\int_S (\nabla\varphi) q \leq \left\|\nabla\varphi\right\|_{\frac{1}{2}}\left\|q\right\|_{-\frac{1}{2}} \sim \left|\nabla\ h\right|_{L^2}\left\|q\right\|_{-\frac{1}{2}}.
$$
We defined diffeomorphism so that at initial time, $\partial_z\varphi(0,y,z) \geq 1$. $\partial_z\varphi$ should be positive during our estimates, so our estimate is valid during on $[0,T^\varepsilon]$ such that 
\begin{equation} \label{diffeo cond}
\partial_z\varphi(t,y,z) \geq c_0,\,\,\,\,\,\forall t\in [0,T^\varepsilon]
\end{equation}
for some $c_0$, where $T^{\varepsilon}$ is time interval of solution with fixed $\varepsilon>0$.

\begin{proposition} \label{eta regularity}
For $\eta$, defined in (\ref{extension}), we obtain the following estimates.
\begin{eqnarray*}
\left\|\nabla\eta\right\|_{H^s(S)} &\leq& C_s\left|h\right|_{H^{s+\frac{1}{2}}(\p S)},  \\
\left\|\nabla\eta\right\|_{X^{m,0}} &\leq& C_s\left|h\right|_{X^{m,\frac{1}{2}}}.
\end{eqnarray*}
Moreover, for $L^\infty$ type estimates, we get
\begin{eqnarray*}
\forall s\in \mathbb{N},\quad \left\|\eta\right\|_{W^{s,\infty}} \leq C_s
\left|h\right|_{s,\infty},   \\
\forall s\in \mathbb{N},\quad \left\|\eta\right\|_{Y^{m,0}} \leq C_s\left|h\right|_{Y^{m,0}}.
\end{eqnarray*}
Note that we have standard Sobolev regularity, not co-normal one.
\end{proposition}

\begin{proof}
The first inequality is from \cite{NMFR}, and $\left|\nabla\partial_t^k\eta\right|_{H^s(S)} \leq C_s\left|\partial_t^kh\right|_{s+\frac{1}{2}}$ is also trivial by definition of $\eta$. So, by summing all cases, we get the second inequality.
For $L^\infty$ type estimates, the third inequality is from \cite{NMFR}, and the last inequality is also gained in similar way.
\end{proof}

From the definition of $\p_{i}^{\varphi}$ in Definition \ref{varphipartial}, we should control Sobolev norm of fractional term $\frac{u}{\p_z \varphi}$. The following lemma is useful to estimate such terms.
\begin{lemma} \label{fraction est}
Assume that (\ref{diffeo cond}) holds. Then we have the following estimate for $u\in X^{m,0}\cap Y^{[\frac{m}{2}],0}(S)$ and $h\in X^{m,\frac{1}{2}}\cap Y^{[\frac{m}{2}],1}(S)$.
\begin{equation*}
\left\|Z^m \frac{u}{\partial_z \varphi}\right\|_{L^2} \lesssim 
\Lambda\Big(\frac{1}{c_0},\left\|u\right\|_{Y^{[\frac{m}{2}],0}} + \left|h\right|_{Y^{[\frac{m}{2}],1}} \Big)
( \left\|u\right\|_{X^{m,0}} + \left|h\right|_{X^{m,\frac{1}{2}}} ).
\end{equation*}
\end{lemma}

\begin{proof}
$F(x)=x/(A+x)$ is a smooth function of which all order derivatives are bounded when $A+x \geq c_0 > 0$. Since $\p_z \varphi = A + \p_z \eta$, $\frac{u}{\partial_z \varphi} = \frac{u}{A}-\frac{u}{A}F(\partial_z \eta)$. Therefore,
\begin{eqnarray*}
\left\|Z^m \frac{u}{\partial_z \varphi} \right\|_{L^2}  &=& \left\|Z^m\Big(\frac{u}{A}-\frac{u}{A}F(\partial_z \eta)\Big)\right\|_{L^2}
\lesssim \left\|u\right\|_{X^{m,0}} + \left\|Z^m(uF(\partial_z\eta))\right\|_{L^2}   \\
&\lesssim& \left\|u\right\|_{X^{m,0}} + \left\|u\right\|_{X^{m,0}}\left\|F(\partial_z\eta)\right\|_{Y^{[\frac{m}{2}],0}} + \left\|u\right\|_{Y^{[\frac{m}{2}],0}}\left\|F(\partial_z\eta)\right\|_{X^{m,0}},
\end{eqnarray*}
where we used the first commutator estimate in Proposition \ref{commutator est}. Moreover, using Proposition \ref{commutator est} again, we have,
\begin{eqnarray*}
\left\|F(\partial_z\eta)\right\|_{X^{m,0}} &\lesssim&
\Lambda\Big(\frac{1}{c_0},\left\|\nabla\eta\right\|_{Y^{[\frac{m}{2}],0}}
\Big)\left\|\partial_z\eta\right\|_{X^{m,0}}
\lesssim \Lambda\Big(\frac{1}{c_0},\left|h\right|_{Y^{[\frac{m}{2}],1}}\Big)\left|h\right|_{X^{m,\frac{1}{2}}} ,    \\
\left\|F(\partial_z\eta)\right\|_{Y^{[\frac{m}{2}],0}} &\lesssim&
\Lambda\Big(\frac{1}{c_0},\left\|\nabla\eta\right\|_{Y^{[\frac{m}{2}],0}}
\Big) \lesssim \Lambda\Big(\frac{1}{c_0},\left|h\right|_{Y^{[\frac{m}{2}],1}}
\Big),
\end{eqnarray*}
where we used smoothing property of Proposition \ref{eta regularity}.
\end{proof}

\subsection{Dissipation term control} 
From the definition (\ref{1.5}), we define volume element $\d V_t := \p_{z}\varphi dxdydz$ in $S$, i.e.
\begin{eqnarray} \label{volume integral}
	\int_{\O_t} f \dd V = \int_S (f\circ\Phi) \ \p_z\varphi \ \dd x \dd y \dd z := \int_S (f\circ\Phi) \ \dd V_t .
\end{eqnarray}
First, we compare $L^2$ type norms between $\nabla v $ and $\nabla^{\varphi} v$.  \\

\begin{lemma} \label{compare}
Assume that $\p_z \varphi \geq c_0$ and $\|\nabla \varphi\|_{L^\infty} \leq \frac{1}{c_0}$ for some $c_0 > 0$, then we have,
\begin{eqnarray*}
	\|\nabla f\|^2_{X^{m,0}} \leq \Lambda\left(\frac{1}{c_0}, \|\nabla f\|_{Y^{[\frac{m}{2}],0}} + |h|_{Y^{[\frac{m}{2}],1}} \right) ( \|\nabla^{\varphi} f\|^{2}_{X^{m,0}} + |h|^2_{X^{m,\frac{1}{2}}}) .
\end{eqnarray*}
\end{lemma}
\begin{proof}
	From Definition \ref{varphipartial},
	\begin{eqnarray*}
		\|\p_z f\|_{X^{m,0}} &=& \| \p_z\varphi \ \p_z^{\varphi} f \|_{X^{m,0}} \lesssim \|\p_z\varphi\|_{X^{m,0}} \|\p_z^{\varphi} f\|_{Y^{[\frac{m}{2}],0}} + \|\p_z^{\varphi} f\|_{X^{m,0}} \|\p_z \varphi \|_{Y^{[\frac{m}{2}],0}}   \\
		&\lesssim& \Lambda\left(\frac{1}{c_0}, \|\nabla f\|_{Y^{[\frac{m}{2}],0}} + |h|_{Y^{[\frac{m}{2}],1}} \right) ( \|\p_z^{\varphi} f\|_{X^{m,0}} + |h|_{X^{m,\frac{3}{2}}} ) ,   \\
		\|\p_{i=x,y} f\|_{X^{m,0}} &=& \|\p_i^{\varphi} f \|_{X^{m,0}} + \| \p_i\varphi \ \p_z^{\varphi} f \|_{X^{m,0}} \lesssim \|\p_i^{\varphi} f \|_{X^{m,0}} + \|\p_z\varphi\|_{X^{m,0}} \|\p_i^{\varphi} f\|_{Y^{[\frac{m}{2}],0}} + \|\p_z^{\varphi} f\|_{X^{m,0}} \|\p_i \varphi \|_{Y^{[\frac{m}{2}],0}}   \\
		&\lesssim& \|\p_i^{\varphi} f \|_{X^{m,0}} + \Lambda\left(\frac{1}{c_0}, \|\nabla f\|_{Y^{[\frac{m}{2}],0}} + |\nabla h|_{Y^{[\frac{m}{2}],0}} \right) ( \|\p_z^{\varphi} f\|_{X^{m,0}} + |h|_{X^{m,\frac{1}{2}}} ) ,  \\
	\end{eqnarray*}
	where we used Proposition \ref{commutator est} and \ref{eta regularity}. Combining above two estimate, we get the result.
\end{proof}
\noindent Using Lemma \ref{compare} and (\ref{volume integral}), we have a version of Korn's inequality for $S^\varphi$.
\begin{proposition} \label{korn}
1) If $\partial_z\varphi \geq c_0,\,\,\left\|\nabla\varphi\right\|_{L^\infty} + \left\|\nabla^2\varphi\right\|_{L^\infty} \leq \frac{1}{c_0}$ for some $c_0 > 0$, then there exists $\Lambda_0 = \Lambda(\frac{1}{c_0}) > 0 $ such that for every $v\in H^1(S)$, we have
\begin{equation*}
\left\|\nabla v\right\|^2_{L^2(S)} \leq \Lambda_0\Big( \int_S\left|S^\varphi v\right|^2 dV_t + \|v\|^2\Big),
\end{equation*}
where $S^\varphi v := \frac{1}{2} (\nabla^\varphi v + (\nabla^\varphi v)^T )$.  \\
\noindent 2) For higher order case, we have the following estimate.
\begin{equation*}
\left\|\nabla v\right\|^2_{X^{m,0}} \leq \Lambda_0\Big(\int_S\left|S^\varphi v\right|^2_{X^{m,0}} \dd V_t + \|v\|^2_{X^{m,0}} \Big).
\end{equation*}
\end{proposition}
\begin{proof}
See Proposition 2.9 in \cite{NMFR} for the first estimate 1). For the second estimate 2), we apply same estimate for $\partial^k Z^\alpha v$. 
\end{proof}

\section{Equations of $ (Z^m v,Z^m h,Z^m q )$}
\noindent The aim of this section is to prove the following proposition.
\begin{proposition} \label{high order eq}
	Let positive integer $m \geq 2$. Then applying $Z^m$ to the system (\ref{1.12})-(\ref{1.15}), we get the following equations with respect to $(Z^m v,Z^m h,Z^m q)$ and commutator estimates.
\begin{equation*}
\begin{split}
 & (\partial_t^{\varphi} + v\cdot\nabla^{\varphi})(Z^{m} v) + \nabla^{\varphi}(Z^{m}q) =  2 \varepsilon \nabla^{\varphi} \cdot S^{\varphi} (Z^m v) - T^{m}(v) +  C^{m}(q) + \varepsilon\Theta^m(v) - \varepsilon\mathcal{D}^m(S^{\varphi}v)   ,  \\
 &\nabla^\varphi\cdot(Z^{m}v) = C^{m}(d) ,   \\
 &\p_t(Z^{m} h) = (Z^{m} v^b)\cdot \mathbf{N} + v^b\cdot(Z^{m} \mathbf{N}) + C^{m} (KB) , \\
 &\Big\{ Z^{m}q^b - gZ^{m}h - 2\varepsilon  (S^\varphi(Z^{m} v))^b + \nabla\cdot \frac{\nabla Z^{m} h}{\sqrt{1+|\nabla h|^2}} \Big\} \mathbf{N} = - C^{m}(B) + \varepsilon (\Theta^{m}(v))^b \mathbf{N},   \\
 &\quad\quad\quad\quad\quad\quad\quad\quad\quad\quad\quad\quad\quad\quad\quad\quad\quad\quad\quad\quad\quad + \nabla\cdot \frac{\nabla h \langle\nabla h,\nabla Z^{m} h\rangle}{\sqrt{1+|\nabla h|^2}^3}\mathbf{N} + \nabla\cdot C^{m}(S) \mathbf{N},
\end{split}
\end{equation*}
where 
\begin{equation*}
\begin{split}
\left\| C^{m}(q) \right\|
&\lesssim \Lambda\Big(\frac{1}{c_0},
\left\|\nabla q\right\|_{Y^{[\frac{m}{2}],0}} + \left|h\right|_{Y^{[\frac{m}{2}],1}} \Big)  ( \left\|\nabla q\right\|_{X^{m-1,0}} + \left|h\right|_{X^{m,\frac{1}{2}}}  ),  \\
\left\| C^{m}(d) \right\|
&\lesssim \Lambda\Big(\frac{1}{c_0},
\left\|\nabla v\right\|_{Y^{[\frac{m}{2}],0}} + \left|h\right|_{Y^{[\frac{m}{2}],1}} \Big)  ( \left\|\nabla v\right\|_{X^{m-1,0}} + \left|h\right|_{X^{m,\frac{1}{2}}} ),  \\
\left\|T^{m}(v)\right\| &\lesssim \Lambda\Big(\frac{1}{c_0},
\left\|\nabla v\right\|_{Y^{[\frac{m}{2}],0}} + \left|h\right|_{Y^{[\frac{m}{2}],1}} \Big)  ( \left\|\nabla v\right\|_{X^{m-1,0}} + \left\|v\right\|_{X^{m,0}} + \left|h\right|_{X^{m-1,\frac{1}{2}}} ),   \\
\| \Theta^{m} (v) \|
&\lesssim \Lambda\Big(\frac{1}{c_0},
\left\|\nabla v\right\|_{Y^{[\frac{m}{2}],0}} + \left|h\right|_{Y^{[\frac{m}{2}],1}} \Big) ( \left\|\nabla v\right\|_{X^{m-1,0}} + \left|h\right|_{X^{m,\frac{1}{2}}} ),  \\
\left\| C^{m} (KB) \right\| &\lesssim \Lambda
\Big( \frac{1}{c_0}, \left\|\nabla v\right\|_{Y^{[\frac{m}{2}],0}} + \left|h\right|_{Y^{[\frac{m}{2}],1}} \Big) 
( \left\|v\right\|_{X^{m-1,0}} + \left\|\nabla v\right\|_{X^{m-1,0}} + |h|_{X^{m-1,1}} ),  \\
\|C^{m}(B)\| &\lesssim 2\varepsilon \ \Lambda\Big(\frac{1}{c_0},\left|h\right|_{Y^{[\frac{m}{2}],1}} + \left\|\nabla v\right\|_{Y^{[\frac{m}{2}],0}}\Big) (\left|h\right|_{X^{m-1,1}} + |v^b|_{X^{m-1,1}} ),  \\ 
\|\nabla\cdot C^{m}(S)\mathbf{N}\| &\lesssim \Lambda\Big(\frac{1}{c_0},\left|h\right|_{Y^{[\frac{m}{2}],2}} + \left\|\nabla v\right\|_{Y^{[\frac{m}{2}],1}}\Big) |h|_{X^{m-1,2}} .    \\
\end{split}
\end{equation*}
\end{proposition}
\begin{proof}
Proof of this proposition is given by (\ref{321})-(\ref{332}), and (\ref{335})-(\ref{338}), those are derived throughout following three subsections.
\end{proof}

\subsection{Commutator estimate}
\begin{proposition} \label{eq comm est}
For $i=t,1,2,3$, let us define $C_i^{m} (f)$ by
\begin{equation} \label{def basic commutator}
Z^{m}(\partial_i^{\varphi}f) := \partial_i^{\varphi}(Z^{m} f) + 
C_i^{m} (f).
\end{equation}
Then we have the following commutator estimate for $C_i^{m} (f)$.
\begin{equation*}
\left\|C_i^{m} (f)\right\| \lesssim \Lambda\Big(\frac{1}{c_0},
\left\|\nabla f\right\|_{Y^{[\frac{m}{2}],0}} + \left|h\right|_{Y^{[\frac{m}{2}],1}} \Big) \Big( \left\|\nabla f\right\|_{X^{m-1,0}} + \left|h\right|_{X^{m,\frac{1}{2}}} \Big).
\end{equation*}
\end{proposition}

\begin{proof}
For $i=t,1,2$
\begin{eqnarray*}
Z^{m}\Big(\partial_i f-\frac{\partial_i\varphi}{\partial_z\varphi}\partial_z f\Big)
&=& \partial_i (Z^{m}f )-Z^{m}\Big(\frac{\partial_i\varphi}{\partial_z\varphi}\partial_i^\varphi f \Big)   \\
&=&\partial_i (Z^{m}f)-\Big(\left[Z^{m},\frac{\partial_i\varphi}{\partial_z\varphi},\partial_z f\right] + \Big(Z^{m}\frac{\partial_i\varphi}{\partial_z\varphi}\Big)\partial_z f + \frac{\partial_i\varphi}{\partial_z\varphi}\Big(Z^{m}\partial_z f\Big)\Big)   \\
&=&\partial_i^\varphi (Z^{m}f) \underbrace{ -  \Big(\left[Z^{m},\frac{\partial_i\varphi}{\partial_z\varphi},\partial_z f\right] + \Big(Z^{m}\frac{\partial_i\varphi}{\partial_z\varphi}\Big)\partial_z f + \frac{\partial_i\varphi}{\partial_z\varphi}\left[Z^{m},\partial_z\right]f \Big) }_{:= C^{m}_{i}(f)}
 \\
&:=& \partial_i^\varphi (Z^{m}f) + C_i^{m}(f).
\end{eqnarray*}
Now we estimate three terms of $C^{m}_{i}(f)$ using Propositions \ref{commutator est}, \ref{eta regularity}, and Lemma \ref{fraction est}.
\begin{eqnarray*}
\left\|\left[Z^{m},\frac{\partial_i\varphi}{\partial_z\varphi},\partial_z f\right]\right\| &\lesssim& \left\|\frac{\partial_i\varphi}{\partial_z\varphi}\right\|_{X^{m-1,0}}\left\|\partial_z f\right\|_{Y^{[\frac{m}{2}],0}} + \left\|\frac{\partial_i\varphi}{\partial_z\varphi}\right\|_{Y^{[\frac{m}{2}],0}}\left\|\partial_z f\right\|_{X^{m-1,0}}    \\
&\lesssim& \Lambda\Big(\frac{1}{c_0},\left|h\right|_{Y^{[\frac{m}{2}],1}} + \left\|\partial_z f\right\|_{Y^{[\frac{m}{2}],0}}\Big)\Big(\left|h\right|_{X^{m-1,\frac{1}{2}}} + \left\|\partial_z f\right\|_{X^{m-1,0}}\Big) ,   \\
\left\|\Big(Z^{m}\frac{\partial_i\varphi}{\partial_z\varphi}\Big)\partial_z f\right\|
&\lesssim&
\left|\partial_z f\right|_{L^\infty}\left\|Z^{m}\frac{\partial_i\varphi}{\partial_z\varphi}\right\| \lesssim \Lambda\Big(\frac{1}{c_0},\left|h\right|_{Y^{[\frac{m}{2}],1}} + \left|\partial_z f\right|_{L^\infty} \Big) \left|h\right|_{X^{m,\frac{1}{2}}},    \\
\left\|\frac{\partial_i\varphi}{\partial_z\varphi}\left[Z^{m},\partial_z\right]f\right\|
&\lesssim&
\left|\frac{\partial_i\varphi}{\partial_z\varphi}\right|_{L^\infty}
\|\sum_{k+|\beta|\leq m-1}c_\beta \partial_z(\p_t^{k}Z^{\beta}f)\|
\lesssim \Lambda\Big(\frac{1}{c_0},\left\|\partial_z f\right\|_{Y^{[\frac{m}{2}],0}} + \left|h\right|_{1,\infty} \Big)
\left\|\partial_z f\right\|_{X^{m-1,0}}.
\end{eqnarray*}
By summing these three terms, we get the results. For $i=3$,
\begin{eqnarray*}
Z^{m}\Big(\frac{\partial_z f}{\partial_z\varphi}\Big)
&=& \left[Z^{m},\frac{1}{\partial_z\varphi},\partial_z f\right]
+ \Big(Z^{m}\frac{1}{\partial_z\varphi}\Big)\partial_z f 
+ \frac{1}{\partial_z\varphi}\Big(Z^{m}\partial_z f\Big)    \\
&=& \partial_3^\varphi(Z^{m}) f + \left[Z^{m},\frac{1}{\partial_z\varphi},\partial_z f\right]
+ \Big(Z^{m}\frac{1}{\partial_z\varphi}\Big)\partial_z f + \frac{1}{\partial_z\varphi}\left[Z^{m},\partial_z\right]f.
\end{eqnarray*}
We just replace $\partial_i\varphi$ as 1, so the control is same as $i=t,1,2$ cases.
\end{proof}

\subsection{Interior Equation (\ref{1.12}) and (\ref{1.13})}
Now we apply $Z^m$ to each terms in (\ref{1.12})-(\ref{1.15}) and use Proposition \ref{eq comm est} for each commutators. 
\subsubsection{Pressure} 
\begin{equation} \label{321}
Z^{m}(\nabla^{\varphi}q)=\nabla^{\varphi}(Z^{m}q) - (C^{m}_1(q),C^{m}_2(q),C^{m}_3(q))
:= \nabla^{\varphi}(Z^{m}q)-C^{m}(q),
\end{equation}
where $C^{m}(q) := (C^{m}_1(q),C^{m}_2(q),C^{m}_3(q))$. By Proposition \ref{eq comm est},
\begin{equation} \label{322}
\left\| C^{m}(q) \right\|
\lesssim \Lambda\Big(\frac{1}{c_0},
\left\|\nabla q\right\|_{Y^{[\frac{m}{2}],0}} + \left|h\right|_{Y^{[\frac{m}{2}],1}} \Big) \Big( \left\|\nabla q\right\|_{X^{m-1,0}} + \left|h\right|_{X^{m,\frac{1}{2}}} \Big).
\end{equation}

\subsubsection{Divergence-free}
\begin{equation} \label{323}
Z^{m}(\nabla^{\varphi}\cdot v)
=\nabla^{\varphi}\cdot (Z^{m}v) - \sum_{i=1}^{3}{C^{m} _i (v)}
=\nabla^\varphi\cdot(Z^{m}v) - C^{m}(d),
\end{equation}
where $C^{m}(d) := \sum_{i=1}^{3}{C^{m} _i (v)}$. By Proposition \ref{eq comm est},
\begin{equation} \label{324}
\left\| C^{m}(d) \right\|
\lesssim \Lambda\Big(\frac{1}{c_0},
\left\|\nabla v\right\|_{Y^{[\frac{m}{2}],0}} + \left|h\right|_{Y^{[\frac{m}{2}],1}} \Big) \Big( \left\|\nabla v\right\|_{X^{m-1,0}} + \left|h\right|_{X^{m,\frac{1}{2}}} \Big).
\end{equation}

\subsubsection{Transport term}
Using divergence free condition, we have
\begin{equation*}
\partial^{\varphi}_t+(v\cdot\nabla^{\varphi})
= \partial_t+(v_y\cdot\nabla_y) + \frac{1}{\partial_z\varphi}\Big(v\cdot \mathbf{N} - \partial_t\varphi\Big)\partial_z,
\end{equation*}
where extended normal vector $\mathbf{N}$ is defined in Definition \ref{extendnormal}. Applying $Z^m$,
\begin{equation} \label{325}
Z^{m} (\partial_t^{\varphi} + v\cdot\nabla^{\varphi} )v
=  (\partial_t^{\varphi} + v\cdot\nabla^{\varphi} )(Z^{m} v) + T^{m}(v),
\end{equation}
where 
\begin{equation*}
T^{m}(v) := 
\sum_{i=1}^2\left\{\partial_i v \cdot Z^{m} v_i + [Z^{m},v_i,\partial_i v]\right\} + [Z^{m},V_z,\partial_z v] + (Z^{m} V_z)\cdot \partial_z v + V_z [Z^{m},\partial_z]v,
\end{equation*}
where 
\begin{equation} \label{Vz}
V_z := \frac{1}{\partial_z\varphi} (v\cdot \mathbf{N} - \partial_t\varphi )
\end{equation}
Using Proposition \ref{eq comm est}, we have
\begin{equation} \label{326}
\left\|T^{m}(v)\right\| \lesssim \Lambda\Big(\frac{1}{c_0},
\left\|\nabla v\right\|_{Y^{[\frac{m}{2}],0}} + \left|h\right|_{Y^{[\frac{m}{2}],1}} \Big) \Big( \left\|\nabla v\right\|_{X^{m-1,0}} + \left\|v\right\|_{X^{m,0}} + \left|h\right|_{X^{m-1,\frac{1}{2}}} \Big).
\end{equation}

\subsubsection{Diffusion}
\begin{equation*}
2 \varepsilon Z^{m} \nabla^\varphi \cdot (S^\varphi v)
= 2 \varepsilon \nabla^\varphi \cdot Z^{m} (S^\varphi v) - \varepsilon\mathcal{D}^{m} (S^\varphi v),
\end{equation*}
where $\mathcal{D}^{m} (S^\varphi v)_i := 2C_j^{m} (S^\varphi v)_{ij}$ and
\begin{equation*}
2Z^{m} (S^\varphi v)
=Z^{m} (\partial_i^\varphi v_j + \partial_j^\varphi v_i)
=2S^\varphi(Z^{m} v) + (C_i^{m} (v_j)+C_j^{m} (v_i))
:=2S^\varphi(Z^{m} v) + \Theta^{m} (v),
\end{equation*}
where $\Theta^{m} (v)_{ij} := C_i^{m} (v_j)+C_j^{m} (v_i)$.
Therefore, 
\begin{equation} \label{327}
	2 \varepsilon Z^{m} \nabla^\varphi \cdot (S^\varphi v) = 2 \varepsilon \nabla^{\varphi} \cdot S^{\varphi} (Z^m v) + \varepsilon\Theta^m(v) - \varepsilon\mathcal{D}^m(S^{\varphi}v) 
\end{equation}
the estimate of $ \Theta^{m} (v) $ is same as $ C^{m}(v) $,
\begin{equation} \label{328}
\| \Theta^{m} (v) \|
\lesssim \Lambda\Big(\frac{1}{c_0},
\left\|\nabla v\right\|_{Y^{[\frac{m}{2}],0}} + \left|h\right|_{Y^{[\frac{m}{2}],1}} \Big) \Big( \left\|\nabla v\right\|_{X^{m-1,0}} + \left|h\right|_{X^{m,\frac{1}{2}}} \Big).
\end{equation}


\subsection{Boundary Condition (\ref{1.14}) and (\ref{1.15})} 
Note that $\alpha_3=0$, because we deal functions on boundary. 
\subsubsection{Kinematic boundary condition} 
Applying $Z^m$ to (\ref{1.14}),
\begin{equation} \label{331}
\p_t(Z^{m} h)-(Z^{m} v^b)\cdot \mathbf{N} - v^b\cdot(Z^{m} \mathbf{N}) - C^{m} (KB) = 0
,\quad \text{where} \quad C^{m} (KB) := [Z^{m},v^b,\mathbf{N}].
\end{equation}
\begin{equation} \label{332}
\begin{split}
\left\| C^{m} (KB) \right\| &= \|[Z^{m},v^b,\mathbf{N}]\| \lesssim \Lambda
\Big( |v^b|_{Y^{[\frac{m}{2}],0}} + \left|\mathbf{N}\right|_{Y^{[\frac{m}{2}],0}} \Big) 
\Big( |v^b|_{X^{m-1,0}} + \left|\mathbf{N}\right|_{X^{m-1,0}} \Big)   \\
&\lesssim \Lambda
\Big( \frac{1}{c_0}, \left\|\nabla v\right\|_{Y^{[\frac{m}{2}],0}} + \left|h\right|_{Y^{[\frac{m}{2}],1}} \Big) 
( \left\|v\right\|_{X^{m-1,0}} + \left\|\nabla v\right\|_{X^{m-1,0}} + |h|_{X^{m-1,1}} ) ,  \\
\end{split}
\end{equation}
where we used trace inequality Proposition \ref{aniso embed}.

\subsubsection{Continuity of Stress tensor}
\begin{lemma} \label{lem 3.2}
We have the following estimate to control $\nabla v(\cdot,0):=(\nabla v)^b := (\nabla v)\vert_{z=0}$ by $v^b$.
\begin{equation*}
\left|\nabla v(\cdot,0)\right|_{X^{m,0}} \lesssim \Lambda\Big(
\frac{1}{c_0}, |h|_{Y^{[\frac{m}{2}],1}} + \|v\|_{Y^{[\frac{m}{2}],1}}\Big) \Big( \left|h\right|_{X^{m,1}} + \left|v(\cdot,0)\right|_{X^{m,1}} \Big).
\end{equation*}
\end{lemma}
\begin{proof}
We divide $\nabla v$ into normal part and tangential part. For tangential derivatives $\partial_1 v$ and $\partial_2 v$, results are obvious, since tangential derivatives commute with $\nabla v(\cdot,0)$. So, we only focus on $\partial_z v$.
Firstly, from the divergence free condition $\nabla^\varphi\cdot v=0$,
\begin{equation*}
\partial_z v\cdot {\mathbf{n}} = \frac{1}{|\mathbf{N}|} (A+\partial_z\eta) (\partial_1 v_1 + \partial_2 v_2),\quad \text{where $\mathbf{N}$ is defined in Definition \ref{extendnormal}}.
\end{equation*}
Using Proposition \ref{commutator est}, \ref{eta regularity}, and \ref{fraction est}, we get normal part estimate.
\begin{equation} \label{normal part est}
\begin{split}
\left|\partial_z v\cdot  \mathbf{n}\right|_{X^{m,0}} &\lesssim \Lambda\Big(\frac{1}{c_0}, \left|h\right|_{Y^{[\frac{m}{2}],1}} + \|v\|_{Y^{[\frac{m}{2}],1}}\Big)
\Big(\left|\partial_z\eta(\cdot,0)\right|_{X^{m,0}} + \left|v(\cdot,0)\right|_{X^{m,1}}\Big)  \\
&\lesssim \Lambda\Big(\frac{1}{c_0}, \left|h\right|_{Y^{[\frac{m}{2}],1}} + \left\|v\right\|_{Y^{[\frac{m}{2}],1}}\Big)
\Big(\left|h\right|_{X^{m,1}} + \left|v(\cdot,0)\right|_{X^{m,1}}\Big),\quad\text{on}\quad \p S.
\end{split}
\end{equation}
To estimate tangential part, note that   \\
\begin{equation*}
2S^\varphi v \mathbf{N} = \frac{1}{\partial_z\varphi}\Big(1+|\nabla h|^2\Big)(\partial_z v)
- \Big( \partial_1 h \partial_1 v + \partial_2 h \partial_2 v \Big)
+ \begin{pmatrix}
\partial_1 v \cdot \mathbf{N} \\
\partial_2 v \cdot \mathbf{N} \\
0
\end{pmatrix} \mathbf{N}
+ \frac{1}{\partial_z\varphi}\Big(\partial_z v\cdot \mathbf{N}\Big) \mathbf{N}, \ \text{on} \ \p S.
\end{equation*}
Now, let define $\Pi = I - \mathbf{n}\otimes\mathbf{n}$ (tangential projection). Then, from compatibility condition (\ref{1.15}), we have $\Pi (S^\varphi v \mathbf{N} ) = 0$. Therefore,
\begin{equation*} 
\Pi\partial_z v(\cdot,0) = \Pi\Big( \frac{\partial_z\varphi}{1+|\nabla h|^2}
\Big\{ 
\Big( \partial_1 h \partial_1 v + \partial_2 h \partial_2 v \Big)
- \begin{pmatrix}
\partial_1 v \cdot \mathbf{N} \\
\partial_2 v \cdot \mathbf{N} \\
0
\end{pmatrix} \mathbf{N}
- \frac{1}{\partial_z\varphi}\Big(\partial_z v\cdot \mathbf{N}\Big) \mathbf{N}
\Big\} \Big), \ \text{on} \ \p S.
\end{equation*}
Using Proposition \ref{commutator est}, \ref{eta regularity}, \ref{fraction est}, and (\ref{normal part est}), we get tangential part estimate,
\begin{equation} \label{tangential part est}
\begin{split}
|\Pi\partial_z v(\cdot,0)|_{X^{m,0}} \lesssim \Lambda\Big(\frac{1}{c_0}, \left|h\right|_{Y^{[\frac{m}{2}],1}} + \left\|v\right\|_{Y^{[\frac{m}{2}],1}}\Big)
\Big(\left|h\right|_{X^{m,1}} + \left|v(\cdot,0)\right|_{X^{m,1}}\Big),\quad\text{on}\quad \p S.
\end{split}
\end{equation}
Combining (\ref{normal part est}) and (\ref{tangential part est}), we get the result.
\end{proof}

\noindent Now we return to the Stress-continuity condition (\ref{1.15}). Applying $Z^m$,
\begin{equation} \label{335}
\begin{split}
&\left\{ Z^{m}q^b - gZ^{m}h - 2\varepsilon \Big(S^\varphi(Z^{m} v) \Big)^b - \varepsilon\Big(\Theta^{m}(v)\Big)^b \right\} \mathbf{N}   \\
& + \Big( q^b - gh - 2\varepsilon(S^\varphi v)^b + \nabla\cdot \frac{\nabla h}{\sqrt{1+|\nabla h|^2}} \Big)Z^{m} \mathbf{N} 
+ \Big(\nabla\cdot Z^{m} \frac{\nabla h}{\sqrt{1+|\nabla h|^2}} \Big)\mathbf{N}
+C^m (B)=0,
\end{split}
\end{equation}
where 
\begin{eqnarray*}
C^{m} (B) &:=& -C^{m} (B)_1 + C^{m} (B)_2   \\
&=& -2\varepsilon[Z^{m},(S^\varphi v)^b,\mathbf{N}] + \Big[Z^{m},q^b-gh+\nabla\cdot \frac{\nabla h}{\sqrt{1+|\nabla h|^2}},\mathbf{N}\Big].
\end{eqnarray*}
Using Proposition \ref{commutator est}, \ref{aniso embed}, and Lemma \ref{lem 3.2},
\begin{equation} \label{336}
\begin{split}
\left\| C^{m} (B)_1 \right\| &= 2\varepsilon \| [Z^{m},(S^\varphi v)^b,\mathbf{N}]  \|  \\
&\lesssim 2\varepsilon \ \Lambda\Big(\frac{1}{c_0},\left|h\right|_{Y^{[\frac{m}{2}],1}} + \left\|\nabla v\right\|_{Y^{[\frac{m}{2}],0}} \Big)
( \left|h\right|_{X^{m-1,1}} + \left|\nabla v(\cdot,0)\right|_{X^{m-1,0}} )  \\
&\lesssim 2\varepsilon \ \Lambda\Big(\frac{1}{c_0},\left|h\right|_{Y^{[\frac{m}{2}],1}} + \left\|\nabla v\right\|_{Y^{[\frac{m}{2}],0}}\Big)
\Big(\left|h\right|_{X^{m-1,1}} + |v^b|_{X^{m-1,1}}\Big) ,  \\
\left\| C^{m} (B)_2 \right\| &= 2\varepsilon \| [Z^{m},(S^\varphi v)^b \mathbf{n}\cdot\mathbf{n},\mathbf{N}] \| \\
&\lesssim 2\varepsilon \ \Lambda\Big(\frac{1}{c_0},\left|h\right|_{Y^{[\frac{m}{2}],1}} + \left\|\nabla v\right\|_{Y^{[\frac{m}{2}],0}}\Big)\Big(\left|h\right|_{X^{m-1,1}} + |v^b|_{X^{m-1,1}}\Big) ,
\end{split}
\end{equation}
where the way of estimate for $\left\| C^{m} (B)_2 \right\|$ is same as $\left\| C^{m} (B)_1 \right\|$. 
We also should estimate $Z^{m} \frac{\nabla h}{\sqrt{1+|\nabla h|^2}}$. We define $C^{m}(S)$ as  
\begin{equation} \label{337}
Z^{m} \frac{\nabla h}{\sqrt{1+|\nabla h|^2}} := 
\frac{\nabla Z^{m} h}{\sqrt{1+|\nabla h|^2}}
- \frac{\nabla h\langle\nabla h,\nabla Z^{m} h>}{\sqrt{1+|\nabla h|^2}^3} + C^{m} (S),
\end{equation}
 which is consist of low order polynomials in terms of $h$. Take a term in this $ C^{m}(S)$, then we take $L^2$ norm for the highest order,
and $L^\infty$ to others. For large $m \geq 2$, $L^\infty$ can be controlled by the highest order terms by Sobolev embedding. Therefore,
\begin{equation} \label{338}
\|\nabla\cdot C^{m}(S)\mathbf{N}\| \lesssim \Lambda\Big(\frac{1}{c_0},\left|h\right|_{Y^{[\frac{m}{2}],2}} + \left\|\nabla v\right\|_{Y^{[\frac{m}{2}],1}}\Big) |h|_{X^{m-1,2}}.
\end{equation}

\section{Pressure Estimates}
From Proposition \ref{high order eq}, we should estimate pressure $q$. From (\ref{1.12}), (\ref{1.13}), and (\ref{1.15}), pressure $q(t,x,y,z)$ solves,
\begin{eqnarray*}
-\Delta^{\varphi} q &=& \nabla^{\varphi}\cdot(v\cdot\nabla^{\varphi}v)\quad\text{in}\quad S,\\
q^b &:=& q\vert_{z=0} = q(t,x,y,0) = gh + 2\varepsilon (S^{\varphi} v)^b \mathbf{n}\cdot\mathbf{n} - \nabla\cdot\frac{\nabla h}{\sqrt{1+|\nabla h|^2}}.
\end{eqnarray*}
We decompose pressure $q$ into $ q=q^E+q^{NS}+q^S $, where \\
$q^E$ solves
\begin{equation} \label{q_E}
\begin{split}
 - \Delta^\varphi q^E &= \nabla^{\varphi}\cdot\big( v\cdot\nabla^\varphi v \big)\quad\text{in}\quad S,  \\
 q^E\vert_{z=0} &= gh,
\end{split}
\end{equation} 
and $q^{NS}$ solves
\begin{equation} \label{q_NS}
\begin{split}
- \Delta^\varphi q^{NS} &= 0,\quad\text{in}\quad S,  \\
q^{NS}\vert_{z=0} &= 2\varepsilon(S^\varphi v)\mathbf{n}\cdot\mathbf{n},
\end{split}
\end{equation}
and $q^S$ solves
\begin{equation} \label{q_S}
\begin{split}
- \Delta^\varphi q^S &= 0,\quad\text{in}\quad S,  \\
q^S\vert_{z=0} &= -\nabla\cdot\frac{\nabla h}{\sqrt{1+|\nabla h|^2}}.
\end{split}
\end{equation}
To find expand $\Delta^{\varphi}$, we should calculate gradient form. Using Definition \ref{varphipartial},
\begin{equation} \label{gradient varphi}
\nabla^\varphi f
= \begin{pmatrix}
\partial_1 f - \frac{\partial_1 \varphi}{\partial_z \varphi}\partial_z f \\
\partial_2 f - \frac{\partial_2 \varphi}{\partial_z \varphi}\partial_z f \\
\frac{\partial_z f}{\partial_z \varphi}
\end{pmatrix}
= \begin{pmatrix}
1 & 0 & - \frac{\partial_1 \varphi}{\partial_z \varphi} \\
0 & 1 & - \frac{\partial_2 \varphi}{\partial_z \varphi} \\
0 & 0 & \frac{1}{\partial_z \varphi}
\end{pmatrix}
\begin{pmatrix}
\partial_1 f \\
\partial_2 f \\
\partial_z f \\
\end{pmatrix}
:= \frac{1}{\partial_z \varphi} P^*\nabla f,
\end{equation}
where 
\begin{equation} \label{P form}	
P :=
\begin{pmatrix}
\partial_z \varphi & 0 & 0 \\
0 & \partial_z \varphi & 0 \\
-\partial_1 \varphi & -\partial_2 \varphi & 1
\end{pmatrix}.
\end{equation}
To calculate divergence, using Definition \ref{varphipartial},
\begin{equation} \label{divergent varphi}
\nabla^\varphi\cdot v = \frac{1}{\partial_z \varphi} \nabla\cdot(Pv).
\end{equation}
Therefore, by (\ref{gradient varphi}) and (\ref{divergent varphi}), we get the following form for $\Delta^{\varphi}$.
\begin{equation} \label{laplacian varphi}
\Delta^\varphi f
=\nabla^\varphi\cdot(\nabla^\varphi f)
=\frac{1}{\partial_z \varphi} \nabla\cdot(P\nabla^\varphi f)
=\frac{1}{\partial_z \varphi} \nabla\cdot(E\nabla f),
\end{equation}
where
\begin{equation} \label{E form}	
E := \frac{1}{\partial_z \varphi} PP^* =
\begin{pmatrix}
\partial_z \varphi & 0 & -\partial_1 \varphi \\
0 & \partial_z \varphi & -\partial_2 \varphi \\
-\partial_1 \varphi & -\partial_2 \varphi & \frac{1+(\partial_1 \varphi)^2+(\partial_2 \varphi)^2}{\partial_z \varphi}
\end{pmatrix}.
\end{equation}
We state two lemmas about elliptic equation with Dirichlet boundary condition. 

\begin{lemma} \label{lemma 4.1}
Let $m\geq 6$ and $\rho$ solves the following system,
\begin{equation*}
-\nabla\cdot (E\nabla\rho ) = \nabla\cdot F,\quad\text{in}\quad S,\quad \rho\vert_{z=0} := \rho(t,x,y,0) = 0,
\end{equation*}
where matrix $E$ and $P$ are defined in (\ref{E form}) and (\ref{P form}), respectively. Then we have the following estimate.
\begin{eqnarray*}
  \|\nabla\rho \|_{X^{m-1,0}} +  \|\nabla^2\rho \|_{X^{m-2,0}} &\lesssim&
\Lambda\Big(\frac{1}{c_0}, |h|_{Y^{[\frac{m-1}{2}],1}} + |h|_{2,\infty} + |h|_3  + \left\|\nabla\cdot F\right\|_{H_{tan}^{1}} + \left\|F\right\|_{H_{tan}^{2}} \Big)   \\
&\times& ( \|F\|_{X^{m-1,0}} + \|\nabla\cdot F \|_{X^{m-2,0}} + |h|_{X^{m-2,1}} ) ,
\end{eqnarray*}
where function spaces $X^{m,s}$, $Y^{m,s}$, and $H^{s}_{tan}$ are defined in Definition \ref{space def} and \ref{XY space}, respectively.
\end{lemma} 
\begin{proof}
By variational argument with homogeneous boundary condition, we have
\begin{equation} \label{zero order}
\|\nabla\rho\| \leq \Lambda\Big( \frac{1}{c_0}, |h|_{1,\infty} \Big) \|F\|_{L^2}.
\end{equation}
If we apply $Z^m$ to the equation, then divergence structure is broken, since $Z_3$ and $\partial_z$ do not commute. Therefore, we redefine,
\begin{equation} \label{redef Z_3}
\tilde{Z}_3 f := Z_3 f + (2z+1)f.
\end{equation}
With this new definition, it is easy to check $\tilde{Z}_3\partial_z = \partial_z Z_3$.
We apply 
\[
\tilde{Z}^{m-1} = \p_t^{k}Z_1^{\alpha_1}Z_2^{\alpha_2}\tilde{Z}_3^{\alpha_3}, (k + \alpha_1 + \alpha_2 + \alpha_3 = m-1)
\]
to the equation.
\begin{equation} \label{commu form}
\nabla\cdot(Z^{m-1} (E\nabla\rho)) = 
\nabla\cdot (Z^{m-1} F + (\tilde{Z}^{m-1} - Z^{m-1})F_h - (\tilde{Z}^{m-1} - Z^{m-1})(E\nabla\rho)_h ),
\end{equation}
where sub index $h$ means horizontal component, {i.e.} $ F_h := (F_1,F_2,0) $. Again, considering commutators between $Z^{m-1}$ and $E$,
\begin{equation*}
\begin{split}
\nabla\cdot(E\cdot\nabla(Z^{m-1}\rho)) &= 
\nabla\cdot \Big(Z^{m-1} F + (\tilde{Z}^{m-1}- Z^{m-1})F_h - (\tilde{Z}^{m-1}-Z^{m-1})(E\nabla\rho)_h \Big) \\
&\quad + \nabla\cdot \tilde{C},
\end{split}
\end{equation*}
where 
$$
\tilde{C} := -E[Z^{m-1},\nabla]\rho - \Big( \sum_{ \tiny \begin{array}{c}
	k_{1}+\beta+k_{2}+\gamma=m-1, \\ k_{1}+\beta \neq 0
	\end{array}} c_{k_1,\beta,k_2,\gamma} \ \p_t^{k_{1}}Z^\beta E \ \p_t^{k_{2}} Z^\gamma\rho \Big).
$$
Since $Z^{m-1} \rho$ is zero on the boundary, using (\ref{zero order}),
\begin{equation} \label{high est}
\begin{split}
\left\| \nabla\rho \right\|_{X^{m-1,0}} &\leq
\Lambda\Big( \frac{1}{c_0}, |h|_{1,\infty} \Big) \Big(\tiny \|F\|_{X^{m-1,0}} + \left\|E\nabla\rho\right\|_{X^{m-2,0}} + \left\|E[Z^{m-1},\nabla]\rho\right\| \\
& \quad + \Big\|  \sum_{\tiny \begin{array}{c}
	k_{1}+\beta+k_{2}+\gamma=m-1, \\ k_{1}+\beta \neq 0
	\end{array}} c_{k_1,\beta,k_2,\gamma} \ \p_t^{k_{1}}Z^\beta E \ \p_t^{k_{2}} Z^\gamma\rho \Big\| \Big)  . 
\end{split}
\end{equation}
Three terms on the right hand side can be estimated using Proposition \ref{commutator est}.
\begin{equation} \label{three RHS}
\begin{split}
&\|E\nabla\rho\|_{X^{m-2,0}} \lesssim
\Lambda\Big(\frac{1}{c_0},\left|h\right|_{Y^{[\frac{m-1}{2}],0}} + \left\|\nabla\rho\right\|_{Y^{[\frac{m-1}{2}],0}} \Big) \Big( \left\|\nabla\rho\right\|_{X^{m-2,0}} + \left|h\right|_{X^{m-2,1}} \Big) ,  \\
&\left\|E[Z^{m-1},\nabla]\rho\right\| \lesssim \Big\|E \sum_{|\beta|\leq m-2}c_{\alpha,\beta}\partial_z\Big(Z^\beta\rho\Big) \Big\| \lesssim
\Lambda\Big(\frac{1}{c_0},|h|_{Y^{[\frac{m-1}{2}],0}} \Big) \left\|\nabla\rho \right\|_{X^{m-2,0}} , \\
&\Big\| \sum_{\tiny \begin{array}{c}
	k_{1}+\beta+k_{2}+\gamma=m-1, \\ k_{1}+\beta \neq 0
	\end{array}} c_{k_1,\beta,k_2,\gamma} \ \p_t^{k_{1}}Z^\beta E \ \p_t^{k_{2}} Z^\gamma\rho \Big\| \lesssim \Lambda\Big( |h|_{Y^{[\frac{m-2}{2}],1}} + \left\|\nabla\rho\right\|_{Y^{[\frac{m-1}{2}],0}} \Big) \\
&\quad\quad\quad\quad\quad\quad\quad\quad\quad\quad\quad\quad\quad\quad\quad\quad\quad\quad\quad\quad\quad\quad  \times \Big( \left\|\nabla\rho\right\|_{X^{k-1,0}} + |h|_{X^{k-1,1}} \Big).
\end{split}
\end{equation}
Putting (\ref{three RHS}) to the right hand side of (\ref{high est}), we have

\begin{equation}
\left\|\nabla\rho\right\|_{X^{m-1,0}} \lesssim
\Lambda\Big(\frac{1}{c_0},\left|h\right|_{Y^{[\frac{m-1}{2}],1}} + \left\|\nabla\rho\right\|_{Y^{[\frac{m-1}{2}],0}} \Big) 
\Big( \left\|F\right\|_{X^{m-1,0}} + \left\|\nabla\rho\right\|_{X^{m-2,0}} + \left|h\right|_{X^{m-2,1}} \Big).
\end{equation}
On the right hand side, we use induction for $ \left\|\nabla\rho \right\|_{X^{m-2,0}} $ until zero order term $\|\nabla\rho\|$ and $\|\nabla\rho\|_{L^\infty}$ appear. To estimate $\|\nabla^2 \rho\|_{X^{m-2,0}}$, we are suffice to focus on $\|\p_{zz} \rho\|_{X^{m-2,0}}$. We apply $Z^{m-1}\p_z$ to the equation to the system again, and make the equation with the form of (\ref{commu form}). We follow the same argument as above to get the result.
\begin{eqnarray*}
	\left\|\nabla^2\rho\right\|_{X^{m-2,0}} &\lesssim&
	\Lambda\Big(\frac{1}{c_0}, |h|_{Y^{[\frac{m-1}{2}],1}} + |h|_{2,\infty} + |h|_3  + \left\|\nabla\cdot F\right\|_{H_{tan}^{1}} + \left\|F\right\|_{H_{tan}^{2}} \Big)    \\
	&\times& ( \left\|F\right\|_{X^{m-1,0}} + \left\|\nabla\cdot F\right\|_{X^{m-2,0}} + |h|_{X^{m-2,1}} ) .
\end{eqnarray*}

\end{proof}

\noindent For homogeneous elliptic problem, \textit{i.e} source $F=0$, we can derive standard Sobolev regularity.
\begin{lemma} \label{lemma 4.2}
Let $m\geq 6$ and $\rho$ solves the following system,
\begin{equation*}
-\nabla\cdot (E\nabla\rho ) = 0,\quad\text{in}\quad S,\quad \rho\vert_{z=0} := \rho(t,x,y,0) = f^b,
\end{equation*}
Then we have the following estimate.
\begin{equation*}
\|\nabla\rho\|_{\mathcal{H}^{m-1,0}} \lesssim \Lambda\Big(\frac{1}{c_0}, |h|_{Y^{[\frac{m-1}{2}],1}} + |h|_{2,\infty} + |h|_{3} + |f^b|_{Y^{[\frac{m-1}{2}],0}} \Big) \Big( |h|_{X^{m-1,\frac{1}{2}}} +  |f^b |_{X^{m-1,\frac{1}{2}}} \Big),
\end{equation*}
where function space $\mathcal{H}^{m,s}$ is defined in Definition \ref{XY space}.
\end{lemma} 
\begin{proof}
We decompose $\rho = \rho^H + \rho^r$, where $\rho^H$ absorbs boundary data and $\rho^r$ solves
\begin{equation*}
-\nabla\cdot (E\nabla\rho^r) = \nabla\cdot (E\nabla\rho^H),\quad \rho^r\vert_{z=0} = 0.
\end{equation*}
We choose $\rho^H$ as $\hat{\rho}^H(\xi,z) = \chi(z\xi)\hat{f}^b$, where $\hat{\rho}$ is Fourier transform with respect to horizontal variables $x,y$. $\xi$ is corresponding two dimensional frequency variable. Using Proposition \ref{eta regularity}, we get
\begin{equation} \label{rho H}
\| \nabla\rho^H \|_{\mathcal{H}^{m-1,0}} \lesssim | f^b |_{X^{m-1,\frac{1}{2}}},\quad \text{and}\quad \|\rho^H \|_{\mathcal{K}^{m-1,0}} \lesssim | f^b |_{Y^{m-1,0}},
\end{equation}
where $\mathcal{H}^{m,s}$ and $\mathcal{K}^{m,s}$ are defined in Definition \ref{XY space}. \\
\noindent We use similar approach as we did in Lemma \ref{lemma 4.1} for standard Sobolev regularity and we get
\begin{eqnarray*}	
\|\nabla\rho^r \|_{\mathcal{H}^{m-1,0}} &\lesssim& 
\Lambda\Big(\frac{1}{c_0}, |h|_{Y^{[\frac{m-1}{2}],1}} + |h|_{2,\infty} + |h|_3  + \|\nabla\cdot ( E \nabla\rho^H ) \|_{H_{tan}^1} + \|E \nabla\rho^H \|_{H_{tan}^2} \Big) \\
&&\times ( \|E \nabla\rho^H \|_{\mathcal{H}^{m-1,0}} + |h|_{X^{m-2,1}} )    \\
&\lesssim& \Lambda\Big(\frac{1}{c_0}, |h|_{Y^{[\frac{m-1}{2}],1}} + |h|_{2,\infty} + |h|_3  + \|E \nabla\rho^H \|_{H^2} \Big) ( \|E \nabla\rho^H \|_{\mathcal{H}^{m-1,0}} + |h|_{X^{k-1,1}} ),
\end{eqnarray*}
where $\|E \nabla\rho^H \|_{\mathcal{H}^{m-1,0}}$ on the right hand side can be controlled using (\ref{rho H}) and Proposition \ref{commutator est},
\begin{equation*}
\|E \nabla\rho^H \|_{\mathcal{H}^{m-1,0}} \lesssim
\Lambda\Big(\frac{1}{c_0}, |h|_{Y^{[\frac{m-1}{2}],0}} + |f^b|_{Y^{[\frac{m-1}{2}],0}} \Big) ( |h|_{X^{m-1,\frac{1}{2}}} + |f^b |_{X^{m-1,\frac{1}{2}}} ).
\end{equation*}
Consequently, this yields estimate for $\nabla\rho^{r}$,
\begin{equation} \label{rho r}
\|\nabla\rho^r\|_{\mathcal{H}^{m-1,0}} \lesssim \Lambda\Big(\frac{1}{c_0}, |h|_{Y^{[\frac{m-1}{2}],1}} + |h|_{2,\infty} + |h|_{3} + |f^b|_{Y^{[\frac{m-1}{2}],0}} \Big) ( |h|_{X^{m-1,\frac{1}{2}}} + |f^b|_{X^{m-1,\frac{1}{2}}} ).
\end{equation}
We put (\ref{rho H}) and (\ref{rho r}) together.  \\
\end{proof}

We use above Lemma \ref{lemma 4.1} and \ref{lemma 4.2} to estimate $q^{E}$, $q^{NS}$ and $q^S$, those are defined in (\ref{q_E}), (\ref{q_NS}), and (\ref{q_S}).
\begin{proposition} \label{q_E est}
	Let $m\geq 6$ and $q^{E}$ solves (\ref{q_E}), we have the following estimate.
	\begin{eqnarray*}
		 && \|\nabla q^{E}\|_{X^{m-1,0}} + \|\nabla^2 q^{E}\|_{X^{m-2,0}} \\
		&\lesssim& \Lambda\Big(\frac{1}{c_0}, |h|_{Y^{[\frac{m-1}{2}],1}} + |h|_{2,\infty} + |h|_3  + \|\nabla v\|_{2} + \|v\|_3 \Big) ( \|v\|_{X^{m-1,0}} + \|\nabla v\|_{X^{m-2,0}} + |h|_{X^{m-1,1}}  ).
	\end{eqnarray*}
\end{proposition}
\begin{proof}
	We split $q^E = q^E_{1} + q^E_{2}$, where
	\begin{eqnarray*}
		q^E_1 &=& 0,\quad q^E_1 \vert_{z=0} = gh  \\
		q^E_2 &=& \nabla\cdot(P(v\cdot\nabla^{\varphi} v)) = \p_z\varphi \nabla^{\varphi}v\cdot \nabla^{\varphi} v,\quad q^E_2 \vert_{z=0} = 0.  \\
	\end{eqnarray*}
	Using Lemma \ref{lemma 4.2}
	\begin{eqnarray*}
		\|\nabla q^E_1\|_{\mathcal{H}^{m-1,0}} &\lesssim& \Lambda\Big(\frac{1}{c_0}, |h|_{Y^{[\frac{m-1}{2}],1}} + |h|_{2,\infty} + |h|_{3} \Big) ( |h|_{X^{m-1,\frac{1}{2}}} ).  \\
	\end{eqnarray*}
	For $q^E_2$, using Lemma \ref{lemma 4.1} and $\nabla^\varphi \cdot (v\cdot\nabla^\varphi v) = \nabla^\varphi v : (\nabla^\varphi v)^T $, which comes from divergence free condition,
	\begin{eqnarray*}
		&&\|\nabla q^{E}_2\|_{X^{m-1,0}} + \|\nabla^2 q^{E}_2\|_{X^{m-2,0}} \\ &\lesssim& \Lambda\Big(\frac{1}{c_0}, |h|_{Y^{[\frac{m-1}{2}],1}} + |h|_{2,\infty} + |h|_3  + \|\nabla v \|_{H^1_{tan}} + \|P(v\cdot\nabla^{\varphi} v)\|_{H^2_{tan}} \Big)   \\
		&& \times ( \|\nabla\cdot P(v\cdot\nabla^{\varphi} v)\|_{X^{m-1,0}} + \|P(v\cdot\nabla^{\varphi} v)\|_{X^{m-1,0}} + |h|_{X^{m-2,1}}  )   \\
		&\lesssim& \Lambda\Big(\frac{1}{c_0}, |h|_{Y^{[\frac{m-1}{2}],1}} + |h|_{2,\infty} + |h|_3  + \|\nabla v\|_{2} + \|v\|_3 \Big) ( \|v\|_{X^{m-1,0}} + \|\nabla v\|_{X^{m-2,0}} + |h|_{X^{m-2,1}}  ) .  \\
	\end{eqnarray*}
	This yields estimate for $\|\nabla q^{E}\|_{X^{m-1,0}}$. \\
\end{proof}
\begin{proposition} \label{q_NS est}
Let $m\geq 6$ and $q^{NS}$ solves (\ref{q_NS}), we have the following estimate.
\begin{eqnarray*}
\|\nabla q^{NS}\|_{\mathcal{H}^{m-1,0}} &\lesssim& \Lambda\Big(\frac{1}{c_0}, |h|_{Y^{[\frac{m-1}{2}],1}} + |h|_{2,\infty} + |h|_{3} + \| \nabla v \|_{Y^{[\frac{m-1}{2}],0}} \Big) \\
&& \times ( |h|_{X^{m-1,\frac{1}{2}}} + \varepsilon |h|_{X^{m-1,\frac{3}{2}}} + \varepsilon \|\nabla v \|_{X^{m-1,1}} + \varepsilon \|v \|_{X^{m-1,2}} ) .
\end{eqnarray*}
\end{proposition}
\begin{proof}
Applying Lemma \ref{lemma 4.2} with $ f^b = 2\varepsilon(S^\varphi v)^b \mathbf{n}\cdot \mathbf{n} $,
\begin{eqnarray*}
\|\nabla q^{NS}\|_{\mathcal{H}^{m-1,0}} &\lesssim& \Lambda\Big(\frac{1}{c_0}, |h|_{Y^{[\frac{m-1}{2}],1}} + |h|_{2,\infty} + |h|_{3} + \Big|\varepsilon (S^\varphi v)^b \mathbf{n}\cdot \mathbf{n} \Big|_{Y^{[\frac{m-1}{2}],0}} \Big) \\
&&\times ( \left|h\right|_{X^{m-1,\frac{1}{2}}} + \Big|\varepsilon (S^\varphi v)^b \mathbf{n}\cdot \mathbf{n} \Big|_{X^{m-1,\frac{1}{2}}} )   \\
&\lesssim& \Lambda\Big(\frac{1}{c_0}, |h|_{Y^{[\frac{m-1}{2}],1}} + |h|_{2,\infty} + |h|_{3} + \| \nabla v \|_{Y^{[\frac{m-1}{2}],0}} \Big) \\
&&\times ( \left|h\right|_{X^{m-1,\frac{1}{2}}} + \varepsilon\left|h\right|_{X^{m-1,\frac{3}{2}}} + \varepsilon |(\nabla v)^b |_{X^{m-1,\frac{1}{2}}} )   \\
&\lesssim& \Lambda\Big(\frac{1}{c_0}, |h|_{Y^{[\frac{m-1}{2}],1}} + |h|_{2,\infty} + |h|_{3} + \left\| \nabla v \right\|_{Y^{[\frac{m-1}{2}],0}} \Big) \\
&&\times ( \left|h\right|_{X^{m-1,\frac{1}{2}}} + \varepsilon\left|h\right|_{X^{m-1,\frac{3}{2}}} + \varepsilon\left\|\nabla v \right\|_{X^{m-1,1}} + \varepsilon\left\|v \right\|_{X^{m-1,2}} ) ,
\end{eqnarray*}
where we used Proposition \ref{commutator est} and \ref{aniso embed} in the first step and used Lemma \ref{lem 3.2} in the second step. \\
\end{proof}

\begin{proposition} \label{q_S est}
Let $m\geq 6$ and $q^{S}$ solves (\ref{q_S}), we have the following estimate.
\begin{equation*}
 \|\nabla q^{S}\|_{\mathcal{H}^{m-1,0}} \lesssim \Lambda\Big(\frac{1}{c_0}, |h|_{2,\infty} + |h|_{3} + |h|_{Y^{[\frac{m-1}{2}],2}} \Big) \Big( |h|_{X^{m-1,\frac{1}{2}}} + \Big| \frac{\nabla h}{\sqrt{1+|\nabla h|^2}} \Big|_{X^{m-1,\frac{3}{2}}} \Big) .
\end{equation*}
\end{proposition}
\begin{proof}
Applying Lemma \ref{lemma 4.2}, $ f^b = -\nabla\cdot\frac{\nabla h}{\sqrt{1+|\nabla h|^2}} $,
\begin{equation*}
\begin{split}
\|\nabla q^{S}\|_{\mathcal{H}^{m-1,0}} &\lesssim \Lambda\Big(\frac{1}{c_0}, |h|_{Y^{[\frac{m-1}{2}],1}} + |h|_{2,\infty} + |h|_{3} + \Big| \nabla\cdot\frac{\nabla h}{\sqrt{1+|\nabla h|^2}} \Big|_{Y^{[\frac{m-1}{2}],0}} \Big) \\
&\quad\times \Big( |h|_{X^{m-1,\frac{1}{2}}} + \Big| \frac{\nabla h}{\sqrt{1+|\nabla h|^2}} \Big|_{X^{m-1,\frac{3}{2}}} \Big)  \\
& \lesssim \Lambda\Big(\frac{1}{c_0}, |h|_{2,\infty} + |h|_{3} + |h|_{Y^{[\frac{m-1}{2}],2}} \Big) \Big( |h|_{X^{m-1,\frac{1}{2}}} + \Big| \frac{\nabla h}{\sqrt{1+|\nabla h|^2}} \Big|_{X^{m-1,\frac{3}{2}}} \Big) .
\end{split}
\end{equation*}
\end{proof}

We should also estimate $L^\infty$-type for pressure. Since we have standard Sobolev regularity for $q^{NS}$ and $q^S$, we can use standard Sobolev embedding. For $q^E$, we use anisotropic embedding Proposition \ref{aniso embed}. 
\begin{proposition} \label{q_infty est}
Let $m\geq 6$. We have the following $L^\infty$ type estimate for $q:=q^{E}+q^{NS}+q^{S}$.
\begin{equation}
\left\|\nabla q\right\|_{Y^{[\frac{m}{2}],0}} \lesssim \Lambda\Big(\frac{1}{c_0}, |h|_{Y^{[\frac{m}{2}],1}} + |h|_{2,\infty} + |h|_3  + \|\nabla v\|_{2} + \|v\|_3 \Big) ( \|\nabla v\|_{X^{[\frac{m}{2}]-1,0}} + |h|_{X^{[\frac{m}{2}]-1,1}}  )  .
\end{equation}
\end{proposition} 
\begin{proof}
Using anisotropic Sobolev embedding and Proposition \ref{q_E est},
\begin{equation} \label{nabla q_E}
\begin{split}
\|\nabla q^E\|^2_{Y^{[\frac{m}{2}],0}} &\lesssim \|\nabla^2 q^E\|^2_{X^{[\frac{m}{2}]+1,0}}  + \|\nabla q^E \|^2_{X^{[\frac{m}{2}]+2,0}}  \\
&\lesssim  \Lambda\Big(\frac{1}{c_0}, |h|_{Y^{[\frac{[\frac{m}{2}]-1}{2}],1}} + |h|_{2,\infty} + |h|_3  + \|\nabla v\|_{2} + \|v\|_3 \Big)  \\
&\quad\times ( \|\nabla v\|_{X^{[\frac{m}{2}]-1,0}} + |h|_{X^{[\frac{m}{2}]-1,1}}  ).
\end{split}
\end{equation}
For $q^{NS}$ and $q^S$, we can use standard Sobolev embedding from the regularity of Proposition \ref{q_NS est} and \ref{q_S est}.
\begin{eqnarray*}
\left\|\nabla q\right\|_{Y^{[\frac{m}{2}],0}} &\lesssim& \|\nabla q^E \|_{Y^{[\frac{m}{2}],0}} +  \|\nabla q^{NS} \|_{Y^{[\frac{m}{2}],0}} +  \|\nabla q^S \|_{Y^{[\frac{m}{2}],0}}  \\
&\lesssim& \|\nabla^2 q^E\|^2_{X^{[\frac{m}{2}]+1,0}}  + \|\nabla q^E \|^2_{X^{[\frac{m}{2}]+2,0}} + \|\nabla q^{NS}\|_{X^{[\frac{m}{2}],2}} + \|\nabla q^{S}\|_{X^{[\frac{m}{2}],2}} .
\end{eqnarray*}
By (\ref{nabla q_E}) and Proposition \ref{q_NS est} and \ref{q_S est}, we end the proof with $m\geq 6$.
\end{proof}

\section{Energy Estimates for $Z^m \neq \p_t^m$}
In this section, we perform energy estimate of the system derived in Proposition \ref{high order eq} in $S$. As defined in (\ref{volume integral}), We will use notation $\dd V_t$ for volume element in $S$,
\begin{equation}
\int_{\O_t} f \dd V := \int_{S}f\circ\Phi dV_t,\quad \text{where}\quad dV_t = \partial_z\varphi(t,x,y,z)\dd x \dd y \dd z.
\end{equation}
The following lemma gives integration by parts rule in fixed domain $S$, with our new derivatives in Definition \ref{varphipartial} and volume element $\dd V_t$ 
\begin{lemma} \label{IBP}
Let $f$ and $g$ are functions defined on $S$. Then, we have the following integration by parts rules in $S$.
\begin{eqnarray*}
\int_{S} \partial_i^\varphi fg dV_t &=& -\int_{S} f \partial_i^\varphi g dV_t + \int_{z=0} fg N_i dy, \,\,\,\,\,\,\,\,i=1,2,3,  \\
\int_{S} \partial_t^\varphi fg dV_t &=& \partial_t\int_{S} fg dV_t - \int_{S} f\partial_t^\varphi g dV_t - \int_{z=0}fg\partial_t h,
\end{eqnarray*}
where $\mathbf{N} = (N_1,N_2,N_3)$.
\end{lemma}
\begin{proof}
This can be derived directly from standard integration by parts in $\O_t$. See \cite{NMFR} for more detail.
\end{proof}

\begin{corollary} \label{IBP coro}
Let $v(t,\cdot)$ is a vector field on $S$, such that $\nabla^\varphi\cdot v = 0 $. Then, for every smooth functions $f,g$ and smooth vector field $u,w$, we have the following estimates.
\begin{eqnarray*}
\int_{S} (\partial_t^\varphi f + v\cdot\nabla^\varphi f )f dV_t &=& \frac{1}{2}\partial_t\int_{S}|f|^2 dV_t - \frac{1}{2}\int_{z=0} (\partial_t h-v\cdot \mathbf{N} ) dy,   \\
\int_{S} (\Delta^\varphi f )g dV_t &=& -\int_{S}\nabla^\varphi f\cdot\nabla^\varphi g dV_t + \int_{z=0}\nabla^\varphi f\cdot \mathbf{N} g dy,  \\
\int_{S}\nabla^\varphi\cdot (S^\varphi u )\cdot w dV_t &=& -\int_{S}S^\varphi u \cdot S^\varphi w dV_t + \int_{z=0} (S^\varphi u \mathbf{N} )\cdot w dy.
\end{eqnarray*}
\end{corollary}
\begin{proof}
This comes from Lemma \ref{IBP coro}. See also \cite{NMFR}.
\end{proof}

\begin{lemma} \label{L2 est}
Let $v,h$, and $\varphi$ are smooth solutions of the system (\ref{1.12})-(\ref{1.15}). Then we have the following energy identity.
\begin{equation}
\frac{d}{dt}\Big( \int_S\left|v\right|^2dV_t + g\int_{z=0}\left|h\right|^2dy + 2\int_{\partial S} ( \sqrt{1+|\nabla h|^2} - 1 ) dy \Big)
+ 4\varepsilon\int_S\left|S^\varphi v\right|^2dV_t = 0 .
\end{equation}
\end{lemma}
\begin{proof}
Using (\ref{1.12}), (\ref{1.13}), (\ref{1.14}), and Corollary \ref{IBP coro}, we have,
\begin{eqnarray*}
\frac{d}{dt}\int_S |v|^2 dV_t = 2\int_S \nabla^\varphi\cdot(2\varepsilon S^\varphi v - q)\cdot v dV_t
\end{eqnarray*}
We use the last equality in Corollary \ref{IBP coro} and (\ref{1.15}) to get,
\begin{eqnarray*}
\frac{d}{dt}\int_S |v|^2 dV_t + 4\varepsilon\int_S |S^\varphi v|^2 dV_t &=& 2\int_{\partial S} (2\varepsilon S^\varphi v - qI)\mathbf{N}\cdot v dy  \\
&=& 2\int_{\partial S} \Big( -gh\mathbf{N}\cdot v + \nabla\cdot\frac{\nabla h}{\sqrt{1+|\nabla h|^2}}\mathbf{N}\cdot v \Big) dy  \\
&=& -g\frac{d}{dt}\int_{\partial S} |h|^2 dy - 2\int_{\partial S}\frac{\nabla h}{\sqrt{1+|\nabla h|^2}} \nabla h_t dy  \\
&=& -g\frac{d}{dt}\int_{\partial S} |h|^2 dy - 2\frac{d}{dt}\int_{\partial S} ( \sqrt{1+|\nabla h|^2} - 1 ) dy .
\end{eqnarray*}
\end{proof}

Before we perform energy estimate, we use define the $\Lambda_{m,\infty}$ which contain $L^{\infty}$ type norms and $L^{2}$ type norms with finite order.
\begin{equation} \label{Lambda infty}
	\Lambda_{m,\infty}(t) := \Lambda \Big(\frac{1}{c_0}, |h(t)|_{Y^{[\frac{m}{2}],1}} + \|\nabla v(t) \|_{ Y^{[\frac{m}{2}],0} } + |h(t)|_{3} + \|v(t)\|_{3} + \|\nabla v(t)\|_{2}  \Big),
\end{equation}

The following proposition gives high order energy estimates for $Z^{m} \neq \p_t^{m}$ in the system of Proposition \ref{high order eq}.
\begin{proposition} \label{energy estimate}
Assume that 
Let $m\geq 6$ and $(v,\varphi,q,h)$ is a smooth solution on $[0,T]$ of the system (\ref{1.12})-(\ref{1.15}). Also assume that
\begin{equation}
\partial_z\varphi \geq c_0,\quad |h|_{2,\infty} \leq \frac{1}{c_0},\quad \forall t\in [0,T].
\end{equation}
Then, from the higher order system Proposition \ref{high order eq}, we get the following high order energy estimate.
\begin{equation*}
\begin{split}
& \|v(T)\|^2_{X^{m-1,1}} + \|h(T)\|^2_{X^{m-1,2}} + \varepsilon\int_0^T \|\nabla v(t)\|^2_{X^{m-1,1}} dt \\
& \leq \Lambda(\frac{1}{c_0}) \sum_{Z^m\neq\p_t^m} ( \|(Z^m v)(0)\|_{L^2} + |(Z^m h)(0)|_{H^1} ) \\
&\quad + \int_0^T \Lambda_{m,\infty}(t) \Big( \|v(t)\|^2_{X^{m,0}} + \|\nabla v(t)\|^2_{X^{m-1,0}} + \varepsilon \|\nabla v(t) \|^2_{X^{m,0}} + \varepsilon \|v(t)\|_{X^{m-1,2}} + |h(t)|^2_{X^{m-1,\frac{5}{2}}} \Big) dt.  \\
\end{split}
\end{equation*}

\end{proposition}
\begin{remark}
	In this estimate, we do not treat $Z^{m} = \p_t^{m}$ case. Therefore, function space of $h$ is $X^{m-1,\frac{5}{2}}(S)$, instead of $X^{m,\frac{3}{2}}(S)$.
\end{remark}
\begin{proof}
Let $Z^{m}\neq \partial_t^m$. From the system in Proposition \ref{high order eq} and Corollary \ref{IBP coro}, we get,
\begin{equation} \label{basic high}
\begin{split}
& \frac{d}{dt}\int_S\left|Z^{m} v\right|^2dV_t
+ 4\varepsilon\int_S\left|S^\varphi (Z^{m} v)\right|^2dV_t  \\
&=2\int_{\p S}(2\varepsilon S^\varphi(Z^{m} v)^b-(Z^{m} q)^bI)\mathbf{N}\cdot (Z^{m} v)^b \dd x \dd y
+ R_S + R_C,
\end{split}
\end{equation}
where
\begin{equation} \label{RS}
R_S := 2\varepsilon \int_S\left\{\nabla^\varphi\cdot\Theta^{m}(v)-D^{m}(S^\varphi v)\right\}\cdot Z^{m} v dV_t ,
\end{equation}
\begin{equation} \label{RC}
R_C := 2\int_S \{C^{m} (q)-C^{m} (T)\}(Z^{m} v) + C^{m}(d)(Z^{m} q) dV_t .
\end{equation}
And we use the third equation in Proposition \ref{high order eq} (the continuity of stress tensor condition) for the boundary integral term on the right hand side. From the boundary condition, 
\begin{equation} \label{BDRY integral}
\begin{split}
& 2\int_{\p S}(2\varepsilon S^\varphi(Z^{m} v)^b-(Z^{m} q)^bI)\mathbf{N}\cdot (Z^{m} v)^b \\
=& 2\int_{\p S}\left\{-gZ^{m} h+ \Big(\nabla\cdot\Big(\frac{\nabla Z^{m} h}{\sqrt{1+|\nabla h|^2}}
- \frac{\nabla h \langle\nabla h,\nabla Z^{m} h\rangle}{\sqrt{1+|\nabla h|^2}^3} 
+ C^{m} (S)\Big) \Big)\right\}\mathbf{N}\cdot(Z^{m} v)^b   \\
& + 2\int_{\p S} \Big\{ (q-gh)I-2\varepsilon(S^\varphi v)^b+ \nabla\cdot\frac{\nabla h}{\sqrt{1+|\nabla h|^2}}I \Big\} 
(Z^{m} \mathbf{N})\cdot(Z^{m} v)^b + R_B,  \\ 
=& 2\int_{\p S}\left\{-gZ^{m} h+ \Big(\nabla\cdot\Big(\frac{\nabla Z^{m} h}{\sqrt{1+|\nabla h|^2}}
- \frac{\nabla h \langle\nabla h,\nabla Z^{m} h\rangle}{\sqrt{1+|\nabla h|^2}^3} 
+ C^{m} (S)\Big) \Big)\right\}    \\
&\times \Big( \p_t Z^m h - v^b\cdot (Z^m \mathbf{N}) + C^m(KB)\Big) + R_{B} \\
=& - g \frac{\dd}{\dd t} \int_{\p S} |Z^m h|^2 - \Big\{ \frac{d}{dt}\int_{\partial S}\frac{{\left|\nabla Z^{m} h\right|}^2}{\sqrt{1+{|\nabla h|}^2}}
-  \int_{\partial S}\frac{{\left|\nabla Z^{m} h\right|}^2}{{\sqrt{1+{|\nabla h|}^2}}^3}\langle\nabla h,\nabla\partial_t h\rangle\Big\}   \\
& + R_S + R_C + R_B + P_1 + P_2 + P_3 ,  
\end{split}
\end{equation}
where $R_S$ and $R_C$ are defined in (\ref{RS}) and (\ref{RC}), and $R_B$, $P_1$, $P_2$, and $P_3$ are defined by 
\begin{equation} \label{high energy comms}
\begin{split}
R_B &:= 2\int_{\p S}(C^{m}(B)-\varepsilon(\Theta^{m}(v))^b \mathbf{N})\cdot (Z^{m} v)^b \\
P_1 &:=
2\int_{\partial S} \Big(
\frac{\nabla Z^{m} h}{\sqrt{1+{|\nabla h|}^2}} - 
\frac{\nabla h\langle\nabla h,\nabla Z^{m} h\rangle}{{\sqrt{1+{|\nabla h|}^2}}^3} + C^{m} (S)
\Big) \cdot\nabla \Big( v^b\cdot(Z^{m} \mathbf{N})+C^{m} (KB)
\Big)  \\
& \quad + 2 \int_{\partial S} \Big(
\frac{\nabla h\langle\nabla h,\nabla Z^{m} h\rangle}{{\sqrt{1+{|\nabla h|}^2}}^3} - C^{m} (S)
\Big) \cdot\nabla \Big( \partial_t(Z^{m} h) \Big)   ,   \\
P_2 &:= 2\int_{\p S}\left\{  q-gh-2\varepsilon(S^\varphi v) +  \nabla\cdot\frac{\nabla h}{\sqrt{1+|\nabla h|^2}}  \right\}
(Z^{m} \mathbf{N})\cdot(Z^{m} v) ,   \\
P_3 &:= 2g\int_{\p S}Z^{m} h \Big(v^b\cdot Z^{m}\mathbf{N} + C^{m}(KB)\Big)   .
\end{split}
\end{equation}

{\color{blue}
} 

\noindent Now we estimate $R_S$, $R_C$, $R_B$, $P_1$, $P_2$, and $P_3$.\\

\noindent\textbf{1) Estimate of $R_B$ }
\begin{equation} \label{RB est}
\begin{split}
|R_B| &= 2\Big|\int_{\p S} (C^{m}(B) - \varepsilon (\Theta^{m}(v) )^b \mathbf{N})\cdot (Z^{m} v)^b \dd x \dd y \Big|   \\
& \lesssim \Big| C^{m}(B)-\varepsilon (\Theta^{m}(v) )^b \mathbf{N}\Big|_{L^2(\partial S)} |({Z^{m} v})^b|_{L^2(\partial S)}   \\
& \lesssim \varepsilon \Lambda\big(\frac{1}{c_0},\left|h\right|_{Y^{[\frac{m}{2}],1}} + \left\|\nabla v\right\|_{Y^{[\frac{m}{2}],0}}\big) 
( \left|h\right|_{X^{m-1,1}} + |v^b|_{X^{m,\frac{1}{2}}} ) \|v\|_{X^m}  \\
& \lesssim \varepsilon \Lambda\big(\frac{1}{c_0},\left|h\right|_{Y^{[\frac{m}{2}],1}} + \left\|\nabla v\right\|_{Y^{[\frac{m}{2}],0}}\big) 
( \left|h\right|^2_{X^{m-1,1}} + |v^b|^2_{X^{m,\frac{1}{2}}} + \|v\|^2_{X^m} ), 
\end{split}
\end{equation}
by Proposition \ref{commutator est} and \ref{high order eq}.

\noindent\textbf{2) Estimate of $P_2$ }
Note that
\begin{equation*}
\begin{split}
P_2 &:= 2\int_{\p S}\Big\{  q-gh-2\varepsilon(S^\varphi v)+ \nabla\cdot\frac{\nabla h}{\sqrt{1+|\nabla h|^2}}  \Big\}
(Z^{m} \mathbf{N})\cdot(Z^{m} v)    \\
&=2\int_{\p S}\left\{ q^{NS}\vert_{z=0} I-2\varepsilon(S^\varphi v) \right\}(Z^{m} \mathbf{N})\cdot(Z^{m} v)    \\
&=4\varepsilon\int_{\p S}\left\{ (S^\varphi v)\mathbf{n}\cdot\mathbf{n}I-(S^\varphi v) \right\}(Z^{m} \mathbf{N})\cdot(Z^{m} v) .
\end{split}
\end{equation*}
Therefore,
\begin{equation} \label{P2 est}
\begin{split}
|P_2| &= \Big|{2\int_{\p S}\Big\{  q-gh-2\varepsilon(S^\varphi v) + \nabla\cdot\frac{\nabla h}{\sqrt{1+|\nabla h|^2}}  \Big\}(Z^{m} \mathbf{N})\cdot(Z^{m} v) \dd x \dd y}\Big|   \\
&\lesssim 2\varepsilon \Big|Z^{m} \mathbf{N} \Big|_{-\frac{1}{2}} 
\Big|\left\{(S^\varphi v)\mathbf{n}\cdot\mathbf{n}I-(S^\varphi v)\right\}^b (Z^{m} v)^b \Big|_{\frac{1}{2}}   \\
&\lesssim \varepsilon \left|h\right|_{X^{m,\frac{1}{2}}}  |v^b |_{X^{m,\frac{1}{2}}} \Big| (S^\varphi v)\mathbf{n}\cdot\mathbf{n}I-(S^\varphi v) \Big|_{1,\infty}   \\
&\lesssim \varepsilon \Lambda\big(\frac{1}{c_0},\left|h\right|_{Y^{[\frac{m}{2}],1}} + \left\|\nabla v\right\|_{Y^{[\frac{m}{2}],0}}\big)  |h|_{X^{m,\frac{1}{2}}} |v^b|_{X^{m,\frac{1}{2}}},
\end{split}
\end{equation}
by Proposition \ref{commutator est}.  \\
\noindent\textbf{3)  Estimate of $P_3$ }
Since,
\begin{equation*}
P_3 := 2g\int_{\p S}Z^{m} h \Big(v^b\cdot Z^{m}\mathbf{N} + C^{m}(KB)\Big) ,
\end{equation*}
\begin{equation} \label{P3 est}
\begin{split}
|P_3| &\leq \Big|\int_{\p S} Z^{m} h (v^b\cdot Z^{m}\mathbf{N}+C^{m}(KB) )   \Big|  \\
&\leq |Z^{m} h(t)|_{L^2(\partial S)}\Big\{ |v^b\cdot Z^{m}\mathbf{N}|_{L^2(\partial S)} + \left|C^{m} (KB)(t)\right|_{L^2(\partial S)} \Big\}  \\
&\lesssim \left|h\right|_{X^{m,0}}\left\{|v^b|_{L^\infty}\left|h\right|_{X^{m,1}} + \left\|C^m(KB)\right\|\right\}  \\
&\lesssim \Lambda\big(\frac{1}{c_0},\left|h\right|_{Y^{[\frac{m}{2}],1}} + \left\|\nabla v\right\|_{Y^{[\frac{m}{2}],0}}\big) \Big(\left\|v\right\|^2_{X^{m,0}} + \left\|\nabla v\right\|^2_{X^{m-1,0}} + \left|h\right|^2_{X^{m,1}} \Big),
\end{split}
\end{equation}
by Proposition \ref{commutator est} and $\|C^m(KB)\|$ estimate in Proposition \ref{high order eq}.

\noindent\textbf{4) Estimate of $R_C$ }
\begin{equation} \label{RC est}
\begin{split}
\left|R_C\right| &\lesssim \Lambda_0
\Big( \left\|{C^{m}(d)}\right\|_{L^2}\left\|{Z^{m} q(t)}\right\|_{L^2} + \left\|{T^{m}(v)}\right\|_{L^2}\left\|{Z^{m} v(t)}\right\|_{L^2} 
+ \left\|{C^{m}(q)}\right\|_{L^2}\left\|{Z^{m} v(t)}\right\|_{L^2}
\Big)    \\
&\lesssim \Lambda_0 
\Big( \left\|{C^{m}(d)}\right\|\left\|{q(t)}\right\|_{X^{m,0}} + \left\|{T^{m}(v)}\right\|\left\|{v(t)}\right\|_{X^{m,0}} + \left\|{C^{m}(q)}\right\|\left\|{v(t)}\right\|_{X^{m,0}}
\Big)   \\
&\lesssim \Lambda\big(\frac{1}{c_0},\left|h\right|_{Y^{[\frac{m}{2}],1}} + \left\|\nabla v\right\|_{Y^{[\frac{m}{2}],0}} + \left\|\nabla q\right\|_{Y^{[\frac{m}{2}],0}}\big) \Big( \left\|v\right\|^2_{X^{m,0}} + \left\|\nabla v\right\|^2_{X^{m-1,0}} + \left|h\right|^2_{X^{m,\frac{1}{2}}} + \left\|\nabla q\right\|^2_{X^{m-1,0}} \Big).
\end{split}
\end{equation}

\noindent\textbf{5) Estimate of $R_S$ }
From definition of $R_S$,
\begin{eqnarray*}
R_S &=& 2\varepsilon \int_S\left\{\nabla^\varphi\cdot\Theta^{m}(v)-D^{m}(S^\varphi v)\right\}\cdot Z^{m} v dV_t   \\
&=& -2\varepsilon\int_S \Theta^{m}(v):\nabla^\varphi(Z^{m} v) 
+ 2\varepsilon\int_{\partial S}\Theta^{m}(v)N\cdot (Z^{m} v)^b 
+ 2\varepsilon\sum_{i,j}\int_S C_j^{m}(S^\varphi v)_{ij}(Z^{m} v)_i dV_t,
\end{eqnarray*}
We have only two types of integrals. ($m = m_1 + m_2$ and both are non-zero.)
\begin{eqnarray*}
I_1 &:=& \int_S \partial_z(Z^{m} v_i)Z^{m_1}(S^\varphi v)_{ij}Z^{m_2}\Big(\frac{\partial_i\varphi}{\partial_z\varphi}\Big),  \\
I_2 &:=& \int_{\partial S} (Z^{m} v_i)Z^{m_1}(S^\varphi v)_{ij}Z^{m_2}\Big(\frac{\partial_i\varphi}{\partial_z\varphi}\Big).
\end{eqnarray*}
By $L^\infty L^2 L^2$ Holder inequality, we get,
\begin{equation*}
\begin{split}
\left|I_1\right| &\lesssim \Lambda\big(\frac{1}{c_0},\left|h\right|_{Y^{[\frac{m}{2}],1}} + \left\|\nabla v\right\|_{Y^{[\frac{m}{2}],0}} \big) \Big(\left\|\nabla v\right\|^2_{X^{m-1,0}} + \left\|S^\varphi v\right\|^2_{X^{m-1,0}} + \left|h\right|_{X^{m-1,\frac{1}{2}}} \Big),  \\
\left|I_2\right| &\lesssim \Lambda\big(\frac{1}{c_0},\left|h\right|_{Y^{[\frac{m}{2}],1}} + \left\|\nabla v\right\|_{Y^{[\frac{m}{2}],0}} \big) \Big( \left\|v\right\|^2_{X^{m}} + \left\|\nabla v\right\|^2_{X^{m-1,0}} + \left|h\right|_{X^{m-1,\frac{1}{2}}} \Big).
\end{split}
\end{equation*}
Therefore, we get
\begin{equation} \label{RS est}
\left|R_S\right| \lesssim \varepsilon \Lambda\big(\frac{1}{c_0},\left|h\right|_{Y^{[\frac{m}{2}],1}} + \left\|\nabla v\right\|_{Y^{[\frac{m}{2}],0}} \big) \Big( \left\|v\right\|^2_{X^{m}} + \left\|\nabla v\right\|^2_{X^{m-1,0}} + \left\|S^\varphi v\right\|^2_{X^{m-1,0}} + \left|h\right|^2_{X^{m-1,\frac{1}{2}}} \Big).
\end{equation}

\noindent\textbf{6) Estimate of $P_1$ } \\
Let us decompose $P_1 = P_{1,1} + P_{1,2}$, where
\begin{equation} \label{P11}
\begin{split}
P_{1,1} &:=
2\int_{\partial S} \Big(
\frac{\nabla Z^{m} h}{\sqrt{1+{|\nabla h|}^2}} - 
\frac{\nabla h\langle\nabla h,\nabla Z^{m} h\rangle}{{\sqrt{1+{|\nabla h|}^2}}^3} + C^{m} (S)
\Big) \cdot\nabla \Big( v^b\cdot(Z^{m} \mathbf{N})+C^{m} (KB) \Big)   \\
P_{1,2} &:= 2 \int_{\partial S} \Big(
\frac{\nabla h\langle\nabla h,\nabla Z^{m} h\rangle}{{\sqrt{1+{|\nabla h|}^2}}^3} - C^{m} (S)
\Big) \cdot\nabla \Big( \partial_t(Z^{m} h) \Big)    .
\end{split}
\end{equation}
Let treat $P_{1,2}$ first. We split $P_{1,2} := P_{1,2,1} + P_{1,2,2}$ again, where
\begin{equation} \label{5.38}
\begin{split}
P_{1,2,1} &:= 2\int_{\partial S} \frac{\langle\nabla h,\nabla Z^{m} h\rangle}{{\sqrt{1+{|\nabla h|}^2}}^3}\langle\nabla h,\nabla Z^{m} \partial_t h\rangle  := P_{1,2,1,1} + P_{1,2,1,2} + P_{1,2,1,3} \\ 
&:= \frac{d}{dt}\int_{\partial S} \frac{\langle\nabla h,\nabla Z^{m} h\rangle^2}{{\sqrt{1+{|\nabla h|}^2}}^3}  
- 2\int_{\partial S} \frac{\langle\nabla h,\nabla Z^{m} h\rangle}{{\sqrt{1+{|\nabla h|}^2}}^3}\langle\nabla\partial_t h,\nabla Z^{m} h\rangle   \\
&\quad + 3\int_{\partial S}\frac{\langle\nabla h,\nabla\partial_t h\rangle}{{\sqrt{1+{|\nabla h|}^2}}^3}\langle\nabla h,\nabla Z^{m}h\rangle^2   ,  \\
P_{1,2,2} &= -2 \int_{\partial S} C^{m} (S) \cdot\nabla \partial_t Z^{m} h    .
\end{split}
\end{equation}
$P_{1,2,1,1} := \frac{d}{dt}\int_{\partial S} \frac{\langle\nabla h,\nabla Z^{m} h\rangle^2}{{\sqrt{1+{|\nabla h|}^2}}^3}  \dd x \dd y $ should be absorbed by the energy on the left hand side under the assumption that $|h|_{1,\infty}$ is bounded. 
Controls of $P_{1,2,1,2}$ and $P_{1,2,1,3}$ are as following,
\begin{equation} \label{P1212 est}
|P_{1,2,1,2} + P_{1,2,1,3}| \leq \Lambda\big(\frac{1}{c_0},\left|h\right|_{Y^{[\frac{m}{2}],1}} + \left\|\nabla v\right\|_{Y^{[\frac{m}{2}],0}} \big) \left|h\right|^2_{X^{m,1}}.
\end{equation}
To estimate $P_{1,2,2}$, we use integration by part in horizontal variables to get,
\begin{equation} \label{P122 est}
\begin{split}
	|P_{1,2,2}| & \leq \Big|\int_{\partial S} C^{m}(S)\cdot \nabla (\partial_t Z^{m} h )  \Big|
	= \Big|\int_{\partial S} \nabla\cdot C^{m}(S) \partial_t Z^{m} h \Big| \\
	& \leq \|\nabla\cdot C^{m}(S)\| |\p_{t} h|_{X^{m-1,1}} \\
	&\lesssim \Lambda\big(\frac{1}{c_0},\left|h\right|_{Y^{[\frac{m}{2}],1}} + \left\|\nabla v\right\|_{Y^{[\frac{m}{2}],0}} \big) |h|_{X^{m,1}} \big( \|v\|_{X^{m,0}} + |h|_{X^{m,1}} \big).
	\end{split}
\end{equation}
For $P_{1,1}$, by definition (\ref{P11}), it can be controlled by 
\begin{equation} \label{P11 est}
\begin{split}
|P_{1,1}| &\leq \Big|\frac{\nabla Z^{m} h}{\sqrt{1+{|\nabla h|}^2}} 
- \frac{\nabla h \langle \nabla h,\nabla Z^{m} h \rangle }{{\sqrt{1+{|\nabla h|}^2}}^3} + C^{m} (S) \Big|_{\frac{1}{2}}
\Big| v^b\cdot(Z^{m} \mathbf{N}) + C^{m} (KB) \Big|_{\frac{1}{2}} \\
& \leq \Lambda\big(\frac{1}{c_0},\left|h\right|_{Y^{[\frac{m}{2}],1}} + \left\|\nabla v\right\|_{Y^{[\frac{m}{2}],0}} \big) \Big( \left\|v\right\|^2_{X^{m,0}} + \left\|\nabla v\right\|^2_{X^{m-1,0}} + \left|h\right|^2_{X^{m-1,\frac{5}{2}}} \Big),
\end{split}
\end{equation}
by Proposition \ref{commutator est}, \ref{aniso embed}, \ref{horizon est}, and $C^{m}(KB)$ estimate in Proposition \ref{high order eq}.   \\
We combine (\ref{P11}), (\ref{5.38}), (\ref{P1212 est}), (\ref{P122 est}), and (\ref{P11 est}), to get 
\begin{equation} \label{P1 est}
\begin{split}
P_1 &\leq  \frac{d}{dt}\int_{\partial S} \frac{\langle\nabla h,\nabla Z^{m} h\rangle^2}{{\sqrt{1+{|\nabla h|}^2}}^3}     \\ 
&\quad  + \Lambda\big(\frac{1}{c_0},\left|h\right|_{Y^{[\frac{m}{2}],1}} + \left\|\nabla v\right\|_{Y^{[\frac{m}{2}],0}} \big) \Big( \left\|v\right\|^2_{X^{m,0}} + \left\|\nabla v\right\|^2_{X^{m-1,0}} + \left|h\right|^2_{X^{m-1,\frac{5}{2}}} \Big).
\end{split}
\end{equation}

\noindent\textbf{7) Estimate of $\Big| \int_{\partial S}\frac{{\left|\nabla Z^{m} h\right|}^2}{{\sqrt{1+{|\nabla h|}^2}}^3}\langle\nabla h,\nabla\partial_t h\rangle \Big|$ . }
\begin{equation} \label{last est}
\begin{split}
\Big| \int_{\partial S}\frac{{\left|\nabla Z^{m} h\right|}^2}{{\sqrt{1+{|\nabla h|}^2}}^3}\langle\nabla h,\nabla\partial_t h\rangle \Big| &\leq \Lambda(\frac{1}{c_0}) \int_{\partial S} \left|\nabla Z^{m} h\right|^2 \Big|\langle\nabla h,\nabla\partial_t h\rangle\Big| dA   \\
&\leq \Lambda(\frac{1}{c_0}) \left| \langle\nabla h,\nabla (v^b\cdot \mathbf{N}) \rangle\right|_{L^\infty}\left|h\right|^2_{X^{m,1}}
\leq \Lambda\Big(\frac{1}{c_0}, |h|_{2,\infty} + \|\nabla v\|_{1,\infty} \Big) |h|^2_{X^{m,1}}.
\end{split}
\end{equation}

We collect (\ref{basic high}), (\ref{BDRY integral}), and (\ref{RB est})-(\ref{last est}), to get, (with small $\varepsilon \ll 1$,)
\begin{equation}
\begin{split}
& \frac{d}{dt}\int_S\left|Z^{m} v\right|^2dV_t
+ 4\varepsilon\int_S\left|S^\varphi (Z^{m} v)\right|^2dV_t + g \frac{\dd}{\dd t} \int_{\p S} |Z^m h|^2    \\
& + \frac{d}{dt} \Big( \int_{\partial S}\frac{{\left|\nabla Z^{m} h\right|}^2}{\sqrt{1+{|\nabla h|}^2}} - \int_{\partial S} \frac{\langle\nabla h,\nabla Z^{m} h\rangle^2}{{\sqrt{1+{|\nabla h|}^2}}^3} \Big) \leq \Lambda\big(\frac{1}{c_0},\left|h\right|_{Y^{[\frac{m}{2}],1}} + \left\|\nabla v\right\|_{Y^{[\frac{m}{2}],0}} + \left\|\nabla q\right\|_{Y^{[\frac{m}{2}],0}}\big)  \\
&\times \Big( \left\|v\right\|^2_{X^{m,0}} + \left\|\nabla v\right\|^2_{X^{m-1,0}} + \varepsilon \|\nabla v \|^2_{X^{m,0}} + \left|h\right|^2_{X^{m-1,\frac{5}{2}}} + \left\|\nabla q\right\|^2_{X^{m-1,0}} \Big) .  \\
\end{split}
\end{equation}
Note that 
\begin{equation} \label{h positivity}
	\int_{\partial S} \frac{{\left|\nabla Z^{m} h\right|}^2}{\sqrt{1+{|\nabla h|}^2}} - \int_{\partial S} \frac{\langle\nabla h,\nabla Z^{m} h\rangle^2}{{\sqrt{1+{|\nabla h|}^2}}^3} \geq \int_{\p S} \frac{{\left|\nabla Z^{m} h\right|}^2}{\sqrt{1+{|\nabla h|}^2}} \Big( \frac{1}{1+|\nabla h|^2} \Big) \geq \Lambda(c_0) \int_{\p S} \frac{{\left|\nabla Z^{m} h\right|}^2}{\sqrt{1+{|\nabla h|}^2}} .
\end{equation}
We consider every $Z^m \neq \p_t^m$ and sum for all the result. Then we use (\ref{h positivity}) and Proposition \ref{korn}, \ref{aniso embed}, and integrate in time for $t\in [0,T]$, to get the result,
\begin{equation*}
\begin{split}
& \|v(T)\|^2_{X^{m-1,1}} + \|h(T)\|^2_{X^{m-1,2}} + \varepsilon\int_{0}^{T} \|\nabla v(t)\|^2_{X^{m-1,1}} dt \\
& \leq \Lambda(\frac{1}{c_0}) \sum_{Z^m\neq\p_t^m} ( \|(Z^m v)(0)\|_{L^2} + |(Z^m h)(0)|_{H^1} ) \\
&\quad + \Lambda(\frac{1}{c_0})\int_{0}^{T} \Lambda\big(\frac{1}{c_0},|h(t)|_{Y^{[\frac{m}{2}],1}} + \|\nabla v(t)\|_{Y^{[\frac{m}{2}],0}} + \|\nabla q(t)\|_{Y^{[\frac{m}{2}],0}}\big) \\
&\quad\quad \times \Big( \|v(t)\|^2_{X^{m,0}} + \|\nabla v(t)\|^2_{X^{m-1,0}} + \varepsilon \|\nabla v(t) \|^2_{X^{m,0}} + |h(t)|^2_{X^{m-1,\frac{5}{2}}} + \|\nabla q(t)\|^2_{X^{m-1,0}} \Big) dt. \\
\end{split}
\end{equation*}
We use pressure estimates Proposition \ref{q_E est}, \ref{q_NS est}, \ref{q_S est}, and \ref{q_infty est} to finish our proof.

\end{proof}
In Proposition \ref{energy estimate}, we should control $|h|_{X^{m-1,\frac{5}{2}}}$. In the next section we claim that this term can be controlled by $|\p_{t} h|_{X^{m-1,1}} \leq | h|_{X^{m,1}}$.  \\

\section{Dirichlet-Neumann Operator estimate on the boundary}
In this section, we claim that, on the boundary $\partial_x^{3/2} h$ can be controlled by $\partial_t h$. We start with this section with a lemma which is needed to prove the next proposition.
\begin{lemma} \label{DN lemma}
	There exists $c>0$ such that for every $h \in W^{1,\infty}(\mathbb{R}^2)$,
	\begin{equation}
	(G[h]f^b,f^b) \geq c (1+\|h\|_{W^{1,\infty}(\mathbb{R}^2)})^{-2} \left\|\frac{|\nabla|}{(1+|\nabla|)^{1/2}}f^b \right\|^2_{L^2(\mathbb{R}^2)}, \,\,\,\,\,\,\, \forall f^b \in H^{\frac{1}{2}}(\mathbb{R}^2),
	\end{equation}
	where $G[h]f^b$ means Dirichlet-Neumann operator,
	$$
	G[h]f^b := (\nabla f)^b \cdot \mathbf{N}, 
	$$
	where $\mathbf{N} = (-\nabla h, 1)$ and $f$ solves harmonics equation, $\Delta f = 0$.
\end{lemma}
\begin{proof}
	See Proposition 3.4 of [2].
\end{proof}
We can apply above lemma for $f = p^S$, since $\Delta p^S = \Delta^\varphi q^S = 0$ in (\ref{q_S}).
\begin{proposition} \label{dirineumann}
	Assume that $(v,\varphi,q,h)$ is a smooth solution on $t\in[0,T]$ of the system (\ref{1.12})-(\ref{1.15}) and
	\begin{equation}
	\partial_z\varphi \geq c_0 ,\quad \forall t\in [0,T].
	\end{equation}
	Then $h$ enjoys the following estimate.
	\begin{equation} \label{DN est}
	\begin{split}
	&\int_0^T | Z^{m-1} \nabla h|^2_{\frac{3}{2}} dt \leq \Lambda\Big( \frac{1}{c_0},\|h(0)\|_{X^{m,1}} \Big) + \Lambda_{m,\infty}(T) \Big\{ {\theta} |Z^{m-1} \p_{t} h(T)|^2_{L^2} + \sqrt{T} |Z^{m-1}\nabla \p_{t} h(T)|^2_{L^2} \Big\}  \\
	&\quad + (1+\sqrt{T}) \int_{0}^{T} \Lambda_{m,\infty}(t)\Big( \|v\|^2_{X^{m-1,1}} + \|\nabla v\|^2_{X^{m-1,0}} + |h|^2_{X^{m,1}} + \varepsilon\|\nabla v\|^2_{X^{m-1,1}} + \varepsilon\|v\|^2_{X^{m-1,2}} \Big) dt, 
	\end{split}
	\end{equation}
	where $\Lambda_{m,\infty}(t)$ is defined in (\ref{Lambda infty}) and $\theta > 0$ is sufficiently small.
\end{proposition}
\begin{proof} 
	From kinematic boundary condition (\ref{1.14}), $h_t = v^b \cdot \mathbf{N}$, we get $\partial_{tt} h = v_t^b\cdot \mathbf{N} + v^b\cdot \mathbf{N}_t$. We apply $Z^{m-1}$ to this equation, where $\alpha_3=0$, because we are on the boundary, i.e. we do not apply any normal derivatives. 
	\begin{equation} \label{DN}
	\begin{split}
	\partial_{tt}\Big(Z^{m-1} h\Big) & = \Big(Z^{m-1}v_t^b\Big)\cdot \mathbf{N} + v_t^b\cdot\Big(Z^{m-1}\mathbf{N}\Big) + \left[Z^{m-1}, v_t^b, \mathbf{N}\right] + \Big(Z^{m-1}v^b\Big)\cdot \mathbf{N}_t \\
	&\quad + v^b\cdot\Big(Z^{m-1} \mathbf{N}_t\Big) + \left[Z^{m-1}, v^b, \mathbf{N}_t\right] \\
	&:= \Big\{ -Z^{m-1}(v\cdot\nabla^\varphi v)^b - Z^{m-1}(\nabla^\varphi q^E + \nabla^\varphi q^{NS} + \nabla^{\varphi} q^{S})^b + 2\varepsilon Z^{m-1}(\nabla^\varphi\cdot S^\varphi v)^b \Big\}
	\cdot\mathbf{N} \\
	&\quad + ( I_1 + I_2 + I_3 + I_4 + I_5 ) \\
	&:= G[h]V^b + ( I_1 + I_2 + I_3 + I_4 + I_5 ) + (J_1 + J_2 + J_3),
	\end{split}
	\end{equation}
	where $V^{b}$, $I_{1,2,3,4,5}$, and $J_{1,2,3}$ are defined by 
	\begin{eqnarray}
	\begin{split}
	V^b & := Z^{m-1} \nabla\cdot\Big(\frac{\nabla h}{\sqrt{1+|\nabla h|^2}}\Big), \\
	I_1 & := v_t^b\cdot\Big(Z^{m-1}\mathbf{N}\Big), \\
	I_2 & := \left[Z^{m-1}, v_t^b, \mathbf{N}\right], \\
	I_3 & := \Big(Z^{m-1}v^b\Big)\cdot \mathbf{N}_t, \\
	I_4 & := v^b\cdot\Big(Z^{m-1} \mathbf{N}_t\Big), \\
	I_5 & := \left[Z^{m-1}, v^b, \mathbf{N}_t\right], \\
	J_1 & := -Z^{m-1}(v\cdot\nabla^\varphi v)^b\cdot\mathbf{N}, \\
	J_2 & := - Z^{m-1}(\nabla^\varphi q^E + \nabla^\varphi q^{NS})^b \cdot\mathbf{N}, \\
	J_3 & :=  2\varepsilon Z^{m-1}(\nabla^\varphi\cdot S^\varphi v)^b \cdot\mathbf{N}.
	\end{split}
	\end{eqnarray}
	We take dot product with $V^b$ and $\int_0^T\int_{\partial S}$ to the equation (\ref{DN}), {i.e},
	\begin{equation} \label{whole DN}
	\begin{split}
	\int_0^T (G[h]V^b,V^b) dt &\leq \int_0^T \int_{\p S} Z^{m-1} h_{tt} V^b dt  \\
	&\quad  + \int_0^T \int_{\p S} ( I_1 + I_2 + I_3 + I_4 + I_5 + J_1 + J_2 + J_3) V^b  dt.
	\end{split}
	\end{equation}
	We estimate upper bound of RHS terms and lower bound of LHS.  \\
	\noindent\textbf{1) Estimate of} $ \int_0^T \int_{\p S} Z^{m-1} h_{tt} V^b $.
	\begin{equation} \label{7.4}
	\int_0^T \int_{\partial S} Z^{m-1} h_{tt}\nabla\cdot Z^{m-1}\Big(\frac{\nabla h}{\sqrt{1+|\nabla h|^2}}\Big) dt = -\int_0^T \int_{\partial S}Z^{m-1}\nabla h_{tt} \cdot Z^{m-1}\Big(\frac{\nabla h}{\sqrt{1+|\nabla h|^2}}\Big) dt	
	\end{equation}
	\begin{equation*}
	\begin{split}
	& = -\Big[\int_{\partial S} Z^{m-1}\nabla h_t \cdot Z^{m-1} \Big(\frac{\nabla h}{\sqrt{1+|\nabla h|^2}}\Big)\Big]_0^T + \int_{0}^{T}\int_{\partial S} Z^{m-1}\nabla h_t\cdot Z^{m-1}\partial_t(\frac{\nabla h}{\sqrt{1+|\nabla h|^2}}\Big) dt \\
	& \leq \Lambda\Big( \frac{1}{c_0},\|h(0)\|_{X^{m,1}} \Big) - \int_{\partial S} Z^{m-1}\nabla h_t\cdot Z^{m-1}\Big(\frac{\nabla h}{\sqrt{1+|\nabla h|^2}}\Big)\biggr\rvert_{t=T} \\
	&\quad + \int_0^T \Lambda(\frac{1}{c_0}, |h|_{Y^{[\frac{m}{2}],1}}) \|\nabla h_t\|^2_{X^{m-1,0}} dt \\
	\end{split}
	\end{equation*}

	To estimate $\int_{\partial S} Z^{m-1}\nabla h_t\cdot Z^{m-1}\Big(\frac{\nabla h}{\sqrt{1+|\nabla h|^2}}\Big)\biggr\rvert_{t=T}$ on the RHS, we compute
	\begin{equation*}
	\begin{split}
	&\int_{\partial S} Z^{m-1}\nabla h_t\cdot Z^{m-1}\Big(\frac{\nabla h}{\sqrt{1+|\nabla h|^2}}\Big)\biggr\rvert_{t=T}  \\
	&= \int_{\partial S} Z^{m-1}\nabla h_t (T) \Big\{ Z^{m-1}\Big(\frac{\nabla h}{\sqrt{1+|\nabla h|^2}}\Big)(0) +  \int_0^T Z^{m-1}\partial_t\Big(\frac{\nabla h}{\sqrt{1+|\nabla h|^2}}\Big) dt \Big\}   \\
	&= -\int_{\partial S} Z^{m-1} h_t(T)  Z^{m-1}\nabla \Big(\frac{\nabla h}{\sqrt{1+|\nabla h|^2}}\Big)(0)  + \int_{\partial S} Z^{m-1}\nabla h_t(T)  \int_0^T Z^{m-1}\partial_t\Big(\frac{\nabla h}{\sqrt{1+|\nabla h|^2}}\Big) dt   \\
	&\leq \Lambda\Big( \frac{1}{c_0},\|h(0)\|_{X^{m,1}} \Big) + {\theta} |Z^{m-1} h_t(T)|^2_{L^2} \\
	&\quad + |Z^{m-1} \nabla \p_{t}h(T)|_{L^{2}(\p S)}  \Big| \int_{0}^{T} Z^{m-1} \nabla \p_{t} \nabla h(t) dt \Big|_{L^{2}(\p S)}.  \\
	\end{split}
	\end{equation*}
	Using Jensen's inequality for the last term, we get
	\begin{equation*}
	\begin{split}
	&\int_{\partial S} Z^{m-1}\nabla h_t\cdot Z^{m-1}\Big(\frac{\nabla h}{\sqrt{1+|\nabla h|^2}}\Big)\biggr\rvert_{t=T}  \\
	&\leq \Lambda\Big( \frac{1}{c_0},\|h(0)\|_{X^{m,1}} \Big) + {\theta} |Z^{m-1} h_t(T)|^2_{L^2} \\
	&\quad + \sqrt{T} |Z^{m-1}\nabla h_t(T)|_{L^2} \Big\{ \int_0^T \Big| Z^{m-1}\partial_t\Big(\frac{\nabla h}{\sqrt{1+|\nabla h|^2}}\Big) \Big|^2_{L^2} dt \Big\}^{1/2}  \\
	&\leq \Lambda\Big( \frac{1}{c_0},\|h(0)\|_{X^{m,1}} \Big) + {\theta} |Z^{m-1} h_t(T)|^2_{L^2} + 2\sqrt{T} |Z^{m-1}\nabla h_t(T)|^2_{L^2} \\
	&\quad + \sqrt{T} \int_0^T \Lambda(\frac{1}{c_0}, |h|_{Y^{[\frac{m}{2}],1}}) |Z^{m-1}\nabla h_t(t) |^2_{L^2} dt,
	\end{split}
	\end{equation*}
	where sufficiently small ${\theta}$ came from Young's inequality.	\\ 

	Now we estimate other terms involving $I_k$ and $J_k$. \\
	
	\noindent \textbf{2) Estimate of} $\int_0^T\int_{\partial S} I_1 \nabla\cdot Z^{m-1}\Big(\frac{\nabla h}{\sqrt{1+|\nabla h|^2}}\Big) dt $. Using \ref{commutator est},
	$$
	\int_0^T\int_{\partial S} v_t^b\cdot (Z^{m-1} \mathbf{N}) \nabla\cdot Z^{m-1}\Big(\frac{\nabla h}{\sqrt{1+|\nabla h|^2}}\Big) dt \leq \int_0^T \Lambda(\frac{1}{c_0},|h|_{Y^{[\frac{m}{2}],1}} + \|\nabla v\|_{2,\infty}) |h|^2_{X^{m-2,2}} dt.
	$$
	\noindent \textbf{3) Estimate of} $\int_0^T\int_{\partial S} I_2 \nabla\cdot Z^{m-1}\Big(\frac{\nabla h}{\sqrt{1+|\nabla h|^2}}\Big) dt$. Using \ref{commutator est},
	\begin{equation*}
	\begin{split}
	&\int_0^T \int_{\partial S} \left[Z^{m-1}, v_t^b, \mathbf{N}\right] \nabla\cdot Z^{m-1}\Big(\frac{\nabla h}{\sqrt{1+|\nabla h|^2}}\Big)dt \\
	&\leq \int_0^T \Lambda(\frac{1}{c_0},|h|_{Y^{[\frac{m}{2}],1}} + \|\nabla v\|_{Y^{[\frac{m}{2}],0}}) \Big( \|v\|^2_{X^{m-1,1}} + \|\nabla v\|^2_{X^{m-1,0}} + |h|^2_{X^{m-2,2}} \Big) dt,
	\end{split}
	\end{equation*}
	We used trace estimate Lemma \ref{aniso embed}. \\
	\noindent \textbf{4) Estimate of} $\int_0^T \int_{\partial S} I_3 \nabla\cdot Z^{m-1}\Big(\frac{\nabla h}{\sqrt{1+|\nabla h|^2}}\Big) dt$.
	\begin{equation*}
	\begin{split}
	&\int_0^T\int_{\partial S} \Big(Z^{m-1}v^b\Big)\cdot \mathbf{N}_t \nabla\cdot Z^{m-1}\Big(\frac{\nabla h}{\sqrt{1+|\nabla h|^2}}\Big)dt \\
	&\leq \int_0^T \Lambda(\frac{1}{c_0},|h|_{Y^{[\frac{m}{2}],1}} + \|\nabla v\|_{2,\infty}) \Big( \|v\|^2_{X^{m-1,1}} + \|\nabla v\|^2_{X^{m-1,0}} + |h|^2_{X^{m-2,2}} \Big)dt,
	\end{split}
	\end{equation*}
	similarly as above. \\
	\noindent \textbf{5) Estimate of} $\int_0^T \int_{\partial S} I_4 \nabla\cdot Z^{m-1}\Big(\frac{\nabla h}{\sqrt{1+|\nabla h|^2}}\Big) dt$.
	$$
	\int_0^T \int_{\partial S} v^b\cdot\Big(Z^{m-1} \mathbf{N}_t\Big) \nabla\cdot Z^{m-1}\Big(\frac{\nabla h}{\sqrt{1+|\nabla h|^2}}\Big) dt \leq \int_0^T \Lambda(\frac{1}{c_0},|h|_{Y^{[\frac{m}{2}],1}} + \|\nabla v\|_{2,\infty}) |h|^2_{X^{m-2,2}} dt.
	$$
	\noindent \textbf{6) Estimate of} $\int_0^T \int_{\partial S} I_5 \nabla\cdot Z^{m-1}\Big(\frac{\nabla h}{\sqrt{1+|\nabla h|^2}}\Big) dt$.
	\begin{equation*}
	\begin{split}
	& \int_0^T \int_{\partial S} \left[Z^{m-1}, v_t^b, \mathbf{N} \right] \nabla\cdot Z^{m-1}\Big(\frac{\nabla h}{\sqrt{1+|\nabla h|^2}}\Big) dt \\
	&\leq \int_0^T \Lambda(\frac{1}{c_0},|h|_{Y^{[\frac{m}{2}],1}} + \|\nabla v\|_{2,\infty}) \Big( \|v\|^2_{X^{m-1,1}} + \|\nabla v\|^2_{X^{m-1,0}} + |h|^2_{X^{m-2,2}} \Big) dt.
	\end{split}
	\end{equation*}
	\noindent \textbf{7) Estimate of} $\int_0^T\int_{\partial S} J_1 \nabla\cdot Z^{m-1}\Big(\frac{\nabla h}{\sqrt{1+|\nabla h|^2}}\Big) dt$. Using Lemma \ref{aniso embed}, we can replace $v^b$ into $v$ by giving $|\cdot|_{H^{-1/2}}$. Hence,
	\begin{equation*}
	\begin{split}
	&\int_0^T\int_{\partial S} -Z^{m-1}(v\cdot\nabla^\varphi v)^b\cdot\mathbf{N} \nabla\cdot Z^{m-1}\Big(\frac{\nabla h}{\sqrt{1+|\nabla h|^2}}\Big) dt \\
	& \leq \int_0^T \Big|Z^{m-1}(v\cdot\nabla^\varphi v)^b \Big|_{-\frac{1}{2}} \Big|\mathbf{N} \nabla\cdot Z^{m-1}\Big(\frac{\nabla h}{\sqrt{1+|\nabla h|^2}}\Big) \Big|_{\frac{1}{2}} dt  \\
	&\leq \int_0^T \Lambda(\frac{1}{c_0},|h|_{Y^{[\frac{m}{2}],1}} + \|\nabla v\|_{2,\infty}) (\|v\|_{X^{m-1,1}} + \|\nabla v\|_{X^{m-1,0}})|h|_{X^{m-1,\frac{5}{2}}} dt \\
	&\lesssim \int_0^T \Lambda(\frac{1}{c_0},|h|_{Y^{[\frac{m}{2}],1}} + \|\nabla v\|_{2,\infty})(\|v\|^2_{X^{m-1,1}} + \|\nabla v\|^2_{X^{m-1,0}}) dt \\
	&\quad + \theta\int_0^T \Lambda(\frac{1}{c_0},|h|_{Y^{[\frac{m}{2}],1}} + \|\nabla v\|_{2,\infty})|\nabla h|^2_{X^{m-1,\frac{3}{2}}} dt,
	\end{split}
	\end{equation*}
	where $\theta$ is sufficiently small constant from Young's inequality. \\
	\noindent \textbf{8) Estimate of} $\int_0^T\int_{\partial S} J_2 \nabla\cdot Z^{m-1}\Big(\frac{\nabla h}{\sqrt{1+|\nabla h|^2}}\Big) dt$. We perform estimate same as $J_1$ case. Hence, using pressure estimate in section 4,
	\begin{equation*}
	\begin{split}
	& \int_0^T \int_{\partial S} - Z^{m-1}(\nabla^\varphi q^E + \nabla^\varphi q^{NS})^b \cdot\mathbf{N} \nabla\cdot Z^{m-1}\Big(\frac{\nabla h}{\sqrt{1+|\nabla h|^2}}\Big) dt \\
	& \lesssim  \int_0^T \Lambda(\frac{1}{c_0},|h|_{Y^{[\frac{m}{2}],1}} + \|\nabla v\|_{2,\infty}) ( |h|^2_{X^{m-1,2}} + \|v\|^2_{X^{m-1,1}} + \varepsilon\|\nabla v\|^2_{X^{m-1,1}} + \varepsilon\|v\|^2_{X^{m-1,2}}) dt \\
	& + \theta\int_0^T \Lambda(\frac{1}{c_0},|h|_{Y^{[\frac{m}{2}],1}} + \|\nabla v\|_{2,\infty}) |\nabla h|^2_{X^{m-1,\frac{3}{2}}} dt.
	\end{split}
	\end{equation*}
	\noindent \textbf{9) Estimate of} $\int_0^T \int_{\partial S} J_3 \nabla\cdot Z^{m-1}\Big(\frac{\nabla h}{\sqrt{1+|\nabla h|^2}}\Big) dt$. Similarly as above estimate, we get
	\begin{equation*}
	\begin{split}
	& \int_0^T \int_{\partial S} 2\varepsilon Z^{m-1}(\nabla^\varphi\cdot S^\varphi v)^b \cdot\mathbf{N} \nabla\cdot Z^{m-1}\Big(\frac{\nabla h}{\sqrt{1+|\nabla h|^2}}\Big) dt \\
	& \leq \int_0^T \Lambda(\frac{1}{c_0},|h|_{Y^{[\frac{m}{2}],1}} + \|\nabla v\|_{Y^{[\frac{m}{2}],0}}) ( \varepsilon\|\nabla v\|^2_{X^{m-1,1}} + \varepsilon\|v\|^2_{X^{m-1,2}}) dt \\
	& + \theta\int_0^T \Lambda(\frac{1}{c_0},|h|_{Y^{[\frac{m}{2}],1}} + \|\nabla v\|_{Y^{[\frac{m}{2}],0}}) |\nabla h|^2_{X^{m-1,\frac{3}{2}}} dt.
	\end{split}
	\end{equation*}
	\noindent\textbf{10) Estimate of} Dirichlet-Neumann operator term. The left hand side of (\ref{whole DN}),
	\begin{equation} \label{DN term 1}
	\Big(1+ |h|_{W^{1,\infty}}\Big)^{-2}\left|\frac{|\nabla|}{\Big(1+|\nabla|\Big)^{1/2}} V^b \right|_{L^2}^2 \leq \Big(G[h]V^b,V^b\Big).
	\end{equation}
	Note that, ($\mathcal{F}$ means Fourier Transform with respect to horizontal direction.)
	\begin{equation} \label{DN term 2}
	\begin{split}
	\left| \frac{|\nabla|}{ (1+|\nabla| )^{1/2}} V^b \right|_{L^2(\partial S)} 
	&\geq \left| \frac{|\xi|^2}{ (1+|\xi| )^{1/2}} \mathcal{F}\Big( Z^{m-1} \frac{\nabla h}{\sqrt{1+|\nabla h|^2}} \Big) \right|_{L^2(\partial S)}   \\
	&= \Big| \Big( \frac{1+|\xi|^2}{ (1+|\xi| )^{1/2}} - \frac{1}{ (1+|\xi| )^{1/2}} \Big) \mathcal{F}\Big( Z^{m-1} \frac{\nabla h}{\sqrt{1+|\nabla h|^2}} \Big)  \Big|_{L^2(\partial S)} \\
	&\gtrsim \Big| Z^{m-1} \frac{\nabla h}{\sqrt{1+|\nabla h|^2}} \Big|_{H^{\frac{3}{2}}(\partial S)} - \Big| Z^{m-1} \frac{\nabla h}{\sqrt{1+|\nabla h|^2}} \Big|_{H^{-\frac{1}{2}}(\partial S)}.
	\end{split}
	\end{equation}
	Now, let $\Lambda^{s}$ be horizontal Fourier multiplier by $(1+|\xi|^2)^{s/2}$. We split $\Big| Z^{m-1} \frac{\nabla h}{\sqrt{1+|\nabla h|^2}} \Big|_{H^{\frac{3}{2}}(\partial S)}$ into highest order term and other low order terms. Then we can rewrite the first term in the RHS as
	\begin{equation} \label{DN term 3}
	\begin{split}
	&\Big| Z^{m-1} \frac{\nabla h}{\sqrt{1+|\nabla h|^2}} \Big|_{H^{\frac{3}{2}}(\partial S)} \\ &\geq \Big| \frac{Z^{m-1}\nabla h}{\sqrt{1 + |\nabla h|^2}} - \frac{\nabla\langle \nabla h, Z^{m-1}\nabla h \rangle}{\sqrt{1 + |\nabla h|^2}^3} \Big|_{H^{\frac{3}{2}}(\p S)} - |C^{m-1}(S)|_{H^{\frac{3}{2}}(\p S)}  \\
	&\geq \Big| \frac{\Lambda^{3/2} Z^{m-1}\nabla h}{\sqrt{1 + |\nabla h|^2}} - \frac{\nabla\langle \nabla h, \Lambda^{3/2} Z^{m-1}\nabla h \rangle}{\sqrt{1 + |\nabla h|^2}^3} \Big|_{L^{2}(\p S)} \\
	&\quad - \Lambda(\frac{1}{c_0}, |h|_{Y^{[\frac{m}{2}],1}}) |Z^{m-1}\nabla  h|_{H^1} - |C^{m-1}(S)|_{H^{\frac{3}{2}}(\p S)}  \\
	&\geq \Lambda(c_0) \Big\{ \int_{\p S} \Big| (1+|\nabla h|^2)|\Lambda^{3/2}Z^{m-1}\nabla h| - |\nabla h|^2|\Lambda^{3/2}Z^{m-1}\nabla h| \Big| \dd x \dd y \Big\}   \\
	&\quad - \Lambda(\frac{1}{c_0}, |h|_{Y^{[\frac{m}{2}],1}}) |Z^{m-1}\nabla  h|_{H^1}  - |C^{m-1}(S)|_{H^{\frac{3}{2}}(\p S)}  \\
	&\geq \Lambda(c_0) |Z^{m-1} \nabla h|_{\frac{3}{2}} - \Lambda(\frac{1}{c_0}, |h|_{Y^{[\frac{m}{2}],1}}) |Z^{m-1}\nabla  h|_{H^1} - |C^{m-1}(S)|_{H^{\frac{3}{2}}(\p S)} . \\
	\end{split}
	\end{equation}
	From (\ref{DN term 1}), (\ref{DN term 2}), and (\ref{DN term 3}), we get
	\begin{equation} \label{7.6}
	\int_0^T |Z^{m-1} \nabla h|^2_{\frac{3}{2}} dt \lesssim \sup_{t\in[0,T]}\Lambda(\frac{1}{c_0}, |h(t)|_{Y^{[\frac{m}{2}],1}}) \int_0^T \Big( (G[h]V^b,V^b) + |h|^2_{X^{m,1}} \Big) dt.
	\end{equation}
	
	\noindent We apply estimates \textbf{1) - 9)} to (\ref{whole DN}) with sufficiently small $\theta$ (from Young's inequaltiy) and Cauchy inequality to finish the proof.
\end{proof}
In the next section, we estimate for $Z^m = \partial_t^m$ case. By summing with above estimate, we get the estimate for norm $\left\|\cdot\right\|_{X^{m,0}}$.\\

\section{Energy Estimate for $Z^m = \p_t^m$ case}
In this section, we treat $Z^m = \partial_t^m$ case. If we apply Proposition \ref{energy estimate}, we will see $|\p_t^m h|_{\frac{3}{2}}$ on the right hand side of energy estimate. However, from elliptic problem of pressure, $\|\p_{t}^{m}q\|$ is not controllable. Therefore, we perform integration by part for time and $z$ to change $\p_{t}^{m}q$ into $\p_{t}^{m-1}\p_{z} q$. Boundary integral from these integration by parts vanish by cancellation. To see this, we will investigate two terms about pressure and surface tension, which generate high order commutators. We start with a simple proposition in the case of $Z^m = \partial_t^m$. 


\begin{proposition}
When $Z^{m}=\partial_t^m$, $C^m(f)$ can be estimated as follow.
\begin{equation*}
\left\|C^m(f)\right\| \leq \Lambda\big(\frac{1}{c_0},|h(s)|_{Y^{[\frac{m}{2}],2}} + \|\nabla f(s)\|_{Y^{[\frac{m}{2}],0}} \big) (|h|_{X^{m,\frac{1}{2}}} + \|\p_t^{m-1}\p_z f\|)
\end{equation*}
\end{proposition}
\begin{proof}
Since, $\partial_t$ commutes with $\partial_z$, we get the following.
\begin{equation*}
\partial_t^m (\partial_i^\varphi f ) = \partial_i^\varphi (\partial_t^m f ) + C_i^m(f),
\end{equation*}
$$
C_i^m(f) = -\left[\partial_t^m,\frac{\partial_i\varphi}{\partial_z\varphi},\partial_z f\right] - \Big(\partial_t^m\frac{\partial_i\varphi}{\partial_z\varphi}\Big)\partial_z f.
$$
From Proposition \ref{commutator est},
\begin{equation*}
\begin{split}
\Big\| \left[\partial_t^m,\frac{\partial_i\varphi}{\partial_z\varphi},\partial_z f\right] \Big\| &\leq \Big\|\partial_t^{m-1}\frac{\partial_i\varphi}{\partial_z\varphi}\Big\| \Big\|\partial_t^{[\frac{m}{2}]}\partial_z f\Big\| + \Big\|\partial_t^{[\frac{m}{2}]}\frac{\partial_i\varphi}{\partial_z\varphi}\Big\|\Big\|\partial_t^{m-1}\partial_z f\Big\|   \\
&\leq \Lambda\big(\frac{1}{c_0},|h(s)|_{Y^{[\frac{m}{2}],2}} + \|\nabla f(s)\|_{Y^{[\frac{m}{2}],0}} \big) ( \left|\partial_t^{m-1} h\right|_{\frac{1}{2}} + \left\|\partial_t^{m-1}\partial_z f\right\| ), \\
\Big\|\Big(\partial_t^m\frac{\partial_i\varphi}{\partial_z\varphi}\Big)\partial_z f\Big\| &\leq
\|\partial_z f\|_{L^\infty} \Big\|\partial_t^m \frac{\partial_i\varphi}{\partial_z\varphi}\Big\| \leq \Lambda(\frac{1}{c_0}, |h|_{Y^{[\frac{m}{2}],2}} ) \|\partial_z f\|_{L^\infty} |h|_{X^{m,\frac{1}{2}}}.
\end{split}
\end{equation*}
We use Lemma \ref{fraction est} in the last step.
\end{proof}

\begin{proposition} \label{alltime}
Assume that 
Let $m\geq 6$ and $(v,\varphi,q,h)$ is a smooth solution on $[0,T]$ of the system (\ref{1.12})-(\ref{1.15}). Also assume that
\begin{equation}
\partial_z\varphi \geq c_0,\quad \forall t\in [0,T].
\end{equation}
For $Z^{m}=\partial_t^m $,  we have the following energy estimate. 
\begin{equation}
\begin{split}
& \left\|\partial_t^m v(T)\right\|^2_{L^2(S)} + \left|\partial_t^m h(T)\right|^2_{H^1} + \varepsilon\int_{0}^{T} \left\|\nabla\partial_t^m v(t)\right\|^2_{L^2(S)} dt  \\
&\leq \Lambda_{m,\infty}(0)\big( \|v(0)\|^{2}_{X^{m-1,1}} + |h(0)|^{2}_{X^{m-1,2}} + \varepsilon\|\nabla v(0)\|^{2}_{X^{m-1,1}} \big) \\
&\quad + \Lambda_{m,\infty}(T)\big( \|v(T)\|^{2}_{X^{m-1,1}} + |h(T)|^{2}_{X^{m-1,2}} \big) \\
&\quad + \int_{0}^{T} \Lambda_{m,\infty}(s) \big( \|v\|^{2}_{X^{m,0}} + \|\nabla v\|^{2}_{X^{m-1,0}} + |h|^{2}_{X^{m-1,2}} + \varepsilon\|\nabla v\|^{2}_{X^{m-1,1}} + |\p_{t}^{m-1} h|^{2}_{H^{\frac{5}{2}}} \big) dt. \\
\end{split}
\end{equation}
\end{proposition}
\begin{proof}
From the analysis of previous section, there are two terms which generate highest order terms. \textbf{First term} is $-\int_S C^m(d) \partial_t^m q \partial_z\varphi \dd x \dd y \dd z$ which corresponds to (\ref{RC est}). The \textbf{Second term} is boundary integral 
\[
-\int_{\partial S}\partial_t^m q^{b} C^m(KB)  
\]
of (\ref{P11 est}). Except these two terms, we do not see any high order commutator $|\p_t^m h|_{\frac{3}{2}}$. Therefore, similar as in proof of Proposition \ref{high est}, we have
\begin{equation} \label{all time est form}
\begin{split}
& \frac{d}{dt}\int_S\left|\p_t^m v\right|^2dV_t
+ 4\varepsilon\int_S\left|S^\varphi (\p_t^m v)\right|^2dV_t + g \frac{\dd}{\dd t} \int_{\p S} |\p_t^m h|^2   + \frac{\dd}{\dd t} \int_{\p S} \frac{|\p_t^m \nabla h|^2}{\sqrt{1 + |\nabla h|^2}}  \\
& \leq \Lambda\big(\frac{1}{c_0},\left|h\right|_{Y^{[\frac{m}{2}],1}} + \left\|\nabla v\right\|_{Y^{[\frac{m}{2}],0}} + \left\|\nabla q\right\|_{Y^{[\frac{m}{2}],0}}\big) \Big( \left\|v\right\|^2_{X^{m,0}} + \left\|\nabla v\right\|^2_{X^{m-1,0}} + \varepsilon \|\nabla v \|^2_{X^{m,0}} + |h|^2_{X^{m,1}} \Big) \\
& + \Lambda(\frac{1}{c_0}) \Big( \underbrace{ -\int_S C^m(d) \partial_t^m q \partial_z\varphi  }_{\textbf{First term}} + \underbrace{ \int_{\partial S}\partial_t^m q^{b} C^m(KB)   }_{\textbf{Second term}} \Big) . 
\end{split}
\end{equation}
Let us analyze \textbf{First term}. Using divergence free condition, we can expand $\partial_z\varphi C^m(d)$ as follow.
\begin{equation} \label{Cm expan}
\begin{split}
\partial_z\varphi C^m(d) & = [\partial_t^m,\mathbf{N},\cdot\partial_z v] + [\partial_t^m,\partial_z\eta,\partial_1 v_1 + \partial_2 v_2] \\
& := C^m(d)_1 + C^m(d)_2 + C^m(d)_3 + C^m(d)_4 + C^m(d)_5,
\end{split}
\end{equation}
where
\begin{eqnarray}
\begin{split}
C^m(d)_1 & := m\partial_t \mathbf{N}\cdot\partial_t^{m-1}\partial_z v, \\
C^m(d)_2 & := m\partial_t\partial_z\eta\partial_t^{m-1}(\partial_1 v_1 + \partial_2 v_2), \\
C^m(d)_3 & := m\partial_t^{m-1} \mathbf{N}\cdot\partial_t\partial_z v, \\
C^m(d)_4 & := m\partial_t^{m-1}\partial_z\eta\partial_t(\partial_1 v_1 + \partial_2 v_2), \\
C^m(d)_5 & := \sum_{l=2}^{m-2}C_{m}^{l}\Big(\partial_t^l \mathbf{N}\cdot\partial_t^{m-l}\partial_z v + \partial_t^l\partial_z \eta\cdot\partial_t^{m-l}(\partial_1 v_1 + \partial_2 v_2) \Big).
\end{split}
\end{eqnarray}
We use (\ref{Cm expan}) to integrand $C^{m}(d)\p_{z}\varphi$ for $i=1,2,3,4,5$. Since $\|\p_{t}^{m}q\|$ cannot be controlled, we should perform integration by part in times. And then $\p_{t}C^{m}(d)_{1}, \p_{t}C^{m}(d)_{2}\  \sim \p_{t}^{m}\nabla v$, so we should also perform integration by part in space.  \\

\noindent\textbf{1) Estimate of} $-\int_{0}^{T} \int_{S} C^m(d)_1 \partial_t^m q ds$. We perfrom integration by part in $z$ and also in $t$ again. When we split $q = q^{E} + q^{NS} + q^{S}$, $q^{E}$ and $q^{NS}$ parts are okay, since their estimates contain only low order terms. The only problem is estimate including $q^{S}$. Inspired by harmonic extension, we consider horizontal fourier multiplier $|\nabla_{y}|$ to estimate terms with $q^{S}$. Using pressure estimates (\ref{q_E est}), (\ref{q_NS est}), and $\varepsilon \ll 1$, we get

\begin{equation} \label{al t expan form}
\begin{split}
&-\int_{0}^{T} \int_{S} C^m(d)_1 \partial_t^m q dt = -\int_{0}^{T} \int_S m\p_t \mathbf{N}\cdot\partial_t^{m-1}\partial_z v \partial_t^m q dt. \\
&= -\int_{0}^{T} \int_{\partial S} m\partial_t \mathbf{N}\partial_t^{m-1} v^b \partial_t^m q^b dt + \int_{0}^{T} \int_S m\partial_z\partial_t \mathbf{N}\partial_t^{m-1} v\partial_t^m q dt + \int_{0}^{T} \int_S m\partial_t \mathbf{N}\partial_t^{m-1} v\partial_t^m\partial_z q dt \\
&\leq -  \int_{0}^{T} \int_{\partial S} m\partial_t \mathbf{N}\partial_t^{m-1} v^b \partial_t^m q^b dt \\
&\quad + \Lambda_{m,\infty}(0) \|\p_{t}^{m-1}v(0) \| \| \p_{t}^{m-1}\p_{z} q^{E}(0) \| + \Lambda_{m,\infty}(T) \|\p_{t}^{m-1}v(T) \| \| |\p_{t}^{m-1}\p_{z}q^{E}(T) \| \\
&\quad + \Lambda_{m,\infty}(0) \| |\nabla_{y}|\p_{t}^{m-1}v(0) \| \| |\nabla_{y}|^{-1} \p_{t}^{m-1}\p_{z}(q^{S} + q^{NS})(0) \|  \\
&\quad + \Lambda_{m,\infty}(T) \| |\nabla_{y}|\p_{t}^{m-1}v(T) \| \| |\nabla_{y}|^{-1}\p_{t}^{m-1}\p_{z}(q^{S} + q^{NS})(T) \| \\
&\quad + \int_{0}^{T} \Lambda_{m,\infty}(t) \big( \|v\|^{2}_{X^{m,0}} + \|\nabla v\|^{2}_{X^{m-1,0}} + |h|^{2}_{X^{m-1,2}} + \varepsilon\|\nabla v\|^{2}_{X^{m-1,1}} + |\p_{t}^{m-1} h|^{2}_{H^{\frac{5}{2}}} \big) dt,  \\
\end{split}
\end{equation}
where we performed integration by parts in time in the last step. Using estimates of $q^{E}$, we have 
\begin{equation} \label{Cmd1-1}
\begin{split}
&\leq -  \int_{0}^{T} \int_{\partial S} m\partial_t \mathbf{N}\partial_t^{m-1} v^b \partial_t^m q^b dt \\
&\quad + \Lambda_{m,\infty}(0) \| |\nabla_{y}|\p_{t}^{m-1}v(0) \| \| |\nabla_{y}|^{-1} \p_{t}^{m-1}\p_{z}(q^{S} + q^{NS})(0) \|  \\
&\quad + \Lambda_{m,\infty}(T) \| |\nabla_{y}|\p_{t}^{m-1}v(T) \| \| |\nabla_{y}|^{-1}\p_{t}^{m-1}\p_{z}(q^{S} + q^{NS})(T) \| \\
&\quad + \Lambda_{m,\infty}(0)\big( \|v(0)\|^{2}_{X^{m-1,0}} + |h(0)|^{2}_{X^{m-1,2}} \big) + \Lambda_{m,\infty}(T)\big( \|v(T)\|^{2}_{X^{m-1,0}} + |h(t)|^{2}_{X^{m-1,2}} \big)   \\
&\quad + \int_{0}^{T} \Lambda_{m,\infty}(t) \big( \|v\|^{2}_{X^{m,0}} + \|\nabla v\|^{2}_{X^{m-1,0}} + |h|^{2}_{X^{m-1,2}} + \varepsilon\|\nabla v\|^{2}_{X^{m-1,1}} + |\p_{t}^{m-1} h|^{2}_{H^{\frac{5}{2}}} \big) dt \\
\end{split}
\end{equation}
It is easy to control $\||\nabla_{y}|\p_{t}^{m-1}v\| \leq \|v\|_{X^{m-1,1}}$ by definition. We treat $\| |\nabla_{y}|^{-1} \p_{t}^{m-1}\p_{z} q^{S}\|$ and $\| |\nabla_{y}|^{-1} \p_{t}^{m-1}\p_{z} q^{NS}\|$ using horizontal Fourier multiplier. Let us write $\hat{q}^{S} := (\Lambda_{y}q^{S})(\xi,z)$, where $\Lambda_{y}$ horizontal fourier multiplier. Note that $q^{NS}$ and $q^{S}$ are harmonic,  
\begin{equation}
\begin{split}
|\xi|^{2}\hat{q}^{S} &= \p_{zz} \hat{q}^{S}, \\
\hat{q}^{S} &= \hat{q}^{S}(\xi,\hat{h}) e^{|\xi|z},  \\
\p_{z} \hat{q}^{S} &= |\xi| \mathcal{F}\Big( {\nabla\cdot\frac{\nabla h}{\sqrt{1+|\nabla h|^{2}}}} \Big) e^{|\xi| z} \leq \Lambda_{m,\infty}|\xi|^{3} \hat{h} e^{|\xi| z},
\end{split}
\end{equation}	
where $\mathcal{F}$ means horizontal Fourier transform. This implies $\p_{z}$ is changed into horizontal derivatives in the case of harmonic function. And we do not need full Fourier multiplier $|\nabla|$ to reduce order and
\begin{equation} \label{fourier q S}
\begin{split}
	\||\nabla_{y}|^{-1} \p_{t}^{m-1}\p_{z} {q}^{S}\|^{2} &= \||\nabla_{y}|^{-1} \p_{t}^{m-1}\p_{z} \hat{q}^{S}\|^{2} \\
	&= \int_{-\infty}^{0}\int_{\mathbb{R}^{2}} \frac{1}{|\xi|^{2}} |\xi|^{6}|\p_{t}^{m-1}\hat{q}^{S}|^{2} e^{2|\xi| z}  \\
	&\leq \Lambda_{m,\infty} \int_{\mathbb{R}^{2}} \frac{1}{2\xi|\xi|^{2}} |\xi|^{6}|\p_{t}^{m-1}\hat{h}|^{2}  \\
	&\leq \Lambda_{m,\infty} \| |\xi|^{3/2}\p_{t}^{m-1}\hat{h}\|^{2} \leq \Lambda_{m,\infty} |h|_{X^{m-1,\frac{3}{2}}}^{2}.  \\
\end{split}
\end{equation}

We treat $q^{NS}$ similar as $q^{S}$ case since $q^{NS}$ also solves harmonic equation. First, $\||\nabla_{y}|\p_{t}^{m-1}v\| \leq \|v\|_{X^{m-1,1}}$. Horizontal operator $|\nabla_{y}|^{-1}$ reduce order by $1$ of boundary data. For $t=0$, we get
\begin{equation} \label{fourier q NS zero}
\begin{split}
\||\nabla_{y}|^{-1} \p_{t}^{m-1}\p_{z} {q}^{NS}(0)\|
&\leq \varepsilon \Lambda_{m,\infty}(0) \big( |h(0)|_{X^{m-1,2}} + \|\nabla \p_{t}^{m-1} v(0)\| \big) . \\
\end{split}
\end{equation}
For $t=T$, since the $\varepsilon\|\nabla \p_{t}^{m-1}v\|$ is dissipation type, we derive time integration by 
\begin{equation} \label{fourier q NS T}
\begin{split}
\||\nabla_{y}|^{-1} \p_{t}^{m-1}\p_{z} {q}^{NS}(T)\|
&\leq \varepsilon \Lambda_{m,\infty}(T) \big( |h(T)|_{X^{m-1,2}} + \|\nabla \p_{t}^{m-1} v(T)\| \big)  \\
&\leq \varepsilon \Lambda_{m,\infty}(T) \big( |h(T)|_{X^{m-1,2}} +  \|\nabla \p_{t}^{m-1}v(0)\| + \int_{0}^{T} \|\nabla \p_{t}^{m} v(t)\| dt \big).  \\
\end{split}
\end{equation}
Therefore,
\begin{equation} \label{q NS usage}
\begin{split}
&\Lambda_{m,\infty}(T) \| |\nabla_{y}|^{1} \p_{t}^{m-1}v(T) \| \| |\nabla_{y}|^{-1}\p_{t}^{m-1}\p_{z}q^{NS}(T) \| \\
&\leq \varepsilon\Lambda_{m,\infty}(T) \|\p_{t}^{m-1}v(T) \|_{X^{0,1}} \big( |h(T)|_{X^{m-1,2}} +  \|\nabla \p_{t}^{m-1}v(0)\| + \int_{0}^{T} \|\nabla \p_{t}^{m} v(t)\| dt \big) \\
&\leq \varepsilon\|\nabla \p_{t}^{m-1}v(0)\|^{2} + \varepsilon\Lambda_{m,\infty}(T) \big( \|v(T)\|^{2}_{X^{m-1,1}} + |h(T)|^{2}_{X^{m-1,2}} \big) + T\int_{0}^{T} \|\nabla \p_{t}^{m} v(t)\|^{2} dt .
\end{split}
\end{equation}

We apply (\ref{fourier q S}), (\ref{fourier q NS zero}), (\ref{fourier q NS T}), and (\ref{q NS usage}) to (\ref{Cmd1-1}) to gain
\begin{equation} \label{Cmd1}
\begin{split}
(\ref{Cmd1}) &\leq - \underbrace{ \int_{0}^{T}\int_{\partial S} m\partial_t \mathbf{N}\partial_t^{m-1} v^b \partial_t^m q^b dt }_{(*)}  \\
&\quad + \Lambda_{m,\infty}(0)\big( \|v(0)\|^{2}_{X^{m-1,1}} + |h(0)|^{2}_{X^{m-1,2}} + \varepsilon\|\nabla v(0)\|^{2}_{X^{m-1,1}} \big) \\
&\quad + \Lambda_{m,\infty}(T)\big( \|v(T)\|^{2}_{X^{m-1,1}} + |h(T)|^{2}_{X^{m-1,2}} \big) + T\varepsilon \int_{0}^{T} \|\nabla \p_{t}^{m} v(t)\|^{2} dt   \\
&\quad + \int_{0}^{T} \Lambda_{m,\infty}(t) \big( \|v\|^{2}_{X^{m,0}} + \|\nabla v\|^{2}_{X^{m-1,0}} + |h|^{2}_{X^{m-1,2}} + \varepsilon\|\nabla v\|^{2}_{X^{m-1,1}} + |\p_{t}^{m-1} h|^{2}_{H^{\frac{5}{2}}} \big) dt.  \\
\end{split}
\end{equation}	

\noindent\textbf{2) Estimate of} $-\int_{0}^{T} \int_{S} C^m(d)_2 \partial_t^m q dt$. This estimate is very similar as above $C^{m}(d)_{1}$ case. However, $C^m(d)_2$ contains horizontal derivatives $\p_{1}v, \ \p_{2}v$, instead of $\p_{z}v$, so we perform integration by part in $t$ and horizontal $x$ or $y$. Also we treat pressure term similar as we did in (\ref{Cmd1-1}) and (\ref{Cmd1}). Therefore we get same estimate with (\ref{Cmd1}) but without $(*)$.
\begin{equation} \label{Cmd2}
\begin{split}
&-\int_{0}^{T} \int_{S} C^m(d)_{2} \partial_t^m q dt  \\
&\leq + \Lambda_{m,\infty}(0)\big( \|v(0)\|^{2}_{X^{m-1,1}} + |h(0)|^{2}_{X^{m-1,2}} + \varepsilon\|\nabla v(0)\|^{2}_{X^{m-1,1}} \big) \\
&\quad + \Lambda_{m,\infty}(T)\big( \|v(T)\|^{2}_{X^{m-1,1}} + |h(T)|^{2}_{X^{m-1,2}} \big) + T\varepsilon \int_{0}^{T} \|\nabla \p_{t}^{m} v(t)\|^{2} dt   \\
&\quad + \int_{0}^{T} \Lambda_{m,\infty}(s) \big( \|v\|^{2}_{X^{m,0}} + \|\nabla v\|^{2}_{X^{m-1,0}} + |h|^{2}_{X^{m-1,2}} + \varepsilon\|\nabla v\|^{2}_{X^{m-1,1}} + |\p_{t}^{m-1} h|^{2}_{H^{\frac{5}{2}}} \big) dt. \\
\end{split}
\end{equation}

\noindent \textbf{3) Estimate of} $-\int_{0}^{T} \int_{S} C^m(d)_{i}\partial_t^m q dt$ for $i=3,4,5$. Since $\eta$ is $\frac{1}{2}$ better regularity than $h$, we just perform integration by part in time to derive similar estimate as (\ref{Cmd2}). Of course we treat pressure similar way as we did in (\ref{Cmd1-1}). For $i=3,4,5$,
\begin{equation}
\begin{split}
&-\int_0^t \int_{S} C^m(d)_{i} \partial_t^m q dt \\
&\leq \Lambda_{m,\infty}(0)\big( \|v(0)\|^{2}_{X^{m-1,1}} + |h(0)|^{2}_{X^{m-1,2}} + \varepsilon\|\nabla v(0)\|^{2}_{X^{m-1,1}} \big) \\
&\quad + \Lambda_{m,\infty}(T)\big( \|v(T)\|^{2}_{X^{m-1,1}} + |h(T)|^{2}_{X^{m-1,2}} \big) + T\varepsilon \int_{0}^{T} \|\nabla \p_{t}^{m} v(t)\|^{2} dt   \\
&\quad + \int_{0}^{T} \Lambda_{m,\infty}(s) \big( \|v\|^{2}_{X^{m,0}} + \|\nabla v\|^{2}_{X^{m-1,0}} + |h|^{2}_{X^{m-1,2}} + \varepsilon\|\nabla v\|^{2}_{X^{m-1,1}} + |\p_{t}^{m-1} h|^{2}_{H^{\frac{5}{2}}} \big) dt. \\
\end{split}
\end{equation}
 
Now let us analyze underbraced \textbf{Second term} in (\ref{all time est form}). From (\ref{331}),
\begin{equation} 
\begin{split}
	C^{m}(KB) &= [Z^{m}, v^{b}, \mathbf{N}] \\
	&= m \p_{t}\mathbf{N} \p_{t}^{m-1}v^{b} + m \p_{t}^{m-1}\mathbf{N} \p_{t}v^{b} + \sum_{\ell=2}^{m-2} \beta_{\ell} \p_{t}^{\ell}\mathbf{N}\p_{t}^{m-\ell} v^{b},
\end{split}
\end{equation}
where $\beta_{\ell}$ is some combinatoric positive integers.
\begin{equation} \label{Cmd345}
\begin{split}
& \int_{0}^{T} \int_{\partial S}\partial_t^m q^b C^m(KB)  dt \\
&= \int_{0}^{T} \int_{\partial S}\partial_t^m q^b \Big( m \p_{t}\mathbf{N} \p_{t}^{m-1}v^{b} + m \p_{t}^{m-1}\mathbf{N} \p_{t}v^{b} + \sum_{\ell=2}^{m-2} \beta_{\ell} \p_{t}^{\ell}\mathbf{N}\p_{t}^{m-\ell} v^{b} \Big)  dt \\
&= \underbrace{ \int_{0}^{T} \int_{\partial S}\partial_t^m q^b m \partial_t \mathbf{N} \partial_t^{m-1}v^b dt }_{(*)}  \\
&\quad + \int_{0}^{T} \int_{\partial S}\partial_t^m q^b \Big( m \p_{t}^{m-1}\mathbf{N} \p_{t}v^{b} + \sum_{\ell=2}^{m-2} \beta_{\ell} \p_{t}^{\ell}\mathbf{N}\p_{t}^{m-\ell} v^{b} \Big)  dt . \\
\end{split}
\end{equation}

For the first term in the last line, we perform integration by parts in $z$ and $t$ as we did in estimate of \textbf{First term}. We use simple estimate 
\[
	\int_{\p S} \p_{t}^{m-1}\p_{z}q^{b} \p_{t}^{m-1}\mathbf{N}\p_{t}v^{b} \leq |\p_{t}^{m-1}\mathbf{N}\p_{t}v^{b}|_{\frac{1}{2}} |\p_{t}^{m-1}\p_{z}q^{b}|_{-\frac{1}{2}},
\]
and use Proposition \ref{aniso embed}, \ref{horizon est} and pressure estimates (\ref{q_E est}), (\ref{q_NS est}), (\ref{q_S est}) to gain
\begin{equation}
\begin{split}
	& \int_{0}^{T} \int_{\partial S}\partial_t^m q^b \p_{t}^{m-1}\mathbf{N} \p_{t}v^{b}  \\
	&\leq \Lambda_{m,\infty}(0) |\p_{t}^{m-1}h(0)|_{\frac{3}{2}} \| \p_{t}^{m-1}\nabla q^{E}(0) \| + \Lambda_{m,\infty}(T) |\p_{t}^{m-1}h(T)|_{\frac{3}{2}} \| |\p_{t}^{m-1}\nabla q^{E}(T) \| \\
	&\quad + \Lambda_{m,\infty}(0) \| |\nabla_{y}|\p_{t}^{m-1}h(0) \| \| |\nabla_{y}|^{-1} \p_{t}^{m-1}\p_{z}(q^{S} + q^{NS})(0) \|  \\
	&\quad + \Lambda_{m,\infty}(T) \| |\nabla_{y}|\p_{t}^{m-1}h(T) \| \| |\nabla_{y}|^{-1}\p_{t}^{m-1}\p_{z}(q^{S} + q^{NS})(T) \| \\
	&\quad + \int_{0}^{T} \Lambda_{m,\infty}(t) \big( \|v\|^{2}_{X^{m,0}} + \|\nabla v\|^{2}_{X^{m-1,0}} + |h|^{2}_{X^{m-1,2}} + \varepsilon\|\nabla v\|^{2}_{X^{m-1,1}} + |\p_{t}^{m-1} h|^{2}_{H^{\frac{5}{2}}} \big) dt.  \\
\end{split}
\end{equation}
This form is nearly same as (\ref{al t expan form}), except $v(0), v(T)$ are changed into $h(0), h(T)$. So we can use same estimates of (\ref{Cmd1}).  \\

Estimating the second term in the last line of (\ref{Cmd345}) is exactly same as above because $|\sum_{\ell=2}^{m-2} \beta_{\ell} \p_{t}^{\ell}\mathbf{N}\p_{t}^{m-\ell} v^{b}|_{\frac{1}{2}}$ does not include high order bad commutators. Finally,  \\

\begin{equation} \label{2nd term est} 
\begin{split}
& (\ref{Cmd345}) \leq \underbrace{ \int_{0}^{T} \int_{\partial S}\partial_t^m q^b m \partial_t \mathbf{N} \partial_t^{m-1}v^b dt }_{(*)}  \\
&\quad + \Lambda_{m,\infty}(0)\big( \|v(0)\|^{2}_{X^{m-1,1}} + |h(0)|^{2}_{X^{m-1,2}} + \varepsilon\|\nabla v(0)\|^{2}_{X^{m-1,1}} \big) \\
&\quad + \Lambda_{m,\infty}(T)\big( \|v(T)\|^{2}_{X^{m-1,1}} + |h(T)|^{2}_{X^{m-1,2}} \big) + T\varepsilon \int_{0}^{T} \|\nabla \p_{t}^{m} v(t)\|^{2} dt   \\
&\quad + \int_{0}^{T} \Lambda_{m,\infty}(t) \big( \|v\|^{2}_{X^{m,0}} + \|\nabla v\|^{2}_{X^{m-1,0}} + |h|^{2}_{X^{m-1,2}} + \varepsilon\|\nabla v\|^{2}_{X^{m-1,1}} + |\p_{t}^{m-1} h|^{2}_{H^{\frac{5}{2}}} \big) dt.  \\
\end{split}
\end{equation}    
 \\
Now we apply (\ref{Cmd1}), (\ref{Cmd2}), (\ref{Cmd345}), (\ref{2nd term est}), and (\ref{q_infty est}) to (\ref{all time est form}) to get estimate. Especially uncontrollable boundary integral $\int_{0}^{T} \int_{\partial S}\partial_t^m q^b m \partial_t \mathbf{N} \partial_t^{m-1}v^b dt$ vanishes from $(*)$ terms in (\ref{Cmd1}) and (\ref{2nd term est}) with opposite sign. Dissipation type terms $T\varepsilon\int_{0}^{T}\|\nabla\p_{t}^{m}v(t)\|^{2} dt$ in (\ref{Cmd1}), (\ref{Cmd2}), and (\ref{Cmd345}) are absorbed by dissipation $\varepsilon\int_{0}^{T}\|\nabla\p_{t}^{m}v(t)\|^{2} dt$ which is in the left hand side of (\ref{all time est form}) for sufficiently small $T \ll 1$. This finishes proof.
\end{proof}

\section{Normal derivative estimate}
From Proposition \ref{energy estimate}, we should control $\|\partial_z v\|_{X^{m-1,0}}$, since $H^m \hookrightarrow H^{m}_{co}$. However, it is hard to estimate $\partial_z v$ directly. Instead, we estimate $S_n$, which is tangential part of $S^\varphi v \mathbf{n}$, i.e. 
\begin{equation} \label{def Sn}
S_n := \Pi(S^\varphi v \mathbf{n})\quad\text{where}\quad \Pi = I - \mathbf{n} \otimes \mathbf{n},
\end{equation} 
where $I$ is identity matrix and $\mathbf{n}$ is defined in Definition \ref{extendnormal}. First, we show that $S_n$ is equivalent to $\partial_z v$ in the space $X^{m-1,0}$. It is clear that 
\begin{eqnarray}
	\|S_n\|_{X^{m-1,0}} \leq \Lambda(\frac{1}{c_0}, |h|_{Y^{[\frac{m}{2}],1}} + \|\nabla v\|_{Y^{[\frac{m}{2}],0}}) ( |h|_{X^{m-1,1}} + \|v\|_{X^{m,0}} + \|\nabla v\|_{X^{m-1,0}} ),
\end{eqnarray}
from above definition (\ref{def Sn}) using Proposition \ref{commutator est} and \ref{fraction est}. The following two lemmas show how to control $\p_z v$ using $S_n$.
\begin{lemma} \label{pre normal control}
We have the following normal part estimate of $\partial_z v$.
\begin{equation}
\| \partial_z v\cdot \mathbf{n}  \|_{X^{m-1,0}} \leq
\Lambda(\frac{1}{c_0}, |h|_{Y^{[\frac{m}{2}],1}} + \|\nabla v\|_{Y^{[\frac{m}{2}],0}}) ( |h|_{X^{m-1,1}} + \|v\|_{X^{m-1,1}} ).
\end{equation}
\end{lemma}
\begin{proof}
From divergence free condition, we have,
\begin{equation*}
\partial_z v \cdot \mathbf{n} = \frac{1}{\sqrt{1+|\nabla_y\varphi|^2}}\partial_z\varphi (\partial_1 v_1 + \partial_2 v_2).
\end{equation*}
Applying $Z^{m-1}$ and using Proposition \ref{commutator est}, we get
\begin{equation*}
\|Z^{m-1} (\partial_z v\cdot \mathbf{n} ) \| \leq
\Lambda(\frac{1}{c_0}, |h|_{Y^{[\frac{m}{2}],1}} + \|\nabla v\|_{Y^{[\frac{m}{2}],0}}) ( |h|_{X^{m-1,1}} + \|v\|_{X^{m-1,1}} ).
\end{equation*}
\end{proof}

\begin{lemma} \label{normal control}
We have the following estimate.  
\begin{equation*}
\left\|\partial_z v\right\|_{X^{m-1,0}} \leq \Lambda(\frac{1}{c_0}, |h|_{Y^{[\frac{m}{2}],1}} + \|\nabla v\|_{Y^{[\frac{m}{2}],0}})
( \left\|v\right\|_{X^{m,0}} + \left|h\right|_{X^{m-1,1}} + \left\|S_n\right\|_{X^{m-1,0}}  ).
\end{equation*}
\end{lemma}
\begin{proof}
By definition,
\begin{equation*}
2S^\varphi v\mathbf{n} := (\nabla u)\mathbf{n} + (\nabla u)^T\mathbf{n} = (\nabla u)\mathbf{n} + g^{ij}(\partial_j v\cdot\mathbf{n})\partial_{y^i}.
\end{equation*}
And from divergence free condition, we get
\begin{equation} \label{div free normal}
\partial_N u = \frac{1+\left|\nabla_y\varphi\right|^2}{\partial_z \varphi}\partial_z v - \partial_1\varphi\partial_1 v - \partial_2\varphi\partial_2 v .
\end{equation}
Therefore,
\begin{equation*}
 \|\partial_z v \|_{X^{m-1,0}} \leq \Lambda(\frac{1}{c_0}, |h|_{Y^{[\frac{m}{2}],1}} + \|\nabla v\|_{Y^{[\frac{m}{2}],0}})
( \left\|v\right\|_{X^{m,0}} + \left|h\right|_{X^{m-1,1}} + \left\|\partial_z v\cdot\mathbf{n}\right\|_{X^{m-1,0}} + \left\|S^\varphi v\mathbf{n}\right\|_{X^{m,0}} )
\end{equation*}
and
\begin{equation*}
S^\varphi v \mathbf{n} = S_n +  (\mathbf{n}\otimes\mathbf{n} ) (S^\varphi v \mathbf{n} ).
\end{equation*}
We use Lemma \ref{pre normal control} to get
\begin{equation*}
\left\|\partial_z v\right\|_{X^{m-1,0}} \leq \Lambda(\frac{1}{c_0}, |h|_{Y^{[\frac{m}{2}],1}} + \|\nabla v\|_{Y^{[\frac{m}{2}],0}})
 ( \left\|v\right\|_{X^{m,0}} + \left|h\right|_{X^{m-1,1}} + \left\|S_n\right\|_{X^{m-1,0}}  ).
\end{equation*}
\end{proof}

To get an estimate for $S_n$ we make an equation for $S_n$. First, we take $\nabla^\varphi$ to (\ref{1.12}).
\begin{equation} \label{eq nabla}
\partial_t^\varphi \nabla^\varphi v + (v\cdot\nabla^\varphi)\nabla^\varphi v + (\nabla^\varphi v)^2 + (D^\varphi)^2 q -\varepsilon\Delta^\varphi\nabla^\varphi v =0,
\end{equation}
where $(D^\varphi)^2$ is Hessian matrix. We also take symmetric part of above (\ref{eq nabla}). Then by adding both equations, we get
\begin{equation*}
\partial_t^\varphi S^\varphi v +  (v\cdot\nabla^\varphi )S^\varphi v + \frac{1}{2} ( (\nabla v )^2 +  ( (\nabla v )^T )^2  ) + (D^\varphi)^2 q - \varepsilon\Delta^\varphi(S^\varphi v) =0.
\end{equation*}
By taking tangential operator, $\Pi$, which was defined in (\ref{def Sn}),
\begin{equation} \label{eq Sn}
\partial_t^\varphi S_n +  (v\cdot\nabla^\varphi )S_n - \varepsilon\Delta^\varphi (S_n ) = F_S,
\end{equation}
where $F_S$ is commutator,
\begin{equation} \label{def fs}
\begin{split}
F_S &= F^1_S + F^2_S + F^3_S,  \\
F^1_S &= -\frac{1}{2}\Pi ( (\nabla^\varphi v )^2 +  ( (\nabla^\varphi v )^T )^2 )\mathbf{n} +  (\partial_t\Pi + v\cdot\nabla^\varphi\Pi )S^\varphi v\mathbf{n} + \Pi S^\varphi v (\partial_t \mathbf{n} + v\cdot\nabla^\varphi \mathbf{n} ), \\
F^2_S &= -2\varepsilon\partial_i^\varphi\Pi\partial_i^\varphi (S^\varphi v\mathbf{n} ) - 2\varepsilon\Pi (\partial_i^\varphi (S^\varphi v)\partial_i^\varphi \mathbf{n} ) - 
\varepsilon (\Delta^\varphi\Pi)S^\varphi v\mathbf{n} - \varepsilon\Pi S^\varphi v \Delta^\varphi\mathbf{n} ,  \\
F^3_S &= - \Pi((D^\varphi)^2 q )\mathbf{n}.
\end{split}
\end{equation}
We apply $Z^{k}$ to get higher order estimate. We need $\|S_n\|_{X^{m-1,0}}$, however, from pressure estimate, optimal regularity is $k=m-2$. This is because,
\[
\|F_S^3\|_{m-2} \sim \| D^2 q \|_{m-2} \sim |h|_{X^{m-2,\frac{7}{2}}}.
\]
To estimateFor $Z^{m-2}F^1_S$, using Propositions \ref{commutator est},
\begin{equation} \label{FS1}
\begin{split}
\left\|Z^{m-2} F^1_S\right\|_{L^2(S)} &\leq \Lambda(\frac{1}{c_0},|h|_{Y^{[\frac{m}{2}],1}} + \|\nabla v\|_{Y^{[\frac{m}{2}],0}} )  ( \left\|\nabla v\right\|_{X^{m-2,0}} + \left|h\right|_{X^{m-2,1}} + \left\|v\right\|_{X^{m-2,0}} )  \\
&\leq \Lambda_{m,\infty}  ( \left\|S_n\right\|_{X^{m-2,0}} + \left|h\right|_{X^{m,1}} + \left\|v\right\|_{X^{m-2,0}} ).
\end{split}
\end{equation}
where we used Lemma \ref{normal control}.  \\
Similary, for $Z^{m-2}F^2_S$,
\begin{equation} \label{FS2}
\begin{split}
\left\|Z^{m-2} F^2_S\right\|_{L^2(S)} &\leq \varepsilon \Lambda(\frac{1}{c_0},|h|_{Y^{[\frac{m}{2}],1}} + \|\nabla v\|_{Y^{[\frac{m}{2}],0}} ) ( \left\|\partial_{zz} v\right\|_{X^{m-2,0}} + \left\|\partial_{z} v\right\|_{X^{m-1,0}} + \left|h\right|_{X^{m-2,\frac{5}{2}}} ) \\
&\leq \varepsilon \Lambda_{m,\infty} ( \left\|\nabla  S_n\right\|_{X^{m-2,0}} + \left\|S_n\right\|_{X^{m-1,0}} + \left|h\right|_{X^{m-2,\frac{5}{2}}} ) .
\end{split}
\end{equation}
\\
For $Z^{m-2}F^3_S$, using pressure estimates Proposition \ref{q_E est}, \ref{q_NS est}, and \ref{q_S est}, and Lemma \ref{normal control},
\begin{equation} \label{FS3}
\left\|Z^{m-2} F^3_S\right\|_{L^2(S)} \leq \Lambda(\frac{1}{c_0},|h|_{Y^{[\frac{m}{2}],1}} + \|\nabla v\|_{Y^{[\frac{m}{2}],0}} ) (  \|S_n\|_{X^{m-2,0}} + |h|_{X^{m-2,\frac{7}{2}}} + \|v\|_{X^{m-2,0}} ).
\end{equation}
Now we apply $Z^{m-2}$ to (\ref{eq Sn}) to get,
\begin{equation*}
\partial_t^\varphi Z^{m-2} S_n + (v\cdot\nabla^\varphi)Z^{m-2} S_n - \varepsilon\Delta^\varphi Z^{m-2} S_n = Z^{m-2} (F_S) + C_S ,
\end{equation*}
where $C_S$ is commutator. We divide $C_S$ into,
\begin{equation}
C_S^1 = [Z^\alpha v_y]\cdot\nabla_y S_n + [Z^\alpha,V_z]\partial_z S_n := C_{S_y} + C_{S_z} \,\,\, , \,\,\,\,\,\,\,C_S^2 = -\varepsilon[Z^\alpha,\Delta^\varphi]S_n,
\end{equation}
where $V_z$ is defined in (\ref{Vz}). From (\ref{compatible}), $S_n\vert_{\p S} = 0$ and $\p_{t}^{k}Z^{\alpha} S_n\vert_{\p S} = 0$ hold. Therefore, we get the following energy estimate.
\begin{equation} \label{Sn energy}
\frac{1}{2}\frac{d}{dt}\int_S |Z^{m-2} S_n |^2 dV_t + \varepsilon\int_S |\nabla^\varphi Z^{m-2} S_n |^2 dV_t = \int_S ( Z^{m-2} F_S + C_S ) \cdot Z^{m-2} S_n dV_t.
\end{equation}
Estimate of $C_{S_y}$ is easy. Using Proposition \ref{commutator est},
\begin{equation} \label{CSy}
\left\|C_{S_y}\right\| \leq \Lambda(\frac{1}{c_0},|h|_{Y^{[\frac{m}{2}],1}} + \|\nabla v\|_{Y^{[\frac{m}{2}],0}} ) (\left\|S_n\right\|_{X^{m-2,0}} + \left\|v\right\|_{X^{m-1,0}}  ).
\end{equation}
To control $C_{S_z}$ we give $\partial_z$ to $V_z$ by integration by part. From the commutator, we have to control the terms like,
\begin{equation}
\|Z^\beta V_z\partial_z Z^\gamma S_n \|,
\end{equation}
where $|\beta|+|\gamma| \leq m-2$, $|\gamma| \leq m-3$ or equivalently $|\beta| \neq 0$. We interchange $\partial_z$ and $Z_3$ by
\begin{equation*}
Z^\beta V_z\partial_z Z^\gamma S_n = \frac{1-z}{z} Z^\beta V_z Z_3 Z^\gamma S_n.
\end{equation*}
Then by commutation between $ \frac{1-z}{z} $ and $ Z^\beta $, we encounter the terms of the following terms,
\begin{equation*}
c_{\tilde{\beta}} Z^{\tilde{\beta}} \Big( \frac{1-z}{z} V_z \Big) Z_3 Z^\gamma S_n ,
\end{equation*}
where $c_{\tilde{\beta}}$ is sufficiently nice and bounded function, and $|\tilde{\beta}| \leq |\beta|$.  \\
If $\tilde{\beta}=0$, 
\begin{equation*}
\| c_{\tilde{\beta}} Z^{\tilde{\beta}} \Big( \frac{1-z}{z} V_z \Big) Z_3 Z^\gamma S_n \| \lesssim \Lambda(\frac{1}{c_0},|h|_{Y^{[\frac{m}{2}],1}} + \|\nabla v\|_{Y^{[\frac{m}{2}],0}} ) ( \| S_n \|_{X^{m-2,0}} + |h|_{X^{m-1,0}} ).
\end{equation*}
If $\tilde{\beta} \neq 0$,
\begin{equation*}
\Big\| c_{\tilde{\beta}} Z^{\tilde{\beta}} \Big( \frac{1-z}{z} V_z \Big) Z_3 Z^\gamma S_n  \Big\| \lesssim 
\Big\| Z\Big( \frac{1-z}{z} V_z\Big) \Big\|_{Y^{[\frac{m}{2}],0}}
\| S_n \|_{X^{m-2,0}} + 
\| S_n \|_{Y^{[\frac{m}{2}],0}}
\Big\| Z\Big( \frac{1-z}{z} V_z\Big) \Big\|_{X^{m-3,0}}.
\end{equation*}
First, we see that,
\begin{equation*}
\Big\| Z\Big( \frac{1-z}{z} V_z\Big) \Big\|_{Y^{[\frac{m}{2}],0}} \lesssim
\| V_z \|_{Y^{[\frac{m}{2}]+1,0}} + \| \partial_z V_z \|_{Y^{[\frac{m}{2}]+1,0}}
\end{equation*}
and,
\begin{equation*}
\begin{split}
\Big\| Z\Big( \frac{1-z}{z} V_z\Big) \Big\|_{X^{m-3,0}} &\lesssim
 \| \nabla V_z  \|_{X^{m-3,0}} +
\Big\| \frac{1}{z(1-z)} V_z \Big\|_{X^{m-3,0}} \\
&\lesssim 
\Big\| \frac{1-z}{z} ZV_z \Big\|_{X^{m-3,0}} +
\Big\| \frac{1}{z(1-z)} V_z \Big\|_{X^{m-3,0}} .
\end{split}
\end{equation*}
So we should estimate the terms those have forms of,
\begin{equation} \label{two types}
\Big\| \frac{1-z}{z} Z^\xi ZV_z \Big\|\quad\text{and}\quad \Big\| \frac{1}{z(1-z)} Z^\xi V_z \Big\|,
\end{equation}
where $ |\xi| \leq m-3$. To estimate these two types of terms, we use the following lemma.
\begin{lemma}
If $f(0) =$, we have the following inequalites,
\begin{equation*}
\begin{split}
\int_{-\infty}^0 \frac{1}{z^2 (1-z)^2}\left|f(z)\right|^2 dz &\lesssim \int_{-\infty}^0 \left|\partial_z f(z)\right|^2 dz,  \\
\int_{-\infty}^0 \Big(\frac{1-z}{z}\Big)^2\left|f(z)\right|^2 dz &\lesssim \int_{-\infty}^0 \Big( \left| f(z)\right|^2 + \left|\partial_z f(z)\right|^2 \Big) dz.
\end{split}
\end{equation*}
\end{lemma}
\begin{proof}
This estimate is Lemma 8.4 in \cite{NMFR}. 
\end{proof}

To estimate two types in (\ref{two types}), we use above lemma to get
\begin{equation}
\begin{split}
\Big\| \frac{1-z}{z} Z^\xi ZV_z \Big\|^2 &\lesssim \Big\| Z^\xi Z V_z \Big\|^2 + \Big\|\partial_z  (Z^\xi Z V_z )\Big\|^2,  \\
\Big\| \frac{1}{z(1-z)} Z^\xi ZV_z \Big\|^2 &\lesssim
\Big\|\partial_z  (Z^\xi Z V_z )\Big\|^2.  \\
\end{split}
\end{equation}
Therefore,
\begin{equation} \label{CSz}
\left\| C_{S_z}\right\| \lesssim
\left\| Z V_z\right\|_{X^{m-3,0}} + \left\| \partial_z Z V_z \right\|_{X^{m-3,0}} + \left\| \partial_z V_z \right\|_{X^{m-3,0}}
\end{equation}
Combining (\ref{CSy}) and (\ref{CSz}) we have
\begin{equation} \label{CS1}
\begin{split}
\left\|C_S^1\right\| &\leq \Lambda(\frac{1}{c_0},|h|_{Y^{[\frac{m}{2}],1}} + \|\nabla v\|_{Y^{[\frac{m}{2}],0}} ) \Big( \left\|S_n\right\|_{X^{m-2,0}} + \|v\|_{X^{m-2,0}} + \left\|\partial_z V_z\right\|_{X^{m-2,0}} \Big)  \\
&\leq \Lambda(\frac{1}{c_0},|h|_{Y^{[\frac{m}{2}],1}} + \|\nabla v\|_{Y^{[\frac{m}{2}],0}} ) \Big( \left\|S_n\right\|_{X^{m-2,0}} + \|v\|_{X^{m-1,0}} + \left|h\right|_{X^{m-2,1}} \Big).
\end{split}
\end{equation}
For $C_S^2$, we expand $C_S^2$ as following.
\begin{equation}
\begin{split}
\varepsilon Z^{m-2} (\Delta^\varphi S_n ) &= 
\varepsilon Z^{m-2} (\frac{1}{\partial_z\varphi}\nabla \cdot (E\nabla S_n ) ) = 
\varepsilon\frac{1}{\partial_z\varphi}Z^{m-2} (\nabla\cdot  ( E\nabla S_n ) ) + C_{S_1}^2 \\
& = \varepsilon\frac{1}{\partial_z\varphi}\nabla\cdot Z^{m-2} ( E\nabla S_n ) + C_{S_2}^2 + C_{S_1}^2   \\
& = \varepsilon\frac{1}{\partial_z\varphi}\nabla\cdot ( E\nabla Z^\alpha S_n ) + C_{S_3}^2 + C_{S_2}^2 + C_{S_1}^2   \\
& = \varepsilon\Delta^\varphi (Z^\alpha S_n) + C_S^2,
\end{split}
\end{equation}
where
\begin{equation}
\begin{split}
C_{S}^2 &:= C_{S_1}^2 + C_{S_2}^2 + C_{S_3}^2,  \\
C_{S_1}^2 &:= \varepsilon\left[Z^\beta,\frac{1}{\partial_z\varphi}\right]\nabla \cdot\Big(E\nabla S_n\Big), \\
C_{S_2}^2 &:= \varepsilon\frac{1}{\partial_z\varphi}\left[Z^\alpha,\nabla\right]\cdot\Big(E\nabla S_n\Big) \\
C_{S_3}^2 &:= \varepsilon\frac{1}{\partial_z\varphi}\nabla\cdot\Big(\left[Z^\alpha,E\nabla\right]S_n\Big).
\end{split}
\end{equation}

\noindent\textbf{1)} Estimate of $C_{S_1}^2$ \\
We need to estimate the following types of terms,
\begin{equation*}
\varepsilon\int_S Z^\beta\Big(\frac{1}{\partial_z\varphi}\Big)Z^{\tilde{\gamma}} (\nabla\cdot (E\nabla S_n ) ) \cdot Z^{m-3} S_n dV_t,
\end{equation*}
where $|\beta|+|\tilde{\gamma}|=m-2,\,\,\,\beta \neq 0 $. Then again, by commutator between $Z^{\tilde{\gamma}}$ and $\nabla$, the forms becomes like the following forms.
\begin{equation*}
\varepsilon\int_S Z^\beta\Big(\frac{1}{\partial_z\varphi}\Big)\partial_i Z^{\gamma} ( (E\nabla S_n )_j ) \cdot Z^{m-2} S_n dV_t,
\end{equation*}
where $|\tilde{\gamma}| \leq |\gamma| $. Now we preform integrate by part, to get
\begin{equation*}
\begin{split}
& \left| \varepsilon\int_S Z^\beta\Big(\frac{1}{\partial_z\varphi}\Big)\partial_i Z^{\gamma} ((E\nabla S_n)_j) \cdot Z^{m-2} S_n dV_t \right|  \\
&\leq \left| \varepsilon\int_S \partial_i Z^\beta (\frac{1}{\partial_z\varphi}) Z^{\gamma} ( (E\nabla S_n )_j ) \cdot Z^{m-2} S_n dV_t \right| + \left| \varepsilon\int_S Z^\beta \Big(\frac{1}{\partial_z\varphi}\Big) Z^{\gamma} ( (E\nabla S_n )_j ) \cdot \partial_i Z^{m-2} S_n dV_t \right|.
\end{split}
\end{equation*}
Using Proposition \ref{commutator est} and Holder inequality, we get,
\begin{equation} \label{Cs12}
\left|\int_S C_{S_1}^2 \cdot Z^{m-3} S_n dV_t \right| \lesssim
\varepsilon \Lambda(\frac{1}{c_0}, |h|_{Y^{[\frac{m}{2}],1}} + \|\nabla v\|_{Y^{[\frac{m}{2}],0}})
( \|\nabla Z^{m-2} S_n \|^2 + \|S_n\|^2_{X^{m-2,0}} + |h|^2_{X^{m-1,1}} ).
\end{equation}

\noindent\textbf{2)} Estimate of $C_{S_2}^2$ \\
We need to estimate following types.
\begin{equation*}
\varepsilon\int_S \partial_z Z^\beta \Big(E\nabla S_n\Big)\cdot Z^{m-2} S_n dV_t,
\end{equation*}
where $\beta \leq m-3 $. Then, by integration by parts again, we can get the same estimate as like $C_{S_1}^2$
\begin{equation} \label{Cs22}
\left|\int_S C_{S_2}^2 \cdot Z^{m-2} S_n dV_t \right| \leq
\varepsilon \Lambda(\frac{1}{c_0}, |h|_{Y^{[\frac{m}{2}],1}} + \|\nabla v\|_{Y^{[\frac{m}{2}],0}})
( \|\nabla Z^{m-2} S_n \|^2 + \|S_n \|^2_{X^{m-2,0}} + |h|^2_{X^{m-1,1}} ).
\end{equation}

\noindent\textbf{3)} Estimate of $C_{S_3}^2$ \\
We give $\nabla\cdot$ to $Z^{m-2} S_n$ by integration by parts, then easily,
\begin{equation} \label{Cs32}
\begin{split}
&\Big|\int_S C_{S_3}^2 \cdot Z^{m-2} S_n dV_t \Big| \leq
\varepsilon \|[Z^{m-2},E\nabla]S_n \| \|\nabla Z^{m-2} S_n \|  \\
& \leq 
\varepsilon \Lambda(\frac{1}{c_0}, |h|_{Y^{[\frac{m}{2}],1}} + \|\nabla v\|_{Y^{[\frac{m}{2}],0}})
( \|\nabla Z^{m-2} S_n \|^2 + \|S_n \|^2_{X^{m-2,0}} + |h|^2_{X^{m-1,1}} ).
\end{split}
\end{equation}
Combining above three estimates (\ref{Cs12}), (\ref{Cs22}), and (\ref{Cs32}), we have the following.
\begin{equation} \label{CS2}
\Big|\int_S C_{S}^2 \cdot Z^{m-2} S_n dV_t \Big| \leq
\varepsilon \Lambda(\frac{1}{c_0}, |h|_{Y^{[\frac{m}{2}],1}} + \|\nabla v\|_{Y^{[\frac{m}{2}],0}})
( \|\nabla Z^{m-2} S_n \|^2 +  \|S_n \|^2_{X^{m-2,0}} +  |h |^2_{X^{m-1,1}}  ).
\end{equation}
Now we use (\ref{FS1}), (\ref{FS2}), (\ref{FS3}), (\ref{CS1}), and (\ref{CS2}) to estimate the right hand side of (\ref{Sn energy}). Especially we use Young's inequality to $\varepsilon\|\nabla S_{n}\|_{X^{m-2,0}}$ in (\ref{FS2}), i.e.
\[
	\varepsilon\|Z^{m-2}S_{n}\|\|\nabla S_{n}\|_{X^{m-2,0}} \lesssim \theta \varepsilon \|\nabla S_{n}\|^{2}_{X^{m-2,0}} + \frac{1}{\theta} \|S_{n}\|^{2}_{X^{m-2,0}}.
\] 
\begin{equation}
\begin{split}
& \frac{1}{2}\frac{d}{dt}\int_S  |Z^{m-2} S_n |^2 dV_t + 
2 \varepsilon \int_S  |\nabla^\varphi Z^{m-2} S_n |^2 dV_t \\
&\leq \Lambda(\frac{1}{c_0}, |h|_{Y^{[\frac{m}{2}],1}} + \|\nabla v\|_{Y^{[\frac{m}{2}],0}}) \\
&\quad \times \Big\{ \|S_n\|^2_{X^{m-2,0}} + |h|^2_{X^{m-2,1}} + \|v\|^2_{X^{m-1,0}} + \varepsilon\Big( \theta \|\nabla S_n \|^2_{X^{m-2,0}} + |h|^2_{X^{m-2,\frac{7}{2}}} + |h|^2_{X^{m-1,1}} \Big) \Big\} ,
\end{split}
\end{equation}
for sufficiently small $\theta \ll 1$. We integrate for time from $0$ to $T$ and sum for all $Z^{m-2}$. Dissipation type term $\theta\|\nabla S_{n}\|^{2}_{X^{m-2,0}}$ is absorbed by left hand side for sufficienlty small $\theta\ll 1$, under assumption of $\sup_{t\in[0,T]}\Lambda(t) \leq C$. And then we use Proposition \ref{korn} to get the following result.
\begin{proposition}  \label{prop Sn}
Assume that $m\geq 6$ and $(v,\varphi,q,h)$ is a smooth solution on $[0,T]$ of the system (\ref{1.12})-(\ref{1.15}). Also assume that
\begin{equation}
	\partial_z\varphi \geq c_0,\quad \forall t\in [0,T] \quad \text{and} \quad \sup_{t\in[0,T]} \Lambda_{m,\infty}(t) \leq C,
\end{equation}
for some uniform constant $C$. Then we have the following $L^2$ type estimate for $S_n$.
\begin{equation*}
\begin{split}
&\left\|S_n(T)\right\|^2_{X^{m-2,0}} + \varepsilon\int_0^T \left\|\nabla S_n(t)\right\|^2_{X^{m-2,0}} dt \leq \Lambda(\frac{1}{c_0}) \sum_{\forall Z^{m-2}} \|(Z^{m-2}S_n)(0)\|^2 \\
&\quad  + \int_0^T \Lambda_{m,\infty}(t) \big( \|S_n\|^2_{X^{m-2,0}} + \|v\|^2_{X^{m-1,0}} + |h|^2_{X^{m,1}} + |h|^2_{X^{m-2,\frac{7}{2}}} \big) dt.
\end{split}
\end{equation*}
\end{proposition}

\section{$L^\infty$ type estimate}
In this section, we estimate high order $L^{\infty}$ type norm. In (\ref{Lambda infty}), we can use standard Sobolev embedding for $|h(t)|_{Y^{[\frac{m}{2}],1}}$ for sufficiently large $m$. Howver, $\|\nabla v\|_{Y^{[\frac{m}{2}],0}}$ cannot be controlled by Sobolev embedding, since conormal space is weaker than standard Sobolev space. Therefore, as we did in previous section, we use $S_n$ instead of $\p_z v$ directly. Although $S_n$ was not sufficient smooth to get $m-1$ order, we can use energy terms $\|S_n\|_{X^{m-3,0}}$ to control $L^\infty$ type terms. The following lemmas show that equivalent relation between $S_n$ and $\p_z v$ in $L^\infty$ sense. 
\begin{lemma} \label{pre normal infty control}
We have the following estimate for normal part of $\partial_z v$.
\begin{equation*}
\|\partial_z v\cdot \mathbf{n} \|_{Y^{k,0}} \leq
\Lambda(\frac{1}{c_0},\left|h\right|_{Y^{k,1}}) \|v\|_{Y^{k+1,0}} .
\end{equation*} 
\end{lemma}
\begin{proof}
From divergence free condition, 
\begin{equation*}
\partial_z v \cdot \mathbf{N} = 
\partial_z \varphi  (\partial_1 v_1 + \partial_2 v_2 ).
\end{equation*}
We take $\|\cdot\|_{Y^{k,0}}$ and use Proposition \ref{aniso embed} to get,
\begin{equation*}
\begin{split}
\|\partial_z v\cdot \mathbf{n} \|_{Y^{k,0}} &\leq \Lambda(\frac{1}{c_0},\left|h\right|_{Y^{k,1}}) \|v\|_{Y^{k+1,0}} .  \\
\end{split}
\end{equation*}
\end{proof}

\begin{proposition}
We have the following estimate.
	\begin{equation*}
		\|\p_z v\|_{Y^{k,0}} \leq \Lambda (\frac{1}{c_0},|h|_{Y^{k,1}} ) ( \|v\|_{Y^{k+1,0}} + \left\|S_n\right\|_{Y^{k,0}} ).
	\end{equation*}
\end{proposition}
\begin{proof}
We have,
$$
2S^\varphi v\mathbf{n} = (\nabla u)\mathbf{n} + (\nabla u)^T\mathbf{n} = (\nabla u)\mathbf{n} + g^{ij}(\partial_j v\cdot\mathbf{n})\partial_{y^i}
$$
and divergence free condition gives, 
$$
\partial_N u = \frac{1+\left|\nabla_y\varphi\right|^2}{\partial_z \varphi}\partial_z v - \partial_1\varphi\partial_1 v - \partial_2\varphi\partial_2 v.
$$
So we obtain,
\begin{equation} \label{ttt}
\|\partial_z v\|_{Y^{k,0}} \leq \Lambda(\frac{1}{c_0},\left|h\right|_{Y^{k,1}}) ( \|v\|_{Y^{k+1,0}} + \|S^\varphi v\mathbf{n}\|_{Y^{k,0}}  ) .
\end{equation}
And, since
$$
S^\varphi v\mathbf{n} = S_n + (\mathbf{n} \otimes \mathbf{n}) (S^\varphi v\mathbf{n} ),
$$
we get
\begin{equation*}
 \| S^\varphi v\mathbf{n}  \|_{Y^{k,0}} \leq  \|S_n\|_{Y^{k,0}} + \Lambda (\frac{1}{c_0},|h|_{Y^{k,1}} ) \|\partial_z v\cdot \mathbf{n}\|_{Y^{k,0}} + \left\|v\right\|_{Y^{k,0}}.
\end{equation*}
Using Lemma \ref{pre normal infty control} and (\ref{ttt}),
\begin{equation*}
\begin{split}
 \|\p_z v\|_{Y^{k,0}} &\leq \Lambda (\frac{1}{c_0},|h|_{Y^{k,1}} ) ( \|v\|_{Y^{k+1,0}} + \left\|S_n\right\|_{Y^{k,0}} ).
\end{split}
\end{equation*}
\end{proof}

Above proposition implies that we suffice to estimate $\|S_n\|_{Y^{[\frac{m}{2}],0}}$, instead of $\|\partial_z v\|_{Y^{[\frac{m}{2}],0}}$. So we use equation for $S_n$ with Dirichlet boundary condition. Main difficulty in this section is commutator between $Z_3 = \frac{z}{1-z}\partial_z$ and $\Delta^\varphi$. This commutator was not a problem in basic $L^2$-type energy estimate of $v$ and $S_n$, because the highest order commutator, which looks like $\sim \varepsilon Z^\alpha\partial_z S_n$ can be absorbed into dissipation term in the energy. But, in $L^\infty$ estimate, we use the following maximal principle for convection-diffusion equation. For $S_n$ equation,
\begin{equation*}
\p_t^\varphi S_n + (v\cdot\nabla^\varphi)S_n - \varepsilon\Delta^\varphi S_n = F_S.
\end{equation*}
We have the following $L^\infty$-type estimate which does not have dissipation in energy term.
\begin{equation*}
 \|S_n(t) \|_{L^\infty} \leq \|S_n(0) \|_{L^\infty} + \int_0^t \|F_S(s) \|_{L^\infty} \dd s.
\end{equation*}  

\noindent Note that, standard Sobolev embedding does not hold for co-normal space in general. This is because of behavior of near the boundary. Away from the boundary, however, $\frac{z}{1-z}$ is not zero, and its all order derivative for $z$ is always uniformly bounded. Hence, we divide co-normal function into two parts, one is supported near the boundary and another is supported away from the boundary. Then the second part is easy to control by standard Sobolev embedding. For the first part, we deform the coordinate so that locally $\partial_{zz}^\varphi$ look like $\partial_{zz}$. Then $\partial_z$ commute with $\partial_{zz}$, so it does not generate any harmful ( which has 1-more order than $L^\infty$-type energy) commutator. This clever idea is introduced in \cite{NMFR2}. We introduce this system briefly and use similar arguments to get the result. First, we start with very simple lemma, which means, co-normal Sobolev space is equivalent to standard Sobolev away from the boundary .

\begin{lemma} \label{away from bdry}
For any smooth cut-off function $\bar{\chi}$ such that $\bar{\chi} = 0$ in a vicinity of $z=0$, we have for $m>k+3/2$:
\begin{equation*}
\left\|\bar{\chi} f\right\|_{W^{k,\infty}} \lesssim \left\|f\right\|_{H_{co}^m}.
\end{equation*}
\end{lemma}

\noindent We decompose $Z^k S_n$ as 
\begin{equation*}
\|Z^k S_n\|_{L^\infty} \leq \|\chi Z^k S_n \|_{L^\infty} + \|v\|_{Y^{k+1,0}}.
\end{equation*}

\noindent To estimate second term $\left\|v\right\|_{Y^{k+1,0}}$, we use the following proposition.
\begin{proposition}
We have the following estimate.
\begin{equation}
\left\|v\right\|_{Y^{k,0}} \lesssim \Lambda \Big(\frac{1}{c_0},\|\nabla v \|_{Y^{1,0}} + \|v\|_{X^{k+2,0}} + |h|_{X^{k+1,1}} + \|S_n\|_{X^{k+1,0}} \Big).
\end{equation}
\end{proposition}
\begin{proof}
Using anisotropic Sobolev embedding in Proposition \ref{aniso embed},
$$
\|Z^k v\|_{L^\infty}^2 \lesssim \|\partial_z v\|_{X^{k+1,0}} \|v\|_{X^{k+2,0}} \leq \|\partial_z v\|^2_{X^{k+1,0}} + \|v\|^2_{X^{k+2,0}}
$$
and using lemme of previous section,
$$
\left\|\partial_z v\right\|_{X^{k+1,0}} \lesssim \Lambda\Big(\frac{1}{c_0},\left\|\nabla v\right\|_{Y^{\frac{k+2}{2},0}}\Big)\Big( \|v\|_{X^{k+2,0}} + |h|_{X^{k+1,1}} + \|S_n\|_{X^{k+1,0}} \Big).
$$
Then using induction for $\left\|\nabla v\right\|_{Y^{\frac{k+2}{2},0}}$, until it become 1. And notice that $\left\|v\right\|_{X^{k+2,0}}$ is absorbed by estimate of $\left\|\partial_z v\right\|_{X^{k+1,0}}$.
\end{proof}

Note that above proposition means that we suffice to estimate $Z^k S_n$ only near the boundary, so now we introduce modified coordinate which was introduced in [15] and [1]. Let us define transformation $\Psi$,
\begin{equation}
\Psi(t,\cdot):S = \mathbb{R}^2 \times (-\infty,0) \rightarrow \Omega_t,
\end{equation}
$$
x=(y,z) \mapsto \begin{pmatrix} y\\ h(t,y) \end{pmatrix} + z\mathbf{n}^b(t,y),
$$
where $\mathbf{n}^b$ is unit normal at the boundary, $(-\nabla h,1)/|N|$. To show that this is diffeomorphism near the boundary, we check
$$
D\Psi(t,\cdot) = \begin{pmatrix} 1 & 0 & -\partial_1 h \\ 
                                 0 & 1 & -\partial_2 h \\
                                 \partial_1 h & \partial_2 h & 1 \end{pmatrix} + \begin{pmatrix} -z\partial_{11} h & -z\partial_{12} h & 0\\ 
                                -z\partial_{21} h & -z\partial_{22} h & 0\\
                                 0 & 0 & 1 \end{pmatrix} .
$$
This is diffeomorphism near the boundary since norm of second matrix is controlled by $|h|_{2,\infty}$. So, we restrict $\Psi(t,\cdot)$ on $\mathbb{R}^2 \times (-\delta,0)$ so that it is diffeomorphism.($\delta$ is depend on $c_0$. Of course, think that above support separation was done by $\chi(z) = \kappa(\frac{z}{\delta(c_0)})$. Now we write laplacian $\Delta^\varphi$ with respect to Riemannian metric of above parametrization. Riemannian metric becomes,
\begin{equation}
g(y,z) = \begin{pmatrix}
\tilde{g}(y,z) & 0 \\
0 & 1 
\end{pmatrix},
\end{equation}
where $\tilde{g}$ is $2\times 2 $ block matrix. And with this metric, laplacian becomes,
\begin{equation}
\Delta_g f = \partial_{zz} f + \frac{1}{2}\partial_z(\text{ln}|g|)\partial_z f + \Delta_{\tilde{g}} f,
\end{equation}
where
\begin{equation}
\Delta_{\tilde{g}} f = \frac{1}{|\tilde{g}|^{\frac{1}{2}}}\sum_{1\leq i,j \leq 2} \partial_{y^i}(\tilde{g}^{ij}|\tilde{g}|^{\frac{1}{2}}\partial_{y^j} f),
\end{equation}
where $\tilde{g}^{ij}$ is inverse matrix element of $\tilde{g}$. We now solve problem in domain of $\Psi$. We restrict $S^\varphi v$ near the boundary and parametrize them via $\Psi$. Let us define
\begin{equation}
S^\chi := \chi(z)S^\varphi v,
\end{equation}
where $\chi(z) = \kappa(\frac{z}{\delta(c_0)})\in [0,1]$ where $\kappa$ is smooth and compactly supported near the boundary, taking value 1 there. Equation for $S^\chi$ is
\begin{equation}
\partial_t^\varphi S^\chi + (v\cdot\nabla^\varphi)S^\chi - \varepsilon\Delta^\varphi S^\chi = F_{S^\chi},
\end{equation}
where 
$$
F_{S^\chi} = F^\chi + F_v,
$$
$$
F^\chi = (V_z \partial_z\chi )S^\varphi v - \varepsilon\nabla^\varphi\chi\cdot\nabla^\varphi S^\varphi v - \varepsilon\Delta^\varphi \chi S^\varphi v,
$$
$$
F_v = -\chi(D^\varphi)^2 q - \frac{\chi}{2}((\nabla^\varphi v)^2 + ((\nabla^\varphi v)^t)^2).
$$
Note that $F^\chi$ is supported away from the boundary. We rewrite this function on our new frame by taking $\Phi^{-1}\circ\Psi$. We define
\begin{equation}
S^\Psi(t,y,z) = S^\chi(t,\Phi^{-1}\circ \Psi)
\end{equation}
and $S^\Psi$ solves
\begin{equation}
\partial_t S^\Psi + w\cdot\nabla S^\Psi - \varepsilon(\partial_{zz} S^\Psi + \frac{1}{2}\partial_z(\text{ln}|g|)\partial_z S^\Psi + \Delta_{\tilde{g}} S^\Psi) = F_{S^\chi}(t,\Phi^{-1}\circ\Psi),
\end{equation}
where 
$$
w = \bar{\chi}(D\Psi)^{-1}(v(t,\Phi^{-1}\circ\Psi)-\partial_t\Psi),
$$
where $\bar{\chi}$ is slightly larger support so that $\bar{\chi}S^\Psi = S^\Psi$ and as like $S^\chi$, $S^\Psi$ is also only supported near the boundary. In this frame $S_n$ correspond to $S_n^\Psi$, which is defined as following.
\begin{equation}
S_n^\Psi (t,y,z) = \Pi^b(t,y)S^\Psi\mathbf{n}^b(t,y) = \Pi^b(t,y)S^\chi(t,\Phi^{-1}\circ \Psi)\mathbf{n}^b(t,y),
\end{equation}
where $\Pi^b = I - \mathbf{n}^b\otimes\mathbf{n}^b$.(tangential operator at the boundary, so they are independent to $z$.) Then equation for $S_n^\Psi$ becomes,
\begin{equation} \label{SnPsi}
\partial_t S_n^\Psi + w\cdot \nabla S_n^\Psi - \varepsilon(\partial_{zz} + \frac{1}{2}\partial_z(\text{ln}|g|)\partial_z)S_n^\Psi = F_n^\Psi,
\end{equation}
where
$$
F_n^\Psi = \Pi^b F_{S_\chi}\mathbf{n}^b + F_n^{\Psi,1} + F_n^{\Psi,2},
$$
where
$$
F_n^{\Psi,1} = ((\partial_t + w_y\cdot\nabla_y)\Pi^b)S^\Psi\mathbf{n}^b + \Pi^b S^\Psi(\partial_t + w_y\cdot\nabla_y)\mathbf{n}^b,
$$
$$
F_n^{\Psi,2} = -\varepsilon\Pi^b(\Delta_{\tilde{g}}S^\Psi)\mathbf{n}^b,
$$
with zero-boundary condition at $z=0$. Note that $S_n = S_n^\Psi$ on the 	boundary. We will estimate $S_n^\Psi$ instead of $S_n$. to validate this, we should show that equivalence of these two terms. Firstly, by definition of $S_n^\Psi$,
\begin{equation}
\|S_n^\Psi\|_{Y^{k,0}} \leq \Lambda( |h|_{Y^{k+1,0}} ) \|\Pi^b S^\varphi v \mathbf{n}^b \|_{Y^{k,0}}
\end{equation}
and since  $|\Pi-\Pi^b|+|\mathbf{n}-\mathbf{n}^b| = O(z)$ near the boundary $z=0$,
\begin{equation}
\|S_n^\Psi\|_{Y^{k,0}} \leq \Lambda(\frac{1}{c_0},\|S_n\|_{Y^{k,0}} + \|v\|_{Y^{k+1,0}} ).
\end{equation}
Now, we apply anisotropic Sobolev embedding to the last term,
$$
\|S_n^\Psi\|_{Y^{k,0}} \leq \Lambda(\frac{1}{c_0},\|S_n\|_{Y^{k,0}} + \|\partial_z v\|_{X^{k+2,0}} + \|v\|_{X^{k+3,0}} ).
$$
For $\|\partial_z v\|_{X^{k+2,0}}$, we use Lemma 8.2 inductively,(to reduce the order of $\|\nabla v\|_{Y^{\frac{k}{2},0}}$, to get
\begin{equation}
\|S_n^\Psi\|_{Y^{k,0}} \leq \Lambda(\frac{1}{c_0},\|v\|_{X^{k+3,0}} + \|S_n\|_{X^{k+2,0}} + |h|_{X^{k+2,0}} + \|S_n\|_{Y^{k,0}} ).
\end{equation}
Since we choose sufficiently smaller $k$ than $m$, this estimate is okay. For opposite direction, we can do similarly to get
\begin{equation}
\|S_n\|_{Y^{k,0}} \leq \Lambda(\frac{1}{c_0},\|v\|_{X^{k+3,0}} + \|S_n\|_{X^{k+2,0}} + |h|_{X^{k+2,0}} + \|S_n^\Psi\|_{Y^{k,0}} ).
\end{equation}  
So, we finish equivalence argument.\\
Now we should apply $Z^k$ to the system (9.17). Applying tangential derivative ($Z_1,Z_2$) is not that harmful, but commutator between $Z_3$ and Laplacian is still a problem. Critical observation in \cite{NMFR} is the following Lemma. (Lemma 9.6 in \cite{NMFR}).
\begin{lemma}
Consider $\rho$ a smooth solution of
\begin{equation}
\partial_t\rho + w\cdot\nabla\rho = \varepsilon\partial_{zz}\rho + \mathcal{H},\,\,-1<z<0,\,\,\rho(t,y,0) = 0,\,\,\rho(t=0)=\rho_0
\end{equation}
for some smooth vector field $w$ such that $w_3$ vanishes on the boundary. Assume that $\rho$ and $\mathcal{H}$ are compactly supported in $z$. Then we have the estimate:
$$
\|Z_i\rho(t)\|_{\infty} \lesssim \|Z_i\rho_0\|_{\infty} + \|\rho_0\|_{\infty} + \int_0^t\Big( (\|w\|_{E^{2,\infty}} + \|\partial_{zz}w_3\|_{L^\infty})(\|\rho\|_{1,\infty} + \|\rho\|_4) + \|\mathcal{H}\|_{1,\infty} \Big),\,\,i=1,2,3.
$$
\end{lemma}
We should generalize this to high order, since we need $k$-order $L^\infty$ -type estimate. Let us first introduce rewriting of system (\ref{SnPsi}) to circumvent difficulty. We set
\begin{equation}
\rho(t,y,z) := |g|^{\frac{1}{4}}S_n^\Psi
\end{equation}
and then $\rho$ solves,
\begin{equation} \label{rho eq}
\partial_t\rho + w\cdot\nabla\rho -\varepsilon\partial_{zz}\rho = |g|^{\frac{1}{4}}(F_n^\Psi + F_g) \doteq \mathcal{H},
\end{equation}
where
$$
F_g = \frac{\rho}{|g|^{\frac{1}{2}}}\Big(\partial_t + w\cdot\nabla - \varepsilon\partial_{zz}\Big)|g|^{\frac{1}{4}},
$$
which shows that $\varepsilon\partial_z \text{ln}|g|\partial_z$ is removed. And trivially, $Z_3 S_n^\Psi$ and $\rho$ are equivalent, i.e
\begin{equation}
\|\rho\|_{Y^{k,0}} \leq \Lambda(|h|_{Y^{k+1,0}}) \|S_n^\Psi\|_{Y^{k,0}},\,\,\,\|S_n^\Psi\|_{Y^{k,0}} \leq \Lambda(|h|_{Y^{k+1,0}}) \|\rho\|_{Y^{k,0}}.
\end{equation}
Hence, instead of $S_n^\Psi$, we estimate $\rho$. Also note that equation of $\rho$ is applicable above lemma. Now we extend above lemma to high order. Since we split horizontal components and vertical component, we use $y$ to denote two horizontal coodinates. 
\begin{lemma} \label{high order version}
(High order version)Consider $\rho$ a smooth solution of
\begin{equation}
\partial_t\rho + w\cdot\nabla\rho = \varepsilon\partial_{zz}\rho + \mathcal{H},\,\,-1<z<0,\,\,\rho(t,y,0) = 0,\,\,\rho(t=0)=\rho_0
\end{equation}
for some smooth vector field $w$ such that vertical component $w_3$ vanishes on the boundary. Assume that $\rho$ and $\mathcal{H}$ are compactly supported in $z$. Then for any integer $k$, we have the following estimate:
$$
\|Z^k\rho\|_{L^\infty} \lesssim \|\rho_0\|_{Y^{k,0}} + \int_0^t \Big( \Big( \|\partial_z w_{y}\|_{Y^{k,0}} + \|\partial_{zz} w_3\|_{Y^{k,0}} \Big) \|\rho\|_{X^{k+3,0}} + \|\mathcal{H}\|_{Y^{k,0}} \Big) d\tau ,
$$
where $w_{y}$ is horizontal component $(w_{1}, w_{2})$.
\end{lemma}
\begin{proof}
Applying co-normal derivatives to equation generate bad commutator which come from between $Z_3$ and Laplacian. For simplicity, we use notations $\vec{y} = (x,y)$ and horizontal derivative $\nabla_{\vec{y}}:=(\p_{x},\p_{y})$ in this proof. We rewrite equation as
\begin{equation}
\partial_t\rho + z\partial_z w_3(t,\vec{y},0)\partial_z\rho + w_{\vec{y}}(t,\vec{y},0)\cdot\nabla_{\vec{y}}\rho - \varepsilon\partial_{zz}\rho = \mathcal{H} - \mathcal{R} \doteq G,
\end{equation}
where ($w_3 = 0$ on the boundary) 
$$
\mathcal{R} = (w_y(t,y,z) - w_y(t,y,0))\cdot\nabla_y\rho + (w_3(t,y,z) - z\partial_z w_3(t,y,0))\partial_z\rho.
$$
We use evolution operator $S(t,\tau)$ for homogeneous solution of above system. Let,
\begin{equation}
\rho(t,y,z) = S(t,\tau)f_0(y,z),\,\,f(\tau,y,z) = f_0(y,z)\,(initial\,condition)
\end{equation}
solves
\begin{equation}
\partial_t\rho + z\partial_z w_3(t,y,0)\partial_z\rho + w_y(t,y,0)\cdot\nabla_y\rho - \varepsilon\partial_{zz}\rho = 0,\,\,-1<z<0,\,\,t>\tau,\,\,\rho(t,y,0) = 0.
\end{equation}
For the full non-homogeneous system, by Duhamel's formula,
\begin{equation}
\rho(t) = S(t,0)\rho_0 + \int_0^t S(t,\tau)G(\tau)d\tau.
\end{equation}
Now, suppose that $\rho$ is compactly supported in $z$,(near the boundary) and $z<0$, then 
$$
\|Z_3^k\rho\|_{L^\infty} \lesssim \sum_{i,j=1}^k |z^j\partial_z^i\rho|_{L^\infty}.
$$
Since $z$ is near boundary, we don't have to consider relatively small index $j$, which means
\begin{equation}
\|Z_3^k\rho\|_{L^\infty} \lesssim \sum_{j=1}^k \sum_{i=1}^j |z^j\partial_z^i\rho|_{L^\infty}.
\end{equation}
To estimate each terms on right hand side, we should control each $|z^j\partial_z^i\rho|_{L^\infty}$.

\begin{lemma}
For evolution operator $S$ as above, we have following estimate.
\begin{equation}
|z^j\partial_z^i S(t,\tau)\rho_0|_{L^\infty} \lesssim \|\rho_0\|_{L^\infty} + \sum_{i_1+i_2 = i}\|z^{j-i_1}\partial_z^{i_2} \rho_0\|_{L^\infty}.
\end{equation}
\end{lemma}

\begin{proof}
Basically we follow the method of Lemma 9.6 in [1]. Let $\rho(t,y,z) = S(t,\tau)\rho_0(y,z)$ solves homogeneous system of (9.31). We extend this variables to whole space by 
\begin{equation}
\tilde{\rho}(t,y,z) = \rho(t,y,z),\,\,z>0,\,\,\,\,\tilde{\rho}(t,y,z) = -\rho(t,y,-z),\,\,z<0.
\end{equation}
So that $\tilde{\rho}$ solves
\begin{equation}
\partial_t\tilde{\rho} + z\partial_z w_3(t,y,0)\partial_z\tilde{\rho} + w_y(t,y,0)\cdot\nabla_y\tilde{\rho} - \varepsilon\partial_{zz}\tilde{\rho} = 0,\,\, z\in\mathbb{R},
\end{equation}
with initial condition $\tilde{\rho}(\tau,y,z) = \tilde{\rho}_0(y,z)$.\\
By introducing $\mathcal{E}$, which solves,
$$
\partial_t\mathcal{E} = w_y(t,\mathcal{E},0),\,\,\mathcal{E}(\tau,\tau,y) = y
$$
and define
$$
g(t,y,z) = \rho(t,\mathcal{E}(t,y,z),z).
$$
Then $g$ solves,
\begin{equation}
\partial_t g + z\gamma(t,y)\partial_z g - \varepsilon\partial_{zz}g = 0,\,\,z\in\mathbb{R},\,\,g(\tau,y,z) = \tilde{\rho}_0(y,z),
\end{equation}
where
$$
\gamma(t,y) = \partial_z w_3(t,\mathcal{E}(t,y,z),0).
$$
By using Fourier transform, we get explicit form of the solution,
\begin{equation}
g(t,y,z) = \int_{\mathbb{R}} k(t,\tau,y,z-z')\tilde{\rho}_0(y,e^{-\Gamma(t)}z')dz',
\end{equation} 
where
$$
k(t,\tau,y,z-z') = \frac{1}{\sqrt{4\pi\varepsilon\int_\tau^t e^{2\varepsilon(\Gamma(t)-\Gamma(s))}ds}}exp\Big(-\frac{(z-z')^2}{4\varepsilon\int_\tau^t e^{2\varepsilon(\Gamma(t)-\Gamma(s))}ds}\Big),\,\,\int_{\mathbb{R}}k(t,\tau,y,z) dz = 1
$$
$$
\Gamma(t) = \int_\tau^t\gamma(s,y)ds.
$$
We note that,
\begin{equation}
z^j\partial_z^i g = \int_{\mathbb{R}} \Big( z^j\partial_z^i k(t,\tau,z-z')\Big)\tilde{\rho}_0(y,e^{-\Gamma(t)}z')dz' 
\end{equation}
$$
= \int_{\mathbb{R}} \Big( (z^j-z'^j)\partial_z^i k(t,\tau,y,z-z') + (-1)^i z'^j\partial_{z'}^i k(t,\tau,y,z-z') \Big)\tilde{\rho}_0(y,e^{-\Gamma(t)}z')dz' .
$$
Now, since $k$ is Gaussian,
\begin{equation}
\int_{\mathbb{R}} |(z^j-z'^j)\partial_z^i k| dz' \lesssim 1.
\end{equation}
So, using integration by parts on the 2nd term, we can deduce
\begin{equation}
\begin{split}
\|z^j\partial_z^i g\|_{L^\infty} &\lesssim \|\tilde{\rho}_0\|_{L^\infty} + \Big\| \int_{\mathbb{R}}\partial_{z'}^{i-1}k(t,\tau,y,z-z')\{jz'^{j-1}\tilde{\rho}_0(y,e^{-\Gamma(t)}z')\\\
&\quad  + z'^j\partial_{z'}\tilde{\rho}_0(y,e^{-\Gamma(t)}z')e^{-\Gamma(t)}\} dz' \Big\|_{L^\infty} \\
&\lesssim \cdots \\
&\lesssim \|\tilde{\rho}_0\|_{L^\infty} + \Big\|\int_{\mathbb{R}} k(t,\tau,y,z-z') \\
&\quad\quad\quad\quad\quad\quad\quad \times \Big\{ \sum_{i_1+i_2 = i} (z')^{j-i_1} \partial_{z'}^{i_2} \tilde{\rho}_0(y,e^{-\Gamma(t)}z') e^{-{i_2}\Gamma(t)}\} \Big\} dz' \Big\|_{L^\infty}  \\
&\lesssim \|\tilde{\rho}_0\|_{L^\infty} + \sum_{i_1+i_2 = i} \left\|\int_{\mathbb{R}} k(t,\tau,y,z-z')\left\{ (z')^{j-i_1} \partial_{z'}^{i_2} \tilde{\rho}_0(y,e^{-\Gamma(t)}z') \right\} dz'\right\|_{L^\infty}.
\end{split}
\end{equation}
By relation of $\rho$ and $g$, we get
\begin{equation}
\|z^j\partial_z^i \rho\|_{L^\infty} \leq \|z^j\partial_z^i \tilde{\rho}\|_{L^\infty} \lesssim \|\tilde{\rho}_0\|_{L^\infty} + \sum_{i_1+i_2 = i}\|z^{j-i_1}\partial_z^{i_2}\tilde{\rho}_0\|_{L^\infty}
\end{equation}
$$
\leq \|z^j\partial_z^i \tilde{\rho}\|_{L^\infty} \lesssim \|\rho_0\|_{L^\infty} + \sum_{i_1+i_2 = i}\|z^{j-i_1}\partial_z^{i_2} \rho_0\|_{L^\infty}.
$$
\end{proof}
Now we apply $Z_3^k$ to Duhamel's formula to get
\begin{equation}
Z_3^k \rho(t) = Z_3^k(S(t,\tau)\rho_0) + \int_0^t Z_3^k(S(t,\tau)G(\tau))d\tau.
\end{equation}
Using above Lemma 9.7 twice on the right hand side, 
\begin{equation}
\|Z_3^k \rho(t)\|_{L^\infty} \lesssim \sum_{i=1}^k \sum_{j=1}^i \left\{ \|\rho_0\|_{L^\infty} + \sum_{j+1+j_2 = j}\|z^{j-j_1}\partial_z^{j_2} \rho_0\|_{L^\infty} \right\} 
\end{equation}
$$
+ \sum_{i=1}^k \sum_{j=1}^i \int_0^t \left\{ \|G\|_{L^\infty} + \sum_{j+1+j_2 = j}\|z^{j-j_1}\partial_z^{j_2} G\|_{L^\infty} \right\} d\tau.
$$
Using the fact that $\rho$ (also $\rho_0$) and $G$ are compactly supported in $z$, we have
$$
\|Z_3^k \rho(t)\|_{L^\infty} \lesssim \|\rho_0\|_{W^{k,\infty}} + \int_0^t \|G\|_{W^{k,\infty}} d\tau.
$$
We also note that other tangential derivatives cases also holds.(This is easier than $Z_3$ case.) Hence 
\begin{equation}
\|Z^k \rho(t)\|_{L^\infty} \lesssim \|\rho_0\|_{Y^{k,0}} + \int_0^t \|G\|_{Y^{k,0}} d\tau.
\end{equation}
Let us estimate $\|\mathcal{R}\|_{Y^{k,0}}$. Since $\rho$\ is compactly supported in $z$ (near the boundary), using Taylor's series and inserting function $\zeta(z) := z(z+1)$ (inserting this function is very useful, because existence of $\zeta(z)$ enables us to control using co-normal derivatives of $\rho$), we have
\begin{equation}
\begin{split}
\|\mathcal{R}\|_{Y^{k,0}} &\leq \|(w_y(t,y,z) - w_y(t,y,0))\cdot\nabla_y\rho\|_{Y^{k,0}} \\
&\quad + \|(w_3(t,y,z) - z\partial_z w_3(t,y,0))\partial_z\rho\|_{Y^{k,0}} \\
&\leq \|\partial_z w_y\|_{Y^{k,0}} \|\zeta(z)\rho\|_{Y^{k+1,0}} + \|\partial_{zz} w_3\|_{Y^{k,0}} \|\zeta^2(z)\partial_z\rho\|_{Y^{k,0}}  \\
&\leq \Big( \|\partial_z w_y\|_{Y^{k,0}} + \|\partial_{zz} w_3\|_{Y^{k,0}} \Big) \Big( \|\zeta(z)\rho\|_{Y^{k+1,0}} + \|\zeta^2(z)\partial_z\rho\|_{Y^{k,0}} \Big).
\end{split}
\end{equation}
Using anisotropic Sobolev embedding for co-normal derivatives,
\begin{equation}
\begin{split}
\|\zeta(z)\rho\|_{Y^{k+1,0}} &\lesssim \|\zeta(z)\rho\|_{X^{k+1,2}} + \|\partial_z(\zeta(z)\rho)\|_{X^{k+1,1}}  \\
&\leq \|\zeta(z)Z^{k+3}\rho\|_{L^2} + \|\partial_z(\zeta(z)Z^{k+2}\rho)\|_{L^2} \\
&\lesssim \|\rho\|_{X^{k+3,0}} + \|\zeta'(z)\rho\|_{X^{k+2,0}} + \|\rho\|_{X^{k+3,0}} \\
& \lesssim \|\rho\|_{X^{k+3,0}},  \\
\|\zeta^2(z)\partial_z\rho\|_{Y^{k,0}} &= \|\zeta(z) Z_3\rho\|_{Y^{k,0}}.
\end{split}
\end{equation}
and $\zeta(z)$ is nice bounded function for all order of derivatives, so at result,
\begin{equation}
\|\mathcal{R}\|_{Y^{k,0}} \lesssim \Big( \|\partial_z w_y\|_{Y^{k,0}} + \|\partial_{zz} w_3\|_{Y^{k,0}} \Big) \|\rho\|_{X^{k+3,0}}.
\end{equation}
Combining with (9.45), we finish the proof.
\end{proof}
Now, we are ready to get energy estimate for $\|S_n\|_{Y^{k,0}}$.
\begin{proposition} \label{L infty estimate}
Let us define non-dissipation type energy ${E}_{m-1}$ as
\begin{equation} \label{cal A def}
{E}_{m-1} := \Lambda\Big( \frac{1}{c_0}, \|v\|_{X^{m-1,0}} + |h|_{X^{m-1,1}} + \|\partial_z v\|_{X^{m-3,0}} + \|\partial_z v\|_{Y^{[\frac{m-1}{2}],0}}   \Big).
\end{equation}
We have the following estimate for $\|S_n\|_{Y^{[\frac{m}{2}],0}}$ for sufficiently large $m\geq 6$.
\begin{equation}
\|S_n(t)\|_{Y^{[\frac{m}{2}],0}}^2 \leq \|S_n(0)\|_{Y^{[\frac{m}{2}],0}}^2 + \int_{0}^{T} {E}^{2}_{m-1}(t) dt + \varepsilon\int_{0}^{T} \|\nabla S_{n}\|^{2}_{X^{m-2,0}}dt , 
\end{equation}
Note that ${E}_{m-2}$ is equivalent to alternate energy
$$
\mathcal{Q}_{m-1} :=  \|v\|_{X^{m-1,0}} + |h|_{X^{m-1,1}} + \|S_n\|_{X^{m-3,0}} + \|S_n\|_{Y^{[\frac{m-1}{2}],0}}, 
$$ 
when $\sup_{t\in[0,T]} \Lambda_{m,\infty}(t) \leq C$, for some uniform constant $C$.
\end{proposition}
\begin{proof}
We already transformed $S_n$ equation into equivalent-$\rho$ equation system (9.26). From the result of Lemma 9.6, we should estimate the following four terms. Here, $k = [\frac{m}{2}]$, and $m$ is sufficiently large.
$$
\|\rho\|_{X^{k+3,0}},\,\,\,\|\partial_z w_y\|_{Y^{k,0}},\,\,\,\|\partial_{zz} w_3\|_{Y^{k,0}},\,\,\,\|\mathcal{H}\|_{Y^{k,0}}.
$$
1), 2) $\|\rho\|_{X^{k+3,0}}$ and $\|\partial_z w_y\|_{Y^{k,0}}$ are trivially controlled by ${E}_{m-2}$ by definition of $\rho$.\\
3) $\|\partial_{zz} w_3\|_{Y^{k,0}}$ : By definition of $w$,
$$
w = \bar{\chi}(D\Psi)^{-1}(v(t,\Phi^{-1}\circ\Psi)-\partial_t\Psi).
$$
\begin{equation}
\|\partial_{zz} w_3\|_{Y^{k,0}} \leq \|\partial_{zz}\Big(\bar{\chi}(D\Phi)^{-1}\partial_t \Psi\Big)\|_{Y^{k,0}} + \|\partial_{zz}\Big( \bar{\chi}(D\Phi)^{-1} v(t,\Phi^{-1}\circ\Psi) \Big)_3 \|_{Y^{k,0}}.
\end{equation}
For the first term,
$$
\|\partial_{zz}\Big(\bar{\chi}(D\Phi)^{-1}\partial_t \Psi\Big)\|_{Y^{k,0}} \lesssim \Lambda\Big(\frac{1}{c_0}, |h|_{Y^{k,0}} + |\partial_t h|_{Y^{k,0}} \Big) \lesssim {E}_{m-2}.
$$
For the second term, by using commutators,
$$
\|\partial_{zz}\Big( \bar{\chi}(D\Psi)^{-1} v_3(t,\Phi^{-1}\circ\Psi) \Big)\|_{Y^{k,0}} \leq \|\bar{\chi}\partial_{zz}\Big((D\Psi(t,y,0))^{-1} v(t,\Phi^{-1}\circ\Psi)\Big)_3 \|_{Y^{k,0}} + {E}_{m-2} \|\partial_z v\|_{Y^{k,0}}
$$
$$
= \| \bar{\chi}\partial_{zz}\Big( v(t,\Phi^{-1}\circ\Psi) \cdot\mathbf{n}^b \Big) \|_{Y^{k,0}} + {E}_{m-2},
$$
where we used
$$
\Big((D\Psi(t,y,0))^{-1} f \Big)_3 = f \cdot \mathbf{n}^b.
$$
Main difficulty is two normal derivatives. Meanwhile, by definition, $v(t,\Phi^{-1}\circ\Psi) = u(t,\Psi) \doteq u^\Psi$ and $u$ is divergence free is zero. In the local coordinate, divergence free condition gives,
\begin{equation}
\partial_z u^\Psi \cdot \mathbf{n}^b = -\frac{1}{2}\partial_z(\text{ln}|g|)u^\Psi\cdot\mathbf{n}^b - \nabla_{\tilde{g}}\cdot u_y^\Psi,
\end{equation}
which means one normal derivative is replaced by tangential derivative. So,
$$
\| \bar{\chi}\partial_{zz}\Big( v(t,\Phi^{-1}\circ\Psi) \cdot\mathbf{n}^b \Big) \|_{Y^{k,0}} \lesssim \Lambda\Big(\frac{1}{c_0}, \|\partial_z u^\Psi\|_{Y^{k,0}} + |h|_{Y^{k+3,0}} \Big) \lesssim {E}_{m-2}.
$$
In conclusion, $\|\partial_{zz} w_3\|_{Y^{k,0}} \lesssim {E}_{m-2}$.\\
4) $\|\mathcal{H}\|_{Y^{k,0}}$ : We have
\begin{equation}
\|F_n^\Psi\|_{Y^{k,0}} \leq \|F_n^{\Psi,1}\|_{Y^{k,0}} + \|F_n^{\Psi,2}\|_{Y^{k,0}} + \|\Pi^b F_{S^\chi} \mathbf{n}^b \|_{Y^{k,0}}.
\end{equation}
For the first term,
$$
\|F_n^{\Psi,1}\|_{Y^{k,0}} \lesssim \Lambda\Big(\frac{1}{c_0},|\partial_t h|_{Y^{k+1,0}} + |h|_{Y^{k+2,0}} + \|w\|_{Y^{k,0}} \Big) \|\partial_z v\|_{Y^{k,0}} \lesssim {E}_{m-2},
$$
where ${E}_{m-2}$ is defined in (\ref{cal A def}). For second term,
$$
\|F_n^{\Psi,2}\|_{Y^{k,0}} \lesssim \varepsilon{E}_{m-2} \|S_n\|_{Y^{k+2,0}} \lesssim {E}_{m-2} \Big( \varepsilon \|S_n\|_{X^{k+2,0}} + \varepsilon \|\partial_z S_n\|_{X^{k+1,0}} \Big).
$$
Note that right hand side can be bounded by dissipation type energy. For third term,
$$
\|\Pi^b F_{S^\chi} \mathbf{n}^b \|_{Y^{k,0}} \lesssim {E}_{m-2} \Lambda \Big(\frac{1}{c_0}, 1 + \varepsilon\|S_n\|_{Y^{k+2,0}} + \|\Pi^b((D^\varphi)^2 q)\mathbf{n}^b\|_{Y^{k,0}} \Big).
$$
$\varepsilon\|S_n\|_{Y^{k+2,0}}$ was treated as second term, and $\|\Pi^b((D^\varphi)^2 q)\mathbf{n}^b\|_{Y^{k,0}}$ can be estimated since $\Pi^b\partial_z \sim Z_3$.
$$
\|\Pi^b((D^\varphi)^2 q)\mathbf{n}^b\|_{Y^{k,0}} \lesssim \Lambda \Big( \frac{1}{c_0} \Big) \|\nabla q\|_{Y^{k+1,0}} \lesssim {E}_{m-2} \Big( 1 + \varepsilon\|S_n\|_{X^{k+4,0}} \Big).
$$
where last term can be treated similarly as above by anisotropic Sobolev embedding. We performed all the estimates to apply Lemma \ref{high order version} and we get finally
\begin{equation}
\begin{split}
\|S_n(T)\|_{Y^{[\frac{m}{2}],0}}^2 &\leq \|S_n(0)\|_{Y^{[\frac{m}{2}],0}}^2 + \varepsilon\int_0^T {E}_{m-2}(t) \|\nabla S_n(t)\|_{X^{m-2,0}} dt  \\
&\leq \|S_n(0)\|_{Y^{[\frac{m}{2}],0}}^2 +   \int_{0}^{T} {E}^{2}_{m-2}(t) dt +  \varepsilon\int_{0}^{T} \|\nabla S_{n}\|^{2}_{X^{m-2,0}}dt . \\
\end{split}
\end{equation}
\end{proof} 

\section{Vorticity estimate}
From above sections, we checked that regularity for $\partial_z v$ is $L^\infty X^{m-3,0}$. However, we should control $X^{m-1,0}$ in space. Using critical idea of \cite{NMFR}, we claim that we have a regularity of $\partial_z v$ in $L^4_{T} X^{m-1,0}$. We use vorticity to control $\partial_z v$, instead of $S_{n}$. Let us define vorticity:
\begin{equation}
	\omega := \nabla^{\varphi}\times v = (\nabla\times u)\circ\Phi.
\end{equation}
From (\ref{div free normal}), we get
\begin{equation}
\|Z^{\alpha} \partial_z v(t)\| \leq \Lambda_{m,\infty}(t) (\|v(t)\|_{m} + |h(t)|_{m-\frac{1}{2}} + \|\omega(t)\|_{m-1} ),\quad |\alpha|=m-1.
\end{equation}
This means that we are suffice to control vorticity instead of $\partial_z v$. Using vorticity, we can neglect pressure effect, because $\nabla^{\varphi}\times\nabla^{\
	\varphi} q = \nabla\times\nabla p = 0$. However, main difficulty in estimating vorticity is that vorticity $\omega$ does not vanish on the free boundary. To analyze vorcitiy estimate, we divide the equation of vorticity into homogeneous and non-homogeneous parts and apply critical $H^{\frac{1}{4}}$ estimate type estimate from \cite{NMFR}. Applying $\nabla^\varphi\times$ to the Navier-Stokes equation gives, 
\begin{equation} \label{10.2}
\partial_t^\varphi\omega  + (v\cdot\nabla^\varphi)\omega   - (\omega \cdot\nabla^\varphi)v  = \varepsilon\triangle^\varphi\omega.
\end{equation}
Applying $Z^{m-1}$ yields,
\begin{equation} \label{10.3}
\partial_t^\varphi Z^{m-1}\omega + (v\cdot\nabla^\varphi)Z^{m-1}\omega - \varepsilon\triangle^\varphi Z^{m-1}\omega =  F,
\end{equation}
where
\begin{equation} \label{def F}
F := Z^{m-1}(\omega\cdot\nabla^\varphi v) + \mathcal{C}_S, \quad \mathcal{C}_S := \mathcal{C}_S^1 + \mathcal{C}_S^2 + \mathcal{C}_S^3,
\end{equation}
with
\begin{equation} \label{Cs1}
\mathcal{C}_S^1 = [Z^{m-1} v_y]\cdot\nabla_y\omega + [Z^{m-1},V_z]\partial_z\omega := \mathcal{C}^1_{S_y} + \mathcal{C}^1_{S_z},
\end{equation}
\begin{equation} \label{Cs2}
\mathcal{C}_S^2 = \varepsilon[Z^{m-1},\Delta^\varphi]\omega.
\end{equation}
And for boundary data we have, by Lemma \ref{lem 3.2}, 
\begin{equation} \label{10.5}
|(Z^{m-1}\omega)^b (t)| \leq \Lambda_{\infty,m}(t)(|v^b(t)|_{X^{m,0}} + |h(t)|_{X^{m,0}}).
\end{equation}
Using trace theorem,
$$
|(Z^{m-1}\omega)^b (t)| \leq \Lambda_{m,\infty}(t) \Big( \|\nabla v(t)\|_{X^{m,0}}^{\frac{1}{2}} \|v(t)\|_{X^{m,0}}^{\frac{1}{2}} + \|v(t)\|_{X^{m,0}} + |h(t)|_{X^{m,0}} \Big).
$$
Because of the first term on the RHS, the only way to control this boundary terms is to treat this term as dissipation, i.e. 
\begin{equation} \label{10.7}
\sqrt{\varepsilon}\int_{0}^{T} |(Z^{m-1}\omega)^b|^2 dt \leq \sqrt{\varepsilon} \int_{0}^{T} \Lambda_{m,\infty}(t)  \big( \|\nabla v\|_{X^{m,0}}\|v\|_{X^{m,0}} + \|v\|_{X^{m,0}}^2 + |h|_{X^{m,0}}^2 \big) dt.
\end{equation}
Using Young's inequality,
$$
\sqrt{\varepsilon}\int_0^{T} |(Z^{m-1}\omega)^b|^2 \leq \varepsilon\int_0^T \|\nabla v\|_{X^{m,0}}^2 dt + \int_0^T \Lambda_{m,\infty}(t) \big(\|v\|_{X^{m,0}}^2 + |h|_{X^{m,0}}^2 \big) dt. \\
$$

Now, we split vorticity into two parts.
$$
Z^{m-1}\omega := \omega^{m-1}_{h} + \omega^{m-1}_{nh}.
$$
\noindent I) $\omega^{m-1}_{nh}$ solves nonhomogeneous equation:
\begin{equation} \label{10.8}
\partial_t^\varphi \omega^{m-1}_{nh} + (v\cdot\nabla^\varphi)\omega^{m-1}_{nh}  - \varepsilon\triangle^\varphi\omega^{m-1}_{nh} =  F,
\end{equation}
with non-zero initial data and zero-boundary condition,
\begin{equation} \label{10.8-1}
(\omega^{m-1}_{nh})^b = 0,\,\,\,(\omega^{m-1}_{nh})_{t=0} = \omega(0).
\end{equation}
Note that $F$ on the right hand side is defined in (\ref{def F}), (\ref{Cs1}), and (\ref{Cs2}). 	\\
\noindent II) $\omega^{m-1}_{h} $ solves homogeneous equation :
\begin{equation} \label{10.10}
\partial_t^\varphi \omega^{m-1}_{h} + (v\cdot\nabla^\varphi)\omega^{m-1}_{h} - \varepsilon\triangle^\varphi\omega^{m-1}_{h} =  0,
\end{equation}
with zero initial data and nonhomogeneous boundary condition,
\begin{equation} \label{10.10-1}
(\omega^{m-1}_{h})^b = (Z^{m-1}\omega)^{b},\,\,\,(\omega^{m-1}_{h})_{t=0} = 0.
\end{equation}
We state the energy estimates for these two vorticity terms $\omega_{nh}^{m-1}$ and $\omega_{h}^{m-1}$. \\
\subsection{Non-homogeneous estimate}
We estimate $\omega^{m-1}_{nh}$ using (\ref{10.8}) and (\ref{10.8-1}). Similar as Proposition 10.2 in \cite{NMFR}, we get the following Proposition.
\begin{proposition} \label{proposition 10.3}
	Non-homogeneous part $\omega_{nh}^{m-1}$ estimate.
	\begin{equation} \label{10.27}
	\begin{split}
	\|\omega^{m-1}_{nh}(T)\|^2  + \varepsilon\int_0^T \|\nabla^\varphi\omega^{m-1}_{nh}\|^2 dt &\leq \Lambda_0\Big( \|\omega^{m-1}_{nh}(0)\|^2 \Big) + \varepsilon \int_0^T  \Lambda_{m,\infty}(t)  \|\nabla S_n\|_{X^{m-2,0}}^2 dt \\
	&+ \int_0^T \Lambda_{m,\infty}(t) \Big( \|v\|_{X^{m,0}}^2 + \|\omega\|_{m-1}^2 + |h|_{X^{m,1}}^2 \Big) ds.
	\end{split}
	\end{equation}
\end{proposition}
\subsection{Homogeneous estimate}
We estimate $\omega^\alpha_{h}$.  We use critical estimate Proposition 10.4 in \cite{NMFR}. Instead of conormal derivatives $Z^{\alpha}$ of \cite{NMFR}, we can change $Z^{\alpha}$ into space-time conormal derivatives $Z^{m-1}$ version. 
\begin{proposition}  \label{L4 estimate}
	Let us assume that
	$$
	\sup_{[0,T]}\Lambda_{m,\infty}(t) + \int_0^T \Big(\varepsilon\|\nabla v\|_{X^{m,0}}^2  + \varepsilon\|\nabla S_n\|_{X^{m-2,0}}^2  \Big) \leq M.
	$$
	Then there exists $\Lambda(M)$ such that 
	\begin{equation} \label{10.39}
	\|\omega^{m-1}_{h} \|^2_{L^4(0,T;L^2)} \leq \Lambda(M)\int_0^T \big( \|v\|_{X^{m,0}}^2 + \|\nabla v\|_{X^{m-1,0}}^2 + |h|_m^2 \big) dt  + {\varepsilon}\int_0^T   \left\|\nabla v\right\|_{X^{m,0}}^2 dt .
	\end{equation}
\end{proposition}

Finally combining Proposition \ref{proposition 10.3} and \ref{L4 estimate}, we gain the following vorticity estimate.
\begin{proposition}  \label{vorticity est}
	Let us assume that
	$$
	\sup_{[0,T]}\Lambda_{m,\infty}(t) + \int_0^T \Big(\varepsilon\|\nabla v\|_{X^{m,0}}^2  + \varepsilon\|\nabla S_n\|_{X^{m-2,0}}^2  \Big) \leq M.
	$$
	Then there exists $\Lambda(M)$ such that 
	\begin{equation} 
	\begin{split}
	\|\omega^{m-1} \|^2_{L^4(0,T;L^2)} &\leq \Lambda_0\Big( \|\omega^{m-1}(0)\|^2 \Big) +  \Lambda(M)\int_0^T \big( \|v\|_{X^{m,0}}^2 + \|\nabla v\|_{X^{m-1,0}}^2 + |h|_m^2 \big) dt  \\
	&\quad + {\varepsilon}\int_0^T \Lambda_{m,\infty}(t) \big( \|\nabla S_n\|_{X^{m-2,0}}^2 + \|\nabla v \|_{X^{m,0}}^2 \big) dt .  \\
	\end{split}
	\end{equation}
\end{proposition}

\section{Uniform Regularity and Local Existence}
To get local existence of free-boundary Navier-Stokes equations with viscosity $\varepsilon$, we use local existence theory by A.Tani [7], and combine our propositions to get uniform regularity. From Proposition \ref{energy estimate}, \ref{alltime}, \ref{prop Sn}, \ref{L infty estimate}, and \ref{vorticity est}, we get required initial data:
\begin{equation} \label{time initial}
\begin{split}
	&\|v(0)\|_{X^{m,0}} + \varepsilon\|\nabla v(0)\|_{X^{m-1,1}} + |h(0)|_{X^{m,1}} + \|\nabla v(0)\|_{X^{m-1,0}} + \|\p_{z}v(0)\|_{Y^{[\frac{m}{2}],0}}.   \\
\end{split}
\end{equation}
Since above terms include some terms like
\begin{equation} \label{time deriva ini}
	\|(\p_{t}^{k}v)(0)\|_{X^{m-k,0}} \quad\text{or}\quad |(\p_{t}^{k}h)(0)|_{X^{m-k,1}},
\end{equation}
we can control (\ref{time initial}) by initial data $(v_{0}, h_{0})$ with only spatial derivatives. From (\ref{1.12}) and pressure estimate (\ref{q_E est}), (\ref{q_NS est}), (\ref{q_S est}), 
\begin{equation} \label{initial data}
\begin{split}
	\|(\p_{t}^{m} v)(0)\|_{L^{2}(S)} &\leq \Lambda \big( \frac{1}{c_{0}}, \|(\p_{t}^{m-1}v)(0)\|_{H^{1}(S)} + \varepsilon\|(\p_{t}^{m-1} v)(0)\|_{H^{2}(\p S)} + |(\p_{t}^{m-1}h)(0)|_{H^{\frac{3}{2}}(\p S)} \big).  \\
\end{split}
\end{equation}	
Similarly, for initial data of $h$, using (\ref{1.14}),
\begin{equation} \label{initial data h}
\begin{split}
	|(\p_{t}^{m} h)(0)|_{H^{1}(\p S)} &\leq \Lambda \big( \frac{1}{c_{0}}, \|(\p_{t}^{m-1}v)(0)\|_{H^{1}(S)} + |(\p_{t}^{m-1}h)(0)|_{H^{\frac{5}{2}}(\p S)} \big)  \\
\end{split}
\end{equation}
We apply same method for the terms with $\p_{t}^{m-1}$ on the RHS of (\ref{initial data}) and
(\ref{initial data h}) inductively to get
\begin{equation} \label{whole initial data}
\begin{split}
(\ref{time initial}) &\leq \Lambda\big(\frac{1}{c_{0}}, \|v_{0}\|_{H^{m}(S)} + \sum_{k=1}^{m} \varepsilon^{k}\|v_{0}\|_{H^{m+k}(S)} + |h_{0}|_{H^{\frac{3}{2}m + 1}(\p S)} \big)  \\
&:= \Lambda\big(\frac{1}{c_{0}}, I_{m}(0) \big) ,  \\
\end{split}
\end{equation}
where initial data satisfy compatibility conditions, 
\begin{equation} \label{compatibility}
\Pi \big( S^\varphi v^\varepsilon \mathbf{n} \big) \big\vert_{t=0,z=0} = 0,\quad  \Pi := {I} - \mathbf{n}\otimes\mathbf{n},
\end{equation}
for $m \geq 6$ ($\frac{m}{2} + 3 \leq m$). Note that our initial data should be standard Sobolev space not only in conormal space by some terms like (\ref{time deriva ini}). 

First we regularize $v_0^\varepsilon$ by parameter $\delta$, so that $v_0^{\varepsilon,\delta} \in H^{l+1}(S)$, where $l \in (\frac{1}{2},1)$. Then by [7], for this initial condition, there is a time interval $T^{\varepsilon,\delta}$ such that on $[0,T^{\varepsilon,\delta}]$, we have a unique solution $v^{\varepsilon,\delta} \in W_2^{l+2,\frac{l}{2}+1}\Big([0,T^{\varepsilon,\delta}]\times S\Big) = L^2\Big([0,T^{\varepsilon,\delta}],H^{l+2}(S)\Big) \cap H^{\frac{l}{2}+1}\Big([0,T^{\varepsilon,\delta}],L^2(S)\Big)$. For convenience, we use $v$ instead of $v^{\varepsilon,\delta}$. Then by parabolic regularity, for $T\in [0,T^{\varepsilon,\delta}]$,
\begin{equation} \label{energy Nm}
\begin{split}
\mathcal{N}_m (T) &:= \sup_{t\in[0,T]}\Big( \|v(t)\|_{X^{m,0}}^2 + |h(t)|_{X^{m,1}}^2 + \|\partial_z v(t)\|_{X^{m-2,0}}^2 + \|\partial_z v(t)\|_{Y^{[\frac{m}{2}],0}}^2 \Big)  \\
&\quad + \|\partial_z v\|_{L_{T}^4 X^{m-1,0}}^2 + \varepsilon\int_0^T \Big( \|\nabla v (t)\|^2_{X^{m,0}} + \|\nabla\partial_z v(t)\|^2_{X^{m-2,0}} \Big) dt < \infty.
\end{split}
\end{equation}

Let us suppose $\mathcal{N}_m(T_0)<\infty$ for some $T_0$, then using Stoke's regularity on $[\frac{T_0}{2},T_0]$, we know that $v(T_0)\in H^{l+1}(S)$, so by considering this as initial condition again, we know that it can be extended to some $T_1>T_0$. We have to show that this extension is uniform in $\varepsilon$ and $\delta$ using our propositions. \\

We also define equivalent quantity 
\begin{equation} \label{energy Em}
	\mathcal{E}_{m}(T) := \sup_{t\in[0,T]}\mathcal{Q}_{m}(t) + \mathcal{D}_{m}(T) + \|\omega\|_{L^{4}_{T}X^{m-1,0}},
\end{equation}
where 
\begin{equation}
\begin{split}
	\mathcal{Q}_{m}(t) &:= \|v(t)\|_{X^{m,0}}^2 + |h(t)|_{X^{m,1}}^2 + \|S_{n}(t)\|_{X^{m-2,0}}^2 + \|S_{n}(t)\|_{Y^{[\frac{m}{2}],0}}^2,  \\
	\mathcal{D}_{m}(T) &:= \varepsilon\int_0^T \Big( \|\nabla v (t)\|^2_{X^{m,0}} + \|\nabla S_{n} (t)\|^2_{X^{m-2,0}} \Big) dt.
\end{split}
\end{equation}
We know that $\mathcal{N}_{m}(T)$ and $\mathcal{E}_{m}(T)$ are equivalent,
\[
	\mathcal{N}_{m}(T) \leq \Lambda\big( \frac{1}{c_{0}}, \mathcal{E}_{m}(T) \big)\quad \text{and}\quad \mathcal{E}_{m}(T) \leq \Lambda\big( \frac{1}{c_{0}}, \mathcal{N}_{m}(T) \big).
\]

To derive uniform time interval, we choose $R$ and $c_0$ so that, $\frac{1}{c_0} << R$, and define,
\begin{equation} \label{conditions}
T_*^{\varepsilon,\delta} = \sup\left\{ T\in[0,1]\,\,\vert\,\, \mathcal{E}_m(t)\leq R,\,\,\,|h(t)|_{2,\infty} \leq \frac{1}{c_0},\,\,\,\partial_z\varphi(t) \geq c_0,\,\,\,\forall t\in [0,T] \right\}.
\end{equation}

Now, let us combine Propositions to close energy estimate. For convenience we define
\begin{equation} \label{def conv}
\begin{split}
	\mathcal{M}_m (T) &:= \sup_{t\in[0,T]}\Big( \|v(t)\|_{X^{m-1,1}}^2 + |h(t)|_{X^{m-1,2}}^2 + \|\partial_z v(t)\|_{X^{m-2,0}}^2 + \|\partial_z v(t)\|_{Y^{[\frac{m}{2}],0}}^2 \Big)  \\
&\quad + \|\partial_z v\|_{L_{T}^4 X^{m-1,0}}^2 + \varepsilon\int_0^T \Big( \|\nabla v (t)\|^2_{X^{m,0}} + \|\nabla\partial_z v(t)\|^2_{X^{m-2,0}} \Big) dt < \infty.
\end{split}
\end{equation}
For Proposition \ref{energy estimate}, we apply Proposition \ref{dirineumann} to $|h|_{X^{m-1,\frac{5}{2}}}$ for sufficiently small $\theta \ll 1$ and $T \ll 1$ on the right hand side of (\ref{DN est}) so that $\Lambda_{m,\infty}(T) \Big\{ {\theta} |Z^{m-1} \p_{t} h(T)|^2_{L^2} + \sqrt{T} |Z^{m-1}\nabla \p_{t} h(T)|^2_{L^2} \Big\}$ is absorbed by energy of left hand sides except $Z^{m-1}=\p_{t}^{m-1}$ case.
\begin{equation} \label{not all time}
\begin{split}
	& \|v(T)\|^{2}_{X^{m-1,1}} + |h(T)|^{2}_{X^{m-1,2}} + \varepsilon\int_{0}^{T} \|\nabla v (t)\|^2_{X^{m,0}} dt \\
	&\quad \lesssim \Lambda(\frac{1}{c_{0}}, I_{m}(0)) + (\theta + T)|\p_{t}^{m} h|^{2}_{H^{1}} +  \int_{0}^{T} \Lambda(\frac{1}{c_{0}}, \mathcal{Q}_{m}(t)) dt	
\end{split}
\end{equation}
For $Z^{m}=\p_{t}^{m}$ case, we multiply suffciently small constant $\gamma \ll 1$ to Proposition \ref{alltime} and apply Propositino \ref{dirineumann} to $|\p_{t}^{m-1}h|_{H^{\frac{5}{2}}}$ to get
\begin{equation} \label{gamma all time}
\begin{split}
& \gamma\|\p_{t}^{m}v(T)\|^{2} + \gamma|\p_{t}^{m}h(T)|_{H^{1}}^{2} + \gamma\varepsilon\int_{0}^{T} \|\nabla \p_{t}^{m}v (t)\|^2  dt \\
&\quad \lesssim \Lambda(\frac{1}{c_{0}}, I_{m}(0)) + (\theta + \sqrt{T}) |\p_{t}^{m}h(T)|^{2}_{H^{1}}  \\
&\quad\quad  + \gamma\Lambda_{0}(R)\Big( \|v(T)\|^{2}_{X^{m-1,1}} + |\p_{t}^{m-1}h|^{2}_{H^{2}} \Big)  +  \gamma \int_{0}^{T} \Lambda(\frac{1}{c_{0}}, \mathcal{Q}_{m}(t)) dt,	
\end{split}
\end{equation}
where we used notation $\Lambda_{0}(R):=\Lambda(\frac{1}{c_{0}},R)$. We choose $\theta, T, \gamma \ll 1$ and $\theta + \sqrt{T} \ll \gamma$ so that $(\theta + R)|\p_{t}^{m}h|^{2}_{H^{1}}$ on the RHS of (\ref{not all time}) is absorbed by $\gamma|\p_{t}^{m}h(T)|^{2}_{H^{1}}$ in the LHS of (\ref{gamma all time}), and $(\theta + \sqrt{T}) |\p_{t}^{m}h(T)|^{2}_{H^{1}} + \gamma\Lambda_{0}(R)\Big( \|v(T)\|^{2}_{X^{m-1,1}} + |\p_{t}^{m-1}h|^{2}_{H^{2}} \Big)$ in the RHS of (\ref{gamma all time}) is absorbed by $\|v(T)\|^{2}_{X^{m-1,1}} + |h(T)|^{2}_{X^{m-1,2}}$ in the LHS of (\ref{not all time}).  \\

 Now, we combine (\ref{not all time}), (\ref{gamma all time}), and (\ref{vorticity est}) with sufficiently small $\gamma, T \ll 1$. Therefore 
\begin{equation} \label{11111}
\begin{split}
& \|v(T)\|^{2}_{X^{m,0}} + |h(T)|^{2}_{X^{m,1}} + \varepsilon\int_{0}^{T} \|\nabla v (t)\|^2_{X^{m-1,1}} dt  \\
&\quad \lesssim \Lambda(\frac{1}{c_{0}}, I_{m}(0)) + \int_{0}^{T} \Lambda(\frac{1}{c_{0}}, \mathcal{Q}_{m}(t)) dt.	
\end{split}
\end{equation}

Applying Proposition \ref{dirineumann} to Proposition \ref{prop Sn}, 
\begin{equation} \label{22222}
\begin{split}
&\left\|S_n(T)\right\|^2_{X^{m-2,0}} + \varepsilon\int_0^T \left\|\nabla S_n(t)\right\|^2_{X^{m-2,0}} dt \leq \Lambda(\frac{1}{c_{0}}, I_{m}(0)) \\
&\quad + \Lambda_{0}(R) (\theta + \sqrt{T}) |\p_{t}^{m}h|^{2}_{H^{1}} + \int_{0}^{T} \Lambda(\frac{1}{c_{0}}, \mathcal{Q}_{m}(t)) dt,\quad \theta, T \ll 1.	
\end{split}
\end{equation}
From (\ref{vorticity est}),
\begin{equation} \label{33333}
\begin{split}
\|\omega \|^2_{L^4_{T}X^{m-1,0}} &\leq \Lambda(\frac{1}{c_{0}}, I_{m}(0)) +  \Lambda(\frac{1}{c_{0}}, \mathcal{Q}_{m}(t)) dt \\
&\quad + {\varepsilon}\int_0^T \Lambda_{m,\infty}(t) \big( \|\nabla S_n\|_{X^{m-2,0}}^2 + \|\nabla v \|_{X^{m,0}}^2 \big) dt .  \\
\end{split}
\end{equation}
We combine (\ref{11111}), (\ref{22222}), and (\ref{33333}) with sufficiently small $\theta, T \ll 1$ then finally we get
\begin{eqnarray} \label{condition 1}
\mathcal{E}_m (T) \leq \Lambda\Big(\frac{1}{c_0},I_m(0)\Big)  + \Lambda_0(R)T,\quad \Lambda_{0}(R) := \Lambda(\frac{1}{c_{0}},R).
\end{eqnarray}
Now, we calculate conditions in (\ref{conditions}).
\begin{equation} \label{condition 2}
\left|h(t)\right|_{2,\infty} \leq \left|h(0)\right|_{2,\infty} + \Lambda_0(R)T\,\,\,\,\,\,\forall t\in [0,T],
\end{equation}
and since we've chosen $A$ in diffeomorphism to be 1, at initial time,
\begin{equation} \label{condition 3}
\partial_z\varphi(t) \geq 1 - \int_0^t \left\|\partial_t\nabla\eta\right\|_{L^\infty} \geq 1 - \Lambda_0(R)T,\,\,\,\,\,\,\forall t\in [0,T]. 
\end{equation}

\noindent From (\ref{condition 1}),(\ref{condition 2}), and (\ref{condition 3}), we see that right hand side is independent to $\varepsilon$ and $\delta$, so are possible to choose $R = \Lambda\Big(|h_0|_{2,\infty},I_m(0)\Big)$ which satisfies that there exist $T_*$(independent to $\varepsilon,\delta$) such that $\forall t\in [0,T_*]$,
\begin{equation}
{\mathcal{N}}_m(t) \leq \frac{R}{2},\,\,\,\,\,\left|h(t)\right|_{2,\infty} \leq \frac{1}{2c_0},\,\,\,\,\,\partial_z\varphi(t) \geq c_0 + \frac{1-c_0}{2} > c_0.
\end{equation}
This implies $T_* < T_*^{\varepsilon,\delta}$, hence $T_*$ implies there exist uniform time, independent to $\varepsilon,\delta$. Since $ {\Psi}_m(T_*)$ is uniformly bounded in $\delta$, we can pass the limit, $\delta \rightarrow 0$, by strong compactness argument. \\


\section{Uniqueness}
\subsection{Uniqueness for Navier-Stokes.}
We prove uniqueness of Theorem \ref{inviscid limit}. As usual, we consider two solution sets $(v_1^\varepsilon, \varphi_1^\varepsilon, q_1^\varepsilon), (v_2^\varepsilon, \varphi_2^\varepsilon, q_2^\varepsilon)$ with same initial condition and proper compatibility conditions. Then on the interval $[0,T^\varepsilon]$, we have uniform bounds of energy,
$$
\Psi_m^i (T^\varepsilon) \leq R,\,\,\,\,i=1,2.
$$ 
Let us define,
\begin{equation}
\bar{v}^\varepsilon := v_1^\varepsilon - v_2^\varepsilon,\,\,\,\bar{h}^\varepsilon := h_1^\varepsilon - h_2^\varepsilon,\,\,\,\bar{q}^\varepsilon := q_1^\varepsilon - q_2^\varepsilon.
\end{equation}
We will make system of equations for $(\bar{v}^\varepsilon, \bar{h}^\varepsilon, \bar{q}^\varepsilon)$ and do energy estimate. By divergence free condition, $\nabla^{\varphi_i}\cdot v_i^\varepsilon = 0$,
$$
\Big(\partial_t + v_{y,i}^\varepsilon\cdot\nabla_y + V_{z,i}^\varepsilon\partial_z\Big)v_{i}^\varepsilon + \nabla^{\varphi_i}q_i^\varepsilon - \varepsilon\Delta^{\varphi_i}v_i^\varepsilon = 0.
$$
Then we have equation about $(\bar{v}^\varepsilon, \bar{h}^\varepsilon, \bar{q}^\varepsilon)$. First for Navier-Stokes,
\begin{equation}
\Big(\partial_t + v_{y,1}^\varepsilon\cdot\nabla_y + V_{z,1}^\varepsilon\partial_z\Big)\bar{v}^\varepsilon + \nabla^{\varphi_1}\bar{q}^\varepsilon - \varepsilon\Delta^{\varphi_1}\bar{v}^\varepsilon = F,
\end{equation}
where
$$
F = (v_{y,2}^\varepsilon - v_{y,1}^\varepsilon)\cdot \nabla_y v_2^\varepsilon + (V_{z,2}^\varepsilon - V_{z,1}^\varepsilon)\partial_z v_2^\varepsilon - \Big(\frac{1}{\partial_z \varphi_2^\varepsilon} - \frac{1}{\partial_z \varphi_1^\varepsilon}\Big)\Big(P_1^*\nabla q_2^\varepsilon\Big) + \frac{1}{\partial_z \varphi_2^\varepsilon}\Big((P_2 - P_1)^*\nabla q_2^\varepsilon\Big) 
$$
$$
+ \varepsilon\Big(\frac{1}{\partial_z \varphi_2^\varepsilon} - \frac{1}{\partial_z\varphi_1^\varepsilon}\Big)\nabla\cdot(E_1\nabla v_2^\varepsilon) + \varepsilon\frac{1}{\partial_z\varphi_2^\varepsilon}\nabla\cdot\Big((E_2-E_1)\nabla v_2^\varepsilon\Big).
$$
For divergence-free condition,
\begin{equation}
\nabla^{\varphi_1}\cdot \bar{v}^\varepsilon = - \Big(\frac{1}{\partial_z \varphi_2^\varepsilon} - \frac{1}{\partial_z\varphi_1^\varepsilon}\Big)\nabla\cdot\Big(P_1 v_2^\varepsilon\Big) - \frac{1}{\partial_z \varphi_2^\varepsilon} \nabla \cdot \Big((P_2 - P_1) v_2^\varepsilon\Big).
\end{equation}
For Kinematic boundary condition,
\begin{equation}
\partial_t \bar{h}^\varepsilon - (v^\varepsilon)^b_{y,1}\cdot\nabla h + \Big((v_{z,1}^\varepsilon)^b - (v_{z,2}^\varepsilon)^b \Big) = -\Big((v_{y,2}^\varepsilon)^b - (v_{y,1}^\varepsilon)^b \Big)\cdot \nabla h_2^\varepsilon.
\end{equation}
Continuity of stress tensor condition becomes,
\begin{equation}
\bar{q}^\varepsilon\mathbf{n}_1 - 2\varepsilon\Big(S^{\varepsilon_1}\bar{v}^\varepsilon\Big)\mathbf{n}_1 = g\bar{h}^\varepsilon - \eta\nabla\cdot\Big(\frac{\nabla\bar{h}^\varepsilon}{\sqrt{1+\left|\nabla h_1^\varepsilon\right|^2}}\Big)
\end{equation}
$$
+ 2\varepsilon\Big( \Big( S^{\varphi_1} - S^{\varphi_2} \Big) v_2^\varepsilon\Big)\mathbf{n}_1 + 2\varepsilon\Big( S^{\varphi_2} v_2^\varepsilon\Big)\Big( \mathbf{n}_1 - \mathbf{n}_2 \Big) - \eta\nabla\cdot\Big( \left\{ \frac{1}{\sqrt{1+\left|\nabla h_1^\varepsilon\right|^2}} - \frac{1}{\sqrt{1+\left|\nabla h_2^\varepsilon\right|^2}} \right\} \nabla h_2^\varepsilon \Big).
$$
Using above 4 equations, we get $L^2$- energy estimate,(since initial condition is zero, no initial term appear on right hand side)
\begin{equation}
\left\|\bar{v}^\varepsilon(t)\right\|_{L^2}^2 + \left|\bar{h}^\varepsilon(t)\right|_{H^1}^2 + \varepsilon\int_0^t\left\|\nabla\bar{v}^\varepsilon\right\|_{L^2}^2 ds \leq \Lambda(R)\int_0^t \Big( \left\|\bar{v}^\varepsilon(s)\right\|_{L^2}^2 + \left|\bar{h}^\varepsilon(s)\right|_{H^{\frac{3}{2}}}^2 \Big)ds.
\end{equation}
We skip detail calculation, since we can use our previous energy estimates basically. But, in above equations for $(\bar{v}^\varepsilon,\bar{h}^\varepsilon, \bar{q}^\varepsilon)$, right hand side does not have low order than $L^2$ energy. However we have uniform bound of m-order energy, so we can extract bad high order terms into $\Lambda(R)$. To estimate $\left|\bar{h}^\varepsilon(s)\right|_{H^{\frac{3}{2}}}$, we don't have to take time derivatives as like in previous section, since we already have bounded high order energy. And moreover, we don't need uniform estimate in $\varepsilon$, since we are dealing about for fixed $\varepsilon$. So, we estimate $\varepsilon\left|\bar{h}^\varepsilon(s)\right|_{H^{\frac{3}{2}}}^2$.
\begin{lemma}
For every $m\in\mathbb{N}$, $\varepsilon\in (0,1)$, we have the estimate
\begin{equation}
\varepsilon\left|h(t)\right|_{m+\frac{1}{2}}^2 \leq \varepsilon\left|h(0)\right|_{m+\frac{1}{2}}^2 + \varepsilon\int_0^t |v^b|_{m+\frac{1}{2}}^2 + \int_0^t \Lambda_{1,\infty}\Big( \left\|v\right\|_m^2 + \varepsilon|h|_{m+\frac{1}{2}}^2 \Big) ds,
\end{equation}
where
$$
\Lambda_{1,\infty} = \Lambda\Big(|\nabla h|_{L^\infty} + \|v\|_{1,\infty} \Big).
$$
\end{lemma}
\begin{proof}
See proposition 3.4 in [1].
\end{proof}
This is true for our case, since it comes from Kinematic boundary condition. We can also apply this lemma, to $\bar{h}^\varepsilon$ case,(surely, $\bar{h}^\varepsilon(0)=0$) and then combine with above $L^2$ estimate, then we get the following estimate.(non-uniform in $\varepsilon$)
\begin{equation}
\left\|\bar{v}^\varepsilon(t)\right\|_{L^2}^2 + \left|\bar{h}^\varepsilon(t)\right|_{H^1}^2 + \varepsilon\left|\bar{h}^\varepsilon(t)\right|_{H^{\frac{3}{2}}}^2 + \varepsilon\int_0^t\left\|\nabla\bar{v}^\varepsilon\right\|_{L^2}^2 ds \leq \frac{\Lambda(R)}{\varepsilon} \int_0^t \Big( \left\|\bar{v}^\varepsilon(s)\right\|_{L^2}^2 + \varepsilon\left|\bar{h}^\varepsilon(s)\right|_{H^{\frac{3}{2}}}^2 \Big)ds.
\end{equation}
Then we can use Gronwall's inequality to get uniqueness. So finish uniqueness part of theorem \ref{theorem-uniform}.

\subsection{Uniqueness for Euler.}
Since our estimate in above subsection(uniqueness for Navier-Stokes) is not uniform in $\varepsilon$, result cannot be applied to Euler equation. As like in Navier-Stokes case, let we have two solutions $(v_1, h_1, q_1),(v_2, h_2, q_2)$ with same initial condition. Suppose,
\begin{equation}
\sup_{[0,T]}\Big( \|v_i\|_m + \left\|\partial_z v_i\right\|_{m-1} + \left\|\partial_z v_i\right\|_{[\frac{m}{2}],\infty} + |h_i|_{m+1} \Big) \leq R,\,\,\,\,i=1,2.
\end{equation}
(This is true from result in Theorem \ref{theorem-uniform}) Define $\bar{v} \doteq v_1 - v_2,\,\,\,\bar{h} \doteq h_1 - h_2,\,\,\,\bar{q} \doteq q_1 - q_2$ and we write equation of $(\bar{v}, \bar{h}, \bar{q})$, as before. Euler equation becomes,
\begin{equation}
\Big(\partial_t + v_{y,1}\cdot\nabla_y + V_{z,1}\partial_z\Big)\bar{v} + \nabla^{\varphi_1}\bar{q} = F',
\end{equation}
where
$$
F' = (v_{y,2} - v_{y,1})\cdot \nabla_y v_2 + (V_{z,2} - V_{z,1})\partial_z v_2 - \Big(\frac{1}{\partial_z \varphi_2} - \frac{1}{\partial_z \varphi_1}\Big)\Big(P_1^*\nabla q_2\Big) + \frac{1}{\partial_z \varphi_2}\Big((P_2 - P_1)^*\nabla q_2\Big).
$$
For divergence-free condition,
\begin{equation}
\nabla^{\varphi_1}\cdot \bar{v} = - \Big(\frac{1}{\partial_z \varphi_2} - \frac{1}{\partial_z\varphi_1}\Big)\nabla\cdot\Big(P_1 v_2\Big) - \frac{1}{\partial_z \varphi_2} \nabla \cdot \Big((P_2 - P_1) v_2\Big).
\end{equation}
For Kinematic boundary condition,
\begin{equation}
\partial_t \bar{h} - v^b_{y,1}\cdot\nabla h + \Big( v_{z,1}^b - v_{z,2}^b \Big) = -\Big( v_{y,2}^b - v_{y,1}^b \Big)\cdot \nabla h_2.
\end{equation}
Continuity of stress tensor condition becomes,
\begin{equation}
\bar{q}\mathbf{n}_1 = g\bar{h} - \eta\nabla\cdot\Big(\frac{\nabla\bar{h}}{\sqrt{1+\left|\nabla h_1\right|^2}}\Big) - \eta\nabla\cdot\Big( \left\{ \frac{1}{\sqrt{1+\left|\nabla h_1\right|^2}} - \frac{1}{\sqrt{1+\left|\nabla h_2\right|^2}} \right\} \nabla h_2 \Big).
\end{equation}
By performing basic $L^2$-estimate, as similarly above, (we skip detail here)
\begin{equation}
\left\|\bar{v}(t)\right\|_{L^2}^2 + \left|\bar{h}(t)\right|_{H^1}^2 \leq \Lambda(R)\int_0^t \Big( \left\|\bar{v}(s)\right\|_{H^1}^2 + \left|\bar{h}(s)\right|_{H^{\frac{3}{2}}}^2 \Big)ds.
\end{equation}
We should control $\|v\|_1$ on right hand side. But, since there are no dissipation on left hand side, we cannot make it absorbed. Instead, we use vorticity. Let's define vorticity $\omega = \nabla^{\varphi} \times v$ (which is equivalent to $\omega = (\nabla\times u)(t,\Phi)$ ). We have
$$
\omega\times\mathbf{n} = \frac{1}{2}\Big( D^\varphi v\mathbf{n} - (D^\varphi v)^T \mathbf{n} \Big) 
$$
$$
= S^\varphi v\mathbf{n} - (D^\varphi v)^T \mathbf{n} = \frac{1}{2}\partial_{\mathbf{n}}u - g^{ij}\Big(\partial_j v\cdot \mathbf{n}\Big)\partial_{y^i}.
$$
Hence, it is suffice to estimate $\omega$ instead of $\partial_z v$, $i.e$
\begin{equation}
\left\|\partial_z v\right\|_{L^2} \leq \Lambda(R)\Big( \|\omega\|_{L^2} + \|v\|_1 + |h|_{\frac{3}{2}} \Big).
\end{equation}
To estimate $\omega$, we use vorticity equation. 
\begin{equation}
\Big( \partial_t^{\varphi_i} + v_i\cdot\nabla^{\varphi_i} \Big) \omega_i = \Big( \omega_i\cdot\nabla^{\varphi_i} \Big) v_i.
\end{equation}
$L^2$ energy estimate of $\bar{\omega}$ is 
\begin{equation}
\left\|\bar{\omega}(t)\right\|_{L^2}^2 \leq \Lambda(R)\int_0^t \Big( |\bar{h}(s)|_1^2 + \|\bar{v}(s)\|_1^2 + \|\partial_z \bar{v}(s)\|_{L^2}^2 + \|\bar{\omega}(s)\|_{L^2}^2 \Big) ds.
\end{equation}
We also should control $|h|^2_{L^2 H^{\frac{3}{2}}}$. As similar to Dirichlet-Neumann estimate we can control this by $|\partial_t h|_{L^2 L^2}^2$. And, from kinematic boundary condition of $\bar{h}$, we easily get 
\begin{equation}
\left|\partial_t\bar{h}(t)\right|_{L^2}^2 \leq \Lambda(R) \Big( \left\|\nabla\bar{v}\right\|_{L^2}^2 + |\bar{h}|_{H^1}^2 \Big).
\end{equation}
Then we can use Gronwall's inequality to get uniqueness. So finish uniqueness part of theorem \ref{theorem-uniform} and theorem \ref{inviscid limit}.

\section{Inviscid limit}
In this section we send $\varepsilon$ to zero, and get a unique solution of free boundary Euler equation. For $\varepsilon \in (0,1]$ and $T\leq T_*$, we have uniform energy bound,
\begin{equation}
\mathcal{N}_{m} (T) \leq C.
\end{equation}
So we have uniform boundness for $v^\varepsilon$ in $L^\infty\Big([0,T],H_{co}^{m}\Big)$ and $h^\varepsilon$ in $L^\infty\Big([0,T],H^{m+1}\Big)$. And, by Rellich-Kondrachov theorem, we have, for each $t$, compactness of $v^\varepsilon(t)$ in $H_{co,loc}^{m-1}$ and $h^\varepsilon(t)$ in $H_{loc}^{m}$. And from our energy function, we have a uniform boundness of $\partial_t v^\varepsilon(t)$ in $H_{co,loc}^{m-1}$ and of $\partial_t h^\varepsilon(t)$ in $H_{loc}^{m}$ for $\forall t \in [0,T]$. Now, we have subsequence $v^{\varepsilon_n},\,h^{\varepsilon_n}$, such that 
\begin{equation}
v^{\varepsilon_n} \rightarrow v,\,\,strongly\,\,in\,\,C([0,T],H_{co,loc}^{m-1}),
\end{equation}
$$
h^{\varepsilon_n} \rightarrow h,\,\,strongly\,\,in\,\,C([0,T],H_{loc}^{m}).
$$
About pressure, from pressure section, we have boundness of $\nabla q^\varepsilon$ in $L^2([0,T]\times S)$, so get some $q$ such that,
$$
\nabla q^{\varepsilon_n} \rightharpoonup \nabla q,\,\,weakly\,\,in\,\,L^2([0,T]\times S)
$$
and limit functions $(v,h,\nabla q)$ satisfy
\begin{equation}
\sup_{t\in[0,T]}\Big( \|v(t)\|_{H_{co}^{m}}^2 + |h(t)|_{H^{m+1}}^2 + \|\p_z v(t)\|_{H_{co}^{m-2}}^2 + \|\p_z v\|_{W_{co}^{[\frac{m}{2}],\infty}}^2 \Big) + \|\p_z v\|_{L^4_{T} H^{m-1}_{co}}^2 \leq R,
\end{equation}
Now, we can pass to the limit and get the fact that $(v,h,\nabla q)$ is a weak solution of Euler equation.(interior). For boundary condition, first we can assume that the trace(boundary function),
\begin{equation}
v^{\varepsilon_n}(t,y,z=0) \rightharpoonup v^b\,\,\text{weakly in}\,\,L^2([0,T]\times S)
\end{equation}
for some $v^b$. In kinematic boundary condition, $v^b$ is linear and we have strong convergence of $h$, so kinematic boundary condition is satisfied weakly surely. Next, for continuity of stress tensor condition, by bounded lipschitz norm of $v^{\varepsilon_n}$, $2\varepsilon(Su)\mathbf{n} \rightarrow 0$ in weak limit process. And, limit of surface tension part is trivial by stong convergence of $h$. Hence, in the weak sense, 
\begin{equation}
q^b = gh - \eta\nabla\cdot\Big(\frac{\nabla h}{\sqrt{1+|\nabla h|^2}} \Big).
\end{equation}
Hence we have a weak solution $(v,h)$ which is strong $H_{co,loc}^{m-1} \times H_{co,loc}^{m}$ - convergence of $(v^{\varepsilon_n}, h^{\varepsilon_n})$.(and global for weak convergence in $L^2 \times H^1$) Moreover this limit is unique by previous section. Meanwhile, we can get strong convergence (non-local) in $L^2 \times H^1$. To get this, we just investigate norm convergence. For $t\in [0,T]$, using basic $L^2$-energy estimate and uniform boundness of high order energy, 
\begin{equation}
\Big( \left\|v^{\varepsilon_n}J^{\varepsilon_n}(t)\right\|_{L^2}^2 + g\left|h^{\varepsilon_n}(t)\right|_{L^2}^2 + 2\eta |\sqrt{1+|\nabla h^{\varepsilon_n}(t)|^2} - 1 |_{L^1} \Big) 
\end{equation}
$$
- \Big( \left\|v_0^{\varepsilon_n}J_0^{\varepsilon_n}\right\|_{L^2}^2 + g\left|h_0^{\varepsilon_n}\right|_{L^2}^2 + 2\eta | \sqrt{1+|\nabla h_0^{\varepsilon_n}|^2} - 1 |_{L^1} \Big) \leq \varepsilon_n\Lambda(R) \rightarrow 0,\,\,as\,\,\varepsilon_n\rightarrow 0,
$$
where $J^{\varepsilon} \doteq (\partial_z\varphi^\varepsilon)^{1/2}$ and $\varepsilon_n$ on the right hand side come from dissipation of energy estimate. We assume that $\|v_0^\varepsilon - v_0\|_{L^2}\rightarrow 0,\,\,\|h_0^\varepsilon - h_0\|_{H^1}\rightarrow 0$ in statement of theorem 1.3 and 
$$
\left\|\partial_z\varphi^\varepsilon\vert_{t=0} - \partial_z\varphi\vert_{t=0} \right\|_{L^2} \lesssim \left|h_0^\varepsilon - h_0\right|_{H^{\frac{1}{2}}} \leq \Lambda(R)\left|h_0^\varepsilon - h_0\right|_{L^2}^{\frac{1}{2}}.
$$
This implies 
\begin{equation}
\begin{split}
&\lim_{\varepsilon\rightarrow 0} \Big( \left\|v_0^{\varepsilon_n}J_0^{\varepsilon_n}\right\|_{L^2}^2 + g\left|h_0^{\varepsilon_n}\right|_{L^2}^2 + 2\eta | \sqrt{1+|\nabla h_0^{\varepsilon_n}|^2} - 1 |_{L^1} \Big) \\
&= \left\|v_0 J_0\right\|_{L^2}^2 + g\left|h_0\right|_{L^2}^2 + 2\eta | \sqrt{1+|\nabla h_0|^2} - 1 |_{L^1}.
\end{split}
\end{equation}
Lastly, using energy conservation in Euler equation ($\varepsilon=0$, in basic $L^2$-estimate), we get
\begin{equation}
\left\|v_0 J_0\right\|_{L^2}^2 + g\left|h_0\right|_{L^2}^2 + 2\eta | \sqrt{1+|\nabla h_0 |^2} - 1 |_{L^1} = \left\|v J\right\|_{L^2}^2 + g\left|h\right|_{L^2}^2 + 2\eta | \sqrt{1+|\nabla h|^2} - 1 |_{L^1}.
\end{equation}
Finally we get norm convergence.
\begin{equation}
\begin{split}
&\lim_{\varepsilon\rightarrow 0} \Big( \left\|v^{\varepsilon_n}J^{\varepsilon_n}(t)\right\|_{L^2}^2 + g\left|h^{\varepsilon_n}(t)\right|_{L^2}^2 + 2\eta | \sqrt{1+|\nabla h^{\varepsilon_n}(t)|^2} - 1 |_{L^1} \Big) \\
&= \left\|v J\right\|_{L^2}^2 + g\left|h\right|_{L^2}^2 + 2\eta | \sqrt{1+|\nabla h|^2} - 1 |_{L^1}.
\end{split}
\end{equation}
With weak convergence, this implies strong convergence to $(vJ,h)$ in $L^2 \times H^1$. And strong convergence of $h$ means $(v^\varepsilon, h^\varepsilon)\rightarrow (v,h)$  strongly in $L^2 \times H^1$. (without $J$) $L^\infty$-type convergence can be done by $L^2$ convergence, uniform energy boundness, and anisotropic embedding.

\section*{Acknowledgments}
The authors would like to thank Professor N. Masmoudi for helpful discussions. The authors  were  partially supported by NSF grant DMS-1211806. Donghyun Lee also thanks to Yanjin Wang for fruitful discussion and suggestion. This research was performed when the authors were in Courant Institute.

\vspace{1cm}
\indent\indent \author{$\textsc{Donghyun Lee}$}\\
\indent \address{$\textsc{University of Wisconsin-Madison, Van Vleck Hall, 480 Lincoln Drive,}$}  \\
\indent \address{$\textsc{Madison, Wi, 53706}$} \\
\indent \email{$\mathtt{donghyun@math.wisc.edu}$}

\vspace{0.5cm}
\author{$\textsc{Tarek\,\,\,M.\,\,\,Elgindi}$}\\
\indent \address{$\textsc{Princeton\,\,\,University,\,\,\,Fine Hall,\,\,\,Washington Road,\,\,\,Princeton,\,\,\,NJ,\,\,\,08544}$} \\
\indent \email{$\mathtt{tme2@math.princeton.edu}$}

\end{document}